\documentstyle[12pt, epsf, amscd, amssymb]{amsart}
\newcommand\bR{{\bf R}}
\newcommand\bC{{\bf C}}

\newcommand\rp{\bR P}
\newcommand\GL{{\rm GL}}

\newcommand\PSO{{\rm PSO}}

\newcommand\SLt{{\rm SL}}

\newcommand\SO{{\rm SO}}
\newcommand\dev{{\bf dev}}
\newcommand\SI{{\bf S}}

\newcommand\Bd{{\rm bd}}
\newcommand\clo{{\rm Cl}}
\newcommand\bdd{{\bf d}}
\newcommand\ideal[1]{\tilde #1_\infty}
\newcommand\ra{\rightarrow}

\newcommand\che{\check}
\newcommand\emp{\emptyset}
\newcommand\eps{\epsilon}
\newcommand\disk{{\cal D}}

\newcommand\vth{\vartheta}
\newcommand\vpi{\varphi}
\newcommand\Aff{{\rm Aff}}
\newcommand\ovl{\overline}

\newcommand\Idd{{\rm I}}

\newcommand\cO{{\rm ,}\ }

\newtheorem{thm}{Theorem}
\newtheorem{lem}{Lemma}
\newtheorem{cor}{Corollary}
\newtheorem{prop}{Proposition}
\newtheorem{claim}{Claim}

\theoremstyle{definition}
\newtheorem{defn}{Definition}

\theoremstyle{remark}
\newtheorem{rem}{Remark}      
\newtheorem{hypo}{Hypothesis}

\newcommand\ideals[1]{{#1}_\infty}
\newcommand\uco[1]{{({#1})\tilde{}}}
\newcommand\pch[1]{{({#1})\check{}}}
\newcommand\lideal[1]{({#1})\tilde{}_\infty}

\begin{document}
\title[The universal cover of an affine three-manifold]
{The universal cover of an affine three-manifold 
with holonomy \\ of shrinkable dimension $\leq$ two.}
\author{Suhyoung Choi}
\address{Department of Mathematics \\
College of Natural Sciences \\ 
Seoul National University \\ 
151--742 Seoul, Korea}
\email{shchoi@@math.snu.ac.kr}
\date{June 12, 1999} 
\thanks{Research partially supported by GARC-KOSEF
and SNU Research Fund 1995-03-1036} 

\begin{abstract}
An affine manifold is a manifold with an affine structure, i.e.
a torsion-free flat affine connection. 
We show that the universal cover of a closed affine $3$-manifold $M$ 
with holonomy group of shrinkable dimension (or
{\em discompacit\'e}\/ in French) 
less than or equal to two is diffeomorphic to $\bR^3$. 
Hence, $M$ is irreducible.
This follows from two results: (i) a simply connected  
affine $3$-manifold which is $2$-convex is diffeomorphic 
to $\bR^3$, whose proof using the Morse theory takes most of 
this paper; and (ii) a closed affine manifold of holonomy of 
shrinkable dimension less or equal to $d$ is $d$-convex.
To prove (i), we show that $2$-convexity is a geometric form 
of topological incompressibility of level sets. 
As a consequence, we show that 
the universal cover of a closed affine three-manifold 
with parallel volume form is diffeomorphic to $\bR^3$,
a part of the weak Markus conjecture.  
As applications, we show 
that the universal cover of a hyperbolic $3$-manifold 
with cone-type singularity of arbitrarily assigned 
cone-angles along a link removed with the singular locus 
is diffeomorphic to $\bR^3$.
A fake cell has an affine structure as shown by Gromov. 
Such a cell must have a concave point at the boundary. 
 
\end{abstract}  
\maketitle

Our research is to understand the geometrical
and topological properties of manifolds with flat real 
projective or affine structures, particularly in low-dimensions.  
A hyperbolic structure on a manifold can be naturally regarded 
as a real projective structure using the Klein model of 
hyperbolic space. Our interest stems from the fact that 
the eight $3$-dimensional homogeneous Riemannian structures
can be considered naturally as real projective or product real 
projective structures, as observed by Thurston, Kapovich, 
and many others even earlier (see \cite{psconv} for a proof 
and also Moln\'ar \cite{Mol} and Thiel \cite{Thiel}).
(A product real projective 
structure is a geometric structure modelled on a product of
real projective spaces and a product of projective transformation
groups.) The research on these homogeneous Riemannian 
structures led to the explosion of discoveries on 
$3$-manifold topology as initiated by Thurston in the 1980's,
from which we gained a very practical detailed and comprehensive 
understanding of $3$-manifolds even unparalleled in 
the presently ever expanding study of $4$-manifolds for example.
However, we feel that we may need to look at more general geometric 
structures to gain deeper insights. The study of contact 
structures and foliations form a line of such research.

Classical affine and projective geometries form 
a very rich field giving us much insight into Euclidean, spherical,
and hyperbolic geometries. A comprehensive treatment of 
classical affine and projective geometry by Berger \cite{B} 
for example shows the depth, beauty, and fertility of this subject 
quite clearly. We feel that a study of real projective and affine 
geometric structures has many interesting global properties
and are attempting to develop tools of study (see Benoist \cite{Bens}, 
Goldman \cite{G1}, \cite{G2}, Choi \cite{psconv}, and Choi-Goldman \cite{CG}
for references). 
Most of local properties of affine or projective differential
geometry were discovered much earlier while there seems 
to be only a small number of global results.
More importantly the study of real projective and affine 
structures are important in understanding representations of
discrete subgroups of Lie groups. However, it is beyond 
the subject of this paper.

Topological and geometrical properties of real projective 
and affine structures on $3$-manifolds are completely unknown 
at the moment. For example, we do not have an example of a $3$-manifold 
not admitting a real projective structure. Sullivan and Thurston
\cite{ST} made a conjecture that all $3$-manifolds admit 
some analytic geometric structure. 
This conjecture implies the Poincar\'e conjecture for example. 
(It is shown by Smillie \cite{Sim2} that some $3$-manifolds 
that are homeomorphic to connected sums of Lens spaces do not 
admit flat affine structures.) 

In this paper, we study affine structure on $3$-manifolds.
Affine geometry is the study of properties of the vector space 
$\bR^n$, $n \geq 1$ invariant under the action of
the group $\Aff(\bR^n)$ of affine transformations according 
to Felix Klein's Erlangen program. 
An {\em affine structure} on a smooth $n$-manifold is given by 
a maximal atlas of charts to $\bR^n$ where the transition 
functions are affine. That is, an affine structure gives a way to lift 
affine geometry on a manifold locally and consistently. 
An affine structure can also be given by a torsion-free 
flat affine connection. 
An {\em affine manifold} is a manifold with an affine structure. 
A local diffeomorphism between affine manifolds is {\em affine} 
if it preserves the affine structures locally.

A so-called equi-affine manifold is an affine manifold with 
a volume form parallel under its connection. Equivalently, 
an equi-affine manifold is an affine manifold which admits 
a compatible atlas of charts where the transition functions are 
volume-preserving affine maps.
We will prove that a closed equi-affine $3$-manifold must have a universal 
cover diffeomorphic to a $3$-cell using Morse theory 
adapted to our geometric situation. In particular such 
an affine $3$-manifold is irreducible and does not contain a fake cell. 
This answers a question of Carri\`ere \cite{Car} and 
sharpens Theorem 2 of Smillie \cite{Sim2}. 
(A much shorter version of this paper was written 
as a proceedings paper \cite{uaf1}.)

A quotient of the complete affine space $\bR^n$ by the properly 
discontinuous and free action of a discrete subgroup of 
$\Aff(\bR^n)$ is an example of an affine manifold.
A {\em complete}\/ affine manifold is 
an affine manifold affinely diffeomorphic to such a quotient.
For example, a closed manifold with a flat Riemannian metric is 
affine since it is isometric with 
the quotient of the complete Euclidean space by 
a group of Euclidean isometries (which are affine also). 
A closed flat Lorentz manifold is also affine since  
it is modeled locally on $\bR^n$ and the group of Lorentzian 
affine transformations (which are affine) on 
$\bR^n$ equipped with a Lorentzian inner product.  

Assigning an affine structure to a smooth manifold $M$ is 
equivalent to assigning a pair $(\dev, h)$ where 
$\dev$ is an immersion from $\tilde M$ equivariant with 
respect to a homomorphism $h$ from the group 
$\pi_1(M)$ of deck transformations of $\tilde M$ to $\Aff(\bR^n)$.
$\dev$ is said to be a {\em developing map}, and $h$ a {\em holonomy 
homomorphism}, whose image is said to be the {\em holonomy group}.  
$\dev$ is obtained by analytically continuing local charts 
and $h$ by piecing together the transition functions as we 
continue (see Ratcliff \cite[Chapter 8]{Ratc}).

A complete affine manifold is characterized by the fact 
that the developing map is a diffeomorphism onto $\bR^n$
since its universal cover is affinely diffeomorphic to $\bR^n$.
Numerous incomplete examples 
were constructed by Sullivan and Thurston \cite{ST}. They include 
some Seifert spaces (see Section 2), and bundles over surfaces with 
fibers homeomorphic to tori. 

Affine $2$-manifolds were classified by Nagano and Yagi \cite{NY}.
They are always homeomorphic to tori or Klein bottles 
as was shown by Benz\'ecri-Milnor \cite{Ben} and \cite{Mil3}.  
Complete affine $3$-manifolds are rather limited and 
classified (see Fried-Goldman \cite{FG}).

A form on an affine manifold is {\em parallel}\/ if 
it is parallel with respect to the affine coordinate charts, 
i.e., its covariant derivative with respect to the affine 
connection is zero. An affine manifold admits a parallel 
volume form if and only if $h(\pi_1(M))$ lies in the group of 
equi-affine transformations, i.e., affine transformations with 
Jacobian identically $1$.  

One of the conjectures in the study of 
affine manifolds is that of Markus \cite{Mark}: 
A closed affine manifold is complete if and only it 
admits a parallel volume form. 

The above conjecture is verified in cases 
when the holonomy group $\Gamma$
is abelian, nilpotent, and solvable of rank $\leq n$
by Smillie \cite{Sim}, Fried, Goldman, Hirsch
\cite{FGH}, and \cite{GH} and when $\Gamma$ is distal by Fried \cite{F}.  

Carri\`ere introduced a set of notions about holonomy groups
\cite{Car}. In this paper, we will use the word 
``shrinkable dimension" 
for the French ``discompacit\'e" formerly translated 
as discompactedness \cite{rdaf}. 
Given a group of affine transformations, 
the shrinkable dimension 
measures the maximal codimension of the degenerating 
limiting sequences of ellipsoids under the affine action of 
the group (see Section \ref{sec:dconv}). That is, one 
measures the maximal ``shrinking" dimension of 
ellipsoids under the group action (see Figure 1).  

\begin{figure}[h] 
\centerline{\epsfysize=4cm
\epsfbox{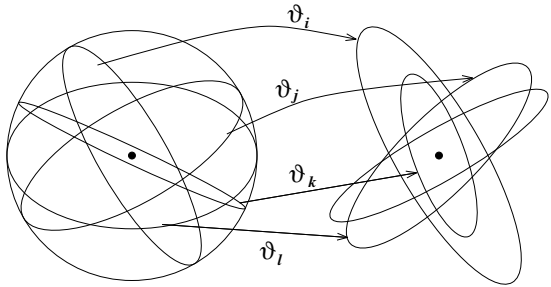}}
\caption{\label{fig:insh} The shrinkable dimension:
Linear maps $\vth_i$ mapping $2$-disks of radius $\eps$ 
in the standard $\eps$-ball to ellipsoids. We are provided 
with a fixed Euclidean metric on $\bR^3$.} 
\end{figure} 

Let $L: \Aff(\bR^n) \ra \GL(n, \bR)$ denote the 
homomorphism taking the linear part of affine transformations.
Carri\`ere \cite{Car} proved that if 
the shrinkable dimension of $\Gamma$,
or $\hbox{\rm sd}(L(\Gamma))$, for 
the holonomy group $\Gamma$ is $\leq 1$, 
then $\tilde M$ is affinely diffeomorphic to 
$\bR^n$ or a half space $\bR^n_+$, given by 
$x_1 > 0$ where $x_1$ is an affine coordinate.
If in addition $L(\Gamma)$ is a subgroup of 
$\SLt(n, \bR)$, then $\tilde M$ is affinely diffeomorphic to $\bR^n$.
In particular, if $L(\Gamma)$ lies in 
$\SO(n-1, 1)$, then $M$ is a complete affine manifold,
which means that any closed flat Lorentzian manifold is complete.

Since the conjecture of Markus is very difficult to attack, 
Carri\`ere proposed a new set of conjectures that if the 
shrinkable dimension of 
$\Gamma$ is $\leq i$ for some some integer $0 < i < n-1$, 
then $\pi_j(\tilde M)$ is trivial for $j \geq i$ \cite{Car2}.
(See also \cite[Section 3.3]{Car}.)

In this paper, we will prove a partial answer to the 
above conjecture as a corollary to our main Theorem \ref{thm:main}: 
\begin{cor}\label{cor:final}
The universal cover of a closed affine
three-manifold with holonomy group of 
shrinkable dimension less than or equal to $2$ is 
diffeomorphic to $\bR^3$. 
\end{cor} 

We need to introduce the powerful notion 
of the (projective) completion that is defined for 
incomplete geometric structures due to Kuiper \cite{Kuiper}
(see Kamishima-Tan \cite{KT}).  
The affine space $\bR^n$ has a Euclidean metric $\mu$ and 
the distance metric $\bdd$ induced from it; we 
choose such a metric standard with respect to an affine coordinate system. 
The developing map $\dev$ of $M$ induces a Riemannian metric $\mu$ 
on $\tilde M$ from the Euclidean metric of $\bR^n$.
(An arc is {\em straight\/} if it is a geodesic for $\mu$.)
$\mu$ induces a distance metric $\bdd$ on $\tilde M$. 
The completion $\che M$ of $M$ is obtained by 
Cauchy completing $\tilde M$ with respect to $\bdd$. 
The set $\ideal{M} = \che M - \tilde M$ is said to 
be the {\em ideal set}. 
$\che M$ and $\ideal{M}$ are topologically independent of 
the choice of $\dev$. The developing map
$\dev$ and every deck transformation $\vth$ extends uniquely on
$\che M$ to a distance-nonincreasing map and a homeomorphism 
respectively. For convenience, the extended maps 
are denoted by the same symbols $\dev$ and $\vth$ 
respectively (see \cite{cdcr1}).

An {\em $m$-simplex \/} in $\che M$ is 
a compact $m$-ball $A$ imbedded in $\che M$ such that 
its manifold interior $A^o$ is included in $\tilde M$ and 
$\dev | A$ is a (topological) imbedding 
onto a closed nondegenerate affine $m$-simplex in 
$\bR^n$. We introduce the ``key" definition given by Carri\`ere 
\cite{Car} adopted in the language of completions. 
(For a definition of $i$-convexity for real projective 
manifolds, see \cite{iconv} or \cite{psconv}.)

\begin{defn}\label{def:mconv}
We say that $M$ is {\em $m$-convex\/}, $0 < m < n$, 
if the following condition holds:
Whenever $T \subset \che M$ is an $(m+1)$-simplex
with sides $F_1, \dots, F_{m+2}$ such that 
$T^o \cup F_2 \cup \cdots \cup F_{m+2}$ does not meet
$\ideal{M}$, then $T$ does not meet 
$\ideal{M}$ and is a subset of $\tilde M$. 
\end{defn}

\begin{figure}[h] 
\centerline{\epsfysize=3.25cm
\epsfbox{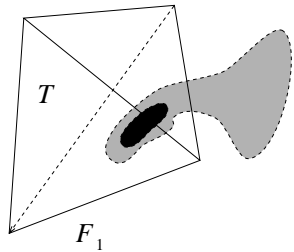}}
\caption{\label{fig:2conv} An example of a tetrahedron $T$ 
detecting not being 2-convex in $\che M$: 
The gray area denotes $\tilde M_\infty$ and 
the dark area $\ideal{M} \cap F_1^o = \ideal{M} \cap T$
in Figures 2 and 3.} 
\end{figure} 

\begin{prop}\label{prop:mapdconv}
Let $T$ be an affine $(m+1)$-simplex in $\bR^n$ 
with sides $F_1, \dots, F_{m+2}$.
Every affine immersion of $T^o \cup F_2 \cup \cdots \cup F_{m+2}$
into $M$ extends to one of $T$ if and only if 
$M$ is $m$-convex. 
\end{prop}

A heuristic idea of $2$-convexity is as follows:
Consider a closed room and remove from it all solid 
objects. Then the completion $\che A$ of the remainder $A$ 
is the union of $A$ and the boundary points.
Then $A$ will be $2$-convex if all solid objects 
do not have ``corners" or convex points or points 
which we can touch by an interior of a side of 
a tetrahedron, which we move around. 

\begin{figure}[h] 
\centerline{\epsfysize=5cm
\epsfbox{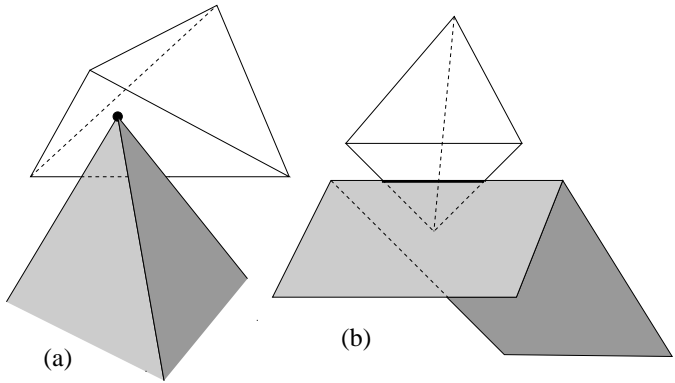}}
\caption{\label{fig:2conv2} The tetrahedron in 
(a) detects the failure of $2$-convexity 
while one as in (b) does not.} 
\end{figure} 

Any convex domain in $\bR^n$ is $m$-convex for every $m$, 
 
and the complement of disjoint lines in $\bR^n$ is $2$-convex. 
However, the complement of a compact convex set in $\bR^n$ is 
not $m$-convex for any $m$. (See Section 2 for explanations 
and more examples.) A last class of examples are given in Section 9 
(see Theorem \ref{thm:sihyp}).


The boundary $\delta M$ of an affine manifold $M^3$ has 
a {\em concave point} $p$ 
if a convex open totally geodesic $2$-disk meets $\delta M$ 
in a compact set containing $p$. A sufficient condition for $\delta M$ 
to have no concave point is that each point of $\delta M$ 
has a neighborhood affinely diffeomorphic to a convex 
domain in $\bR^3$ or to a $3$-dimensional submanifold $\Omega$ 
in $\bR^3$ with $\delta \Omega$ with 
principal curvatures of both signs when $\delta \Omega$ is 
considered as an isometrically imbedded two-dimensional 
submanifold of the Euclidean space $\bR^3$.  

One can find $(m+1)$-simplices in $\che M$ 
inside an affine coordinate neighborhood in $\tilde M$, and 
there are many more $(m+1)$-simplices in $\che M$ 
(see Remark~\ref{rem:magnify}). The following lemma 
gives us our {\em principal}\/ criterion of $m$-convexity. 

\begin{prop}\label{prop:mconv}
$M$ is not $m$-convex if and only if there exists 
an $(m+1)$-simplex $T$ with a side $F_1$ such that 
$T \cap \ideal{M} = F_1^o \cap \ideal{M} \ne \emp$.  
\end{prop}
\begin{pf}
Suppose that every $(m+1)$-simplex $T$ has the property that 
$T$ do not meet $\ideal{M}$; 
$T \cap \ideal{M}$ is a subset of the union of two or more sides
but not less than two sides; or 
if $T \cap \ideal{M}$ is a subset of a side $F$, then $F \cap \ideal{M}$ 
is not a subset of $F^o$. Then one see easily that 
the definition for $m$-convexity is satisfied by $M$. 

Conversely, if $M$ is $m$-convex, and 
$T$ is an $(m+1)$-simplex with sides $F_1, \dots, F_{m+2}$ such that 
$T^o \cup F_2 \cup \dots \cup F_{m+2} \subset \tilde M$,
then $T \subset \tilde M$. Consequently, there is no 
$(m+1)$-simplex $T$ with $T \cap \ideal{M} = F_1^o \cap \ideal{M} 
\ne \emp$.
\end{pf}

Clearly, $M$ is $m$-convex if and only if 
any of the covers of $M$ is $m$-convex.  
It is easy to see that $1$-convexity is 
equivalent to convexity, and $i$-convexity implies 
$j$-convexity whenever $i \leq j \leq n$ 
(see Section 2). 

Our main theorem is:
\begin{thm}\label{thm:main} 
Let $M$ be a two-convex affine three-manifold. 
Then $\tilde M$ is diffeomorphic to $\bR^3$.    
\end{thm}

Note above that $M$ does not need to be a closed $3$-manifold
in the premise of the theorem:
$M$ could be an open, noncompact $3$-manifold. 

\begin{thm}\label{thm:dconv}
Let $M$ be a closed $n$-dimensional affine manifold. 
If the holonomy group of $M$ has 
shrinkable dimension 
less than or equal to $d$\cO then $M$ is $d$-convex. 
\end{thm}
The above theorem almost 
answers another conjecture of Carri\`ere 
\cite[Section 3.3]{Car}. 
The proof is an essentially straightforward generalization of 
Carri\`ere's argument \cite{Car} with 
some modification involving the definition of 
shrinkable dimension. (However, we feel that our definition 
is just as useful in application.) 

The above theorem combined with Theorem \ref{thm:main} 
implies Corollary \ref{cor:final}.

The holonomy group of an affine three-manifold $M$ with parallel 
volume form has shrinkable dimension less than or equal to $2$  
(see equation \ref{eqn:collapse}) so we obtain: 
\begin{cor}\label{cor:para} 
Let $M$ be a closed affine three-manifold with 
parallel volume form. Then the universal cover of $M$ is 
diffeomorphic to $\bR^3$.
\end{cor}

\begin{cor}
Let $M$ be a closed $n$-dimensional affine manifold 
with holonomy of shrinkable dimension 
less than or equal to two $(n \geq 4).$ Then
the universal cover $\tilde M$ is foliated by totally geodesic 
leaves which are diffeomorphic to $\bR^3/\Gamma$ for $\Gamma$ 
depending on the leaves. 
\end{cor}
\begin{pf}
Obviously $\tilde M$ is foliated by 
$3$-dimensional leaves which are defined by 
setting the first $n-3$ coordinates 
of $x_i \circ \dev$, $i=1, \dots, n$, equal to constants. 
Let $l$ be a leaf. Then since $M$ is $2$-convex, 
$l$ is a $2$-convex affine $3$-manifold by Proposition \ref{prop:mapdconv}. 
Hence, the universal cover of $l$ is 
diffeomorphic to $\bR^3$ and $l$ is    
diffeomorphic to $\bR^3/\Gamma$.
\end{pf}

We are curious how much the above tells you 
about the topology of $\tilde M$. In fact, 
we may choose any $n-3$ coordinate directions 
to obtain the above result. Our conjecture is 
that the universal cover of $M$ is diffeomorphic to $\bR^n$.
In particular, our conjecture implies that 
the universal cover of a closed affine $4$-manifold with 
parallel symplectic form is diffeomorphic to $\bR^4$.

As an application of Theorem \ref{thm:main}, 
we see that the universal covers of the submanifolds 
listed in Proposition \ref{prop:examples} are diffeomorphic to 
$\bR^3$ respectively. By the result of Gromov \cite{Gromov}, 
any compact $3$-manifold $M$ with nonempty boundary admits 
an affine structure with 
trivial holonomy. We see that if $M$ is reducible, 
then the boundary of $M$ always has a concave point, 
and given any affine structure on a fake cell, 
the boundary of it always has a concave point.
(Here is an example of geometry and topology interacting.)

Lastly, we have the applications to singular hyperbolic manifolds.
Singular hyperbolic manifolds with cone-type singularity will 
be defined in Section \ref{sec:sihyp} as was done  
by Thurston (see Hodgson \cite{Hodgson}.)

\begin{thm}\label{thm:sihyp}
Let $M$ be a singular hyperbolic $3$-manifold with cone-type 
singular locus along a link $L$. 
Then the universal cover of $M-L$ is 
diffeomorphic to $\bR^3$.  
\end{thm}
This is also implied by Theorem 1.2.1 of Kojima \cite{Koj2}, 
which is a stronger result (see also Kojima \cite{Koj1})
since he shows $M-L$ is atoridal.

We also make the following conjecture with a proof to be supplied 
shortly: Let $M$ be a singular hyperbolic $3$-manifold with cone-type
singular locus along a graph $K$ with cone-angle $\leq 2\pi$. 
Then the universal cover of $M -K$ is diffeomorphic to 
$\bR^3$.

We think that this conjecture has some relation to a recent conjecture 
of Casson announced in the summer of 1998 at Tokyo Institute of 
Technology that the cone-angles of such manifolds can always be 
infinitesimally deformed.


In this paper, we will prove Theorems 1, 2, and 3, as the other
propositions were proved already in this introduction.

In Section 1, we review materials on completions needed 
in this paper. In Section 2, we list examples on 
$d$-convexity and discuss $d$-convexity. We prove 
Propositions \ref{prop:mapdconv} and \ref{prop:examples}.

The purpose of Sections 3 to 7 is to prove Theorem \ref{thm:main},
Section 8 to prove Theorem \ref{thm:dconv} and Section 9 
to prove Theorem \ref{thm:sihyp}.
{\em For simplicity, we assume that $M = \tilde M$
or $M$ is simply connected from Sections 3 to 7, and  
we assume otherwise in Sections 1, 2, 8, and 9.}

Since $M$ is simply connected in Sections 3 to 7, 
the developing map $\dev$ is defined on $M$ itself.
Let $M$ have a coordinate 
function $x_1$ given by composing the developing map
$\dev: M \ra \bR^3$ 
to $\bR^3$ with an affine coordinate function $x_1$ of 
$\bR^3$. This determines level sets, which are open surfaces, 
in $M$ with $x_1^{-1}(t)$ denoted by $M_t$.  

Our plan to prove Theorem \ref{thm:main} is
to show geometric incompressibily, i.e. 
$2$-convexity, implies topological incompressibility
of level sets. That is, $2$-convexity implies that 
each component of the level set in $M$ is 
incompressible. Since $M$ is simply connected, each component 
is homeomorphic to a disk. A result of Palmeira \cite{Pal}
says that a simply connected $3$-manifold foliated by 
leaves diffeomorphic to $\bR^2$ is diffeomorphic to $\bR^3$. 
Thus, we obtain our result that $M$ is diffeomorphic to $\bR^3$.  
To prove that each component surface $F$ of a level set 
is incompressible, we show that there are no compressing disks
of $F$. An {\em amenable disk} is an embedded disk where the height 
function $x_1$ is Morse. We deform any disk with boundary in $F$ to 
an amenable disk and we cancel critical points 
to homotope it to an imbedded disk in $F$.
This is accomplished by showing that $2$-convexity implies  
the following properties:
\begin{itemize}
\item A property associated with a $3$-simplex in $\che M$
(Section 3). 
\item Filling of a drum with the bottom side removed
(see Theorem \ref{thm:fill} in Section 4).   
\item Filling of an amenable disk with only one 
critical point and annuli with two critical points
(Section 5). 
\item A homotopy of an amenable disk with only 
three critical points into a disk in $F$ (Section 6). 
\item A homotopy of an amenable disk with 
arbitrarily many critical points into a disk in $F$ (Section 7). 
\item The incompressibility of each component of $M_t$.
\end{itemize} 

We note that this method of Morse cancellation was already employed 
by Flyod-Hatcher \cite{FH} in their work on punctured-torus
bundles over circles. (Hatcher-Thurston \cite{HT} shortly 
generalized this result to two-bridge knots.)
Our method only differs from theirs slightly and we provide 
some detailed analysis steps missing from their work
but needed in this paper (mainly in the Appendix).

In Section 3, as a preliminary result, we show that 
given a $2$-convex affine $3$-manifold $M$,
if $T \subset \che M$ is a $3$-simplex with 
a side $F_1$ meeting $\ideal{M}$, then there is 
no exposed point of $\ideal{M} \cap F_1$ 
in $F_1^o$. (Note that even when $M$ is $2$-convex 
a $3$-simplex may still meet $\ideal{M}$ in a side, 
but not in the interior of a side.) This will be used to prove that 
a solid-cylinder map $f: \Omega \times [a, b] \ra \che M$ 
can only be a map into $\tilde M$. 

\begin{figure}[h] 
\centerline{\epsfysize=2.75cm
\epsfbox{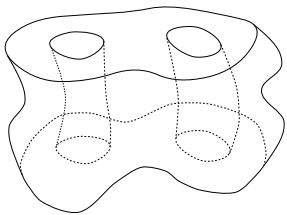}}
\caption{\label{fig:uf3} A figure of a drum} 
\end{figure} 

Section 4 is on the filling of a drum:
An {\em embedding \/} from 
a smooth manifold (perhaps with corners) to another 
smooth manifold is a topological imbedding,
which is a restriction of a smooth immersion in 
some ambient manifolds.
(The word ``imbedding" will be reserved for topological 
imbeddings.)
A {\em drum with the bottom side removed}
is the union of an embedded image of a surface $\Omega$ 
into a level $b$ and the image of embedding from 
$\delta \Omega \times [a, b]$ to $M$ mapping $(x, t)$ into 
level $t$ for each $t \in [a, b]$ which is transversal 
to the level sets. 
The following process describes the process of 
the {\em filling} of a drum with the bottom side removed. 
Let $\Omega$ be the surface with smooth boundary in the plane $\bR^2$. 
We show that if there is an embedding 
from $\delta \Omega \times [a, b]$ to $M$  
preserving levels, that is, mapping 
$(x, t)$ to a point of $M_t$ for each $t \in [a, b]$, 
and if the embedding restricted to $\delta \Omega \times \{b\}$
extends to one from $\Omega \times \{b\}$ 
to $M_b$, then the embedding extends to 
one from $\Omega \times [a, b]$ into $M$ preserving levels.   

The proof that the filling of a drum is always possible 
uses the flow argument based on noticing 
the flow argument can be applied to ideal boundary 
in some limited sense.

In Sections 5, 6, and 7,
the result of Section 4 will be modified further to show that 
if an embedded disk $D$ in $M$ has the boundary in 
a level set $M_t$, then the boundary bounds a disk in $M_t$. 
Any embedding $F$ of a surface $\Omega$ into $M$
where $F(\delta \Omega) \subset M_t$ for some $t$ can be 
modified to a so-called amenable embedding $F': \Omega \ra M$
of level $t$ such that $x_1 \circ F'$ is a Morse function and 
the boundary $\delta \Omega$ maps into $M_t$.
Our strategy is to reduce the number of critical points by 
cancellation.

By a {\em filling} of an amenable imbedding 
$F: \Omega \ra M$ of level $t$, we mean the process of 
finding a subsurface $\Omega'$ in $M_t$ and a $3$-dimensional 
submanifold $N$ of $M$ such that $\Omega$ and 
$\Omega'$ meet in their common boundary and 
$\Omega \cup \Omega'$ forms the boundary of $N$.
Here $\Omega'$ is unique due to Lemma \ref{mt:open}. 

In Section 5, we show that given an amenable imbedding 
of a disk with one critical point, then the boundary 
bounds a disk in $M_t$. Then we study amenable imbeddings of 
annuli with two critical points. We show that the amenable disk 
and annuli are fillable.  

In Section 6, we will show that given an amenable disk $D$
of level $t$ with three critical points, $\delta D$ bounds 
a disk in $M_t$. This is the most technical part of the paper 
and the reason is that we want to cancel singularities 
without introducing new ones, and hence, we need to construct 
deformations carefully. However, the manner in which we 
do is entirely similar to what one would do in 
p.l. $3$-manifold theory consisting of various disk moves 
and isotopies. 

In Section 7, we do an inductive argument by 
reducing the number of critical points of amenable disks,
showing that given an amenable disk $D$ of level $t$, $\delta D$ bounds
a disk in $M_t$. 
This result, the incompressible surface theory, and 
Corollary 3 of Palmeira \cite{Pal}   
shows that $M$ is diffeomorphic to $\bR^3$. 
This completes the proof of Theorem \ref{thm:main}.

In Section 8, we discuss shrinkable dimension. 
We will give examples of shrinkable dimensions for 
some Lie groups. The basic idea behind this definition is that 
along the direction of real maximal tori the affine action
of the given Lie group can shrink and stretch indefinitely 
and in the other directions, the Lie group action can do so 
only finitely. {\em We now assume that $M$ is a closed affine $n$-manifold}. 
We give the proof of Theorem 
\ref{thm:dconv}. The main idea is to follow 
Carri\`ere's argument \cite{Car} based on infinite rays.
However, we need to perturb the $(d+1)$-simplices sometimes
to control the size of the ideal set with rather involved argument. 
(We point out that Carri\`ere's argument \cite{Car} does not 
carry out this step and hence we believe it
contains a negligible gap, which we fill here.)

In Section 9, we prove Theorem \ref{thm:sihyp}. 
Since the universal cover $(M-L)\tilde{}$ has a hyperbolic structure, 
$(M-L)\tilde{}$ develops into the standard disk in $\rp^3$. 
Since the standard disk is in the affine patch of $\rp^3$, 
it follows that $(M-L)\tilde{}$ has an affine structure. 
The set of ideal points $(M-L)_\infty\tilde{}$ consists 
of lines corresponding to $L$ and points of the sphere of 
infinity $\SI^2_\infty$. A large part of Section 9 is 
devoted to proving that $(M-L)\tilde{}$ is $2$-convex.
By Theorem \ref{thm:main}, we obtain Theorem \ref{thm:sihyp}. 

In the appendix, we prove technical Lemma \ref{lem:coordfun} needed in
Section \ref{sec:hdthree}. That is, we show that there exists 
a suitable coordinate system on a connected neighborhood 
including a saddle and a relative maximum point so that 
our Morse function takes an integral form.

\section*{Acknowledgment} 
I wish to thank Boris Apanasov, 
Yves Carri\`ere, Dongho Chae, David Gabai, Bill Goldman, 
Morris Hirsch, Michael Kapovich, Hyuk Kim, Sadayoshi Kojima, 
Shigenori Matsumoto, Darryl McCullough, Peter Scott and 
Bill Thurston for many interesting discussions. 
In particular, this research was started by a conversation 
with Yves Carri\`ere in an AMS meeting in Kansas, 1989.
I thank our department computer assistant Sangeun Lee 
for numerous hours and hard work spent on setting up my Linux 
system. This paper is prepared with a Linux system and 
the figures were drawn by xfig. 


\section{Properties of Completions}

In this section, we will capture the necessary facts about 
the completions due to Kuiper \cite{Kuiper} needed in this paper.
Let $M$ be an affine $n$-manifold, which is not necessarily 
simply connected.
Proposition \ref{prop:extmap} shows that we can tell precisely 
about how two balls meet in $\che M$ by looking 
at the way that their developed images meet. 
Remark \ref{rem:magnify} is used 
in many occasions in this paper. The two are basic tools 
of this paper and other papers (see \cite{cdcr1}, \cite{cdcr2}, 
\cite{cdcr3}, and \cite{psconv}). 

Associated with an affine structure on 
an $n$-manifold $M$, there is a development pair
$(\dev, h)$ where $\dev: \tilde M \rightarrow \bR^n$ is an immersion 
and $h: \pi_1(M) \rightarrow \Aff(\bR^n)$ is a homomorphism such that 
$h(\vth) \circ \dev = \dev \circ \vth$ for $\vth \in \pi_1(M)$. 
$\dev$ is obtained by an analytic continuation of 
charts of $\tilde M$. 
Given every other such pair $(\dev', h')$ for the given 
affine structure on $M$, we have
\[ (\dev', h'(\cdot)) = (\vpi \circ \dev, \vpi\circ h(\cdot) 
\circ \vpi^{-1}) \mbox{ for } \vpi \in \Aff(\bR^n). \] 
The {\em holonomy group} is the image of a holonomy homomorphism.
The pair $(\dev, h)$ determines the affine structure on $M$.

A {\em convex subset\/} of $\tilde M$ 
is a subset $A$ such that every two points of $A$ 
can be connected by a segment in $A$.  
Obviously, if $A$ is a convex subset of $\tilde M$, 
then $\dev|A$ is an isometry with respect to 
the metrics both denoted by $\bdd$ on $\tilde M$ and 
$\bR^n$. Hence, $\dev|A$ is a imbedding onto 
a convex subset of $\bR^n$. 

\begin{prop}\label{prop:extconv}
Let $A$ be a convex subset of $\tilde M$.  
Then for the closure $\clo(A)$ of $A$ in $\che M$
$\dev|\clo(A)$ is an imbedding onto 
a closed convex set $\clo(\dev(A))$. 
Moreover\cO if $A$ is $\bdd$-bounded, then 
$\clo(A)$ is a $\bdd$-bounded compact subset of $\che M$. 
\end{prop} 
\begin{pf}
The map $\dev|\clo(A)$ is a $\bdd$-isometry. 
The lemma follows. 
\end{pf}

A {\em tame} subset of $\che M$ is 
a compact ball $A$ imbedded in $\che M$ such that 
$\dev| A$ is an imbedding onto a convex 
subset of $\bR^n$ and $A^o \subset \tilde M$.  
For example, a simplex is tame.

\begin{prop}\label{prop:tame}
The closure of a $\bdd$-bounded 
convex subset of $\tilde M$ is a tame subset of $\che M$. 
Conversely\cO the interior of a tame 
subset is a convex subset of $\tilde M$. 
\end{prop}

\begin{defn}\label{def:kball}
A {\em $k$-ball} in $\che M$ is a subspace homeomorphic 
to a $k$-ball such that its interior is a subset of $\tilde M$
and $\dev| A$ is an imbedding onto a $k$-ball in 
a $k$-dimensional affine subspace $P^k$ of $\bR^n$. 
\end{defn}

\begin{prop}\label{prop:extmap}
Let $A$ be a $k$-ball in $\che M$\cO and $B$ and $l$-ball.
Suppose that $A^o \cap B^o \ne \emp$\cO 
$\dev(A) \cap \dev(B)$ is a compact manifold 
in $\bR^n$ with interior equal to $\dev(A)^o \cap \dev(B)^o$\cO 
and $\dev(A)^o \cap \dev(B)^o$ is arcwise-connected.  
Then $\dev| A \cup B$ is a homeomorphism 
onto $\dev(A) \cup \dev(B)$. 
\end{prop}
\begin{pf}
First, we prove injectivity: 
Let $x \in A$, $y \in B$, and $z \in A^o \cap B^o$. 
Suppose that $\dev(x) = \dev(y)$. 
There is a path $\gamma$ in $\dev(A) \cap \dev(B)$ 
from $\dev(z)$ to $\dev(x)$ such that 
$\gamma| [0, 1)$ maps into $\dev(A)^o \cap \dev(B)^o$.
Since $\dev| A$ is an imbedding onto $\dev(A)$, 
there is a lift $\gamma_A$ of $\gamma$
into $A^o$ from $z$ to $x$. Similarly, there is 
a lift $\gamma_B$ of $\gamma$ into $B^o$ 
from $z$ to $y$. Note that $\gamma_A| [0, 1)$ 
agrees with $\gamma_B| [0, 1)$ by the uniqueness of 
lifts of arcs for immersions. Thus, $x=y$ 
since $\che M$ is a complete metric space.

By injectivity, there is a well-defined 
inverse function $f$ to $\dev| A \cup B$. 
$f$ restricted to $\dev(A)$ equals 
the inverse map of $\dev| A$, 
and $f$ restricted to $\dev(B)$ 
the inverse map of $\dev| B$. 
Since both inverse maps are continuous,
and $\dev(A)$ and $\dev(B)$ are closed, 
$f$ is continuous. (See also Proposition 1.3.1 in \cite{Car}.)  
\end{pf}

\begin{prop}\label{prop:extseq}
Let $\{K_i\}$ be a sequence of uniformly 
$\bdd$-bounded tame subsets of 
$\che M$ where $K_i \subset K_{i+1}$ for each $i,$ $i=1, 2, \dots$.  
Then the closure $K'$ of 
$\bigcup_{i=1}^\infty K_i$ is a tame subset of 
$\che M$\cO and $\dev(K')$ equals 
the closure of $\bigcup_{i=1}^\infty \dev(K_i)$. 
\end{prop}
\begin{pf}
Straightforward.
\end{pf}     

\begin{rem}\label{rem:magnify}
We use the above proposition in the following form:
(See Carri\`ere \cite[Section 3]{Car} also.)
Let $T'$ be an affine $m$-simplex in $\bR^n$
and $x$ a point of $\tilde M$ such that 
$\dev(x)$ is a point of $\delta T'$. 
Let $T'_t$ denote the image of $T'$ under the action 
of dilatation with center $\dev(x)$ by a factor $t$, $0 < t \leq 1$.
Let $A$ be the subset of $(0, 1]$ such that for $t \in A$
there exists a tame set $T_t$ in $\tilde M$ 
such that $x \in \delta T_t$ and $\dev(T_t) = T'_t$.  
Then we can show 
\begin{enumerate}
\item[(a)] $A$ includes $(0, \eps)$ for a small positive 
number $\eps$.  
\item[(b)] $T_t \subset T_{t'}$ if $t < t'$, $t, t' \in A$, 
and for each $t$, $T_t$ is unique.  
\item[(c)] $A$ is an open subset of $(0, 1]$. 
\item[(d)] $A$ is of form $(0, \tau)$ for 
a positive number $\tau$ or equals $(0, 1]$.  
\item[(e)] The closure $T_A$ of $\bigcup_{t \in A} T_t$ is a tame 
set such that $\dev| T_A$ is an imbedding onto 
the closure of $\bigcup_{t \in A} T'_t$. 
Thus $T_A$ is an $m$-simplex in $\che M$
(however, not in $\tilde M$ in general).
\end{enumerate} 
\end{rem}
\begin{pf}
(a) For a coordinate neighborhood $U$ of $x$, $\dev| U$ is 
a coordinate chart of $x$. Assume for convenience, 
that $\dev(U)$ is a convex open ball.  
Hence for $\eps$ small enough, there exists 
$T_t$ in $U$ for $t < \eps$. 

(b) Since each $T_t$ contains $x$ in $\delta T_t$,  
$T_t$ includes $T_k$ for small enough $k$ so that $T_k \subset U$.
This follows since the maximal totally geodesic connected
submanifold including $T_t$ or $T_k$ is determined by 
its tangent space at $x$.
Hence, $T_t$ and $T_{t'}$ includes $T^o_k$ for $k < \eps, t, t'$. 
By Proposition \ref{prop:extmap}, this means $T_t \subset T_{t'}$
since $\dev(T_t) \subset \dev(T_{t'})$.  
The uniqueness follows from this statement.

(c) Since $T_t$ is a compact subset of $\tilde M$, and 
$\dev| T_t$ is an imbedding, there exists a compact 
neighborhood $U$ of $T_t$ such that $\dev| U$ is an imbedding 
to its image. Hence, for $t'$ sufficiently near $t$, 
$U$ includes $T_{t'}$, and $t' \in A$.

(d) $A$ is connected since $T_t \subset T_{t'}$ if $t < t'$ 
for $t, t' \in A$. Thus (d) follows from (a) and (c).  

(e) follows from Proposition \ref{prop:extseq}.
\end{pf}

\section{$d$-convexity: Examples}

We start this section by showing that $1$-convexity is 
equivalent to convexity. Next, we give some examples of 
$m$-convex manifolds and prove Propositions \ref{prop:1convex}
\ref{prop:examples}. We also prove Proposition \ref{prop:mapdconv}
from the introduction.
Lastly, we show that $i$-convexity implies 
$j$-convexity if $j \geq i$.

The affine manifold $M$ is said to be {\em convex} if 
its universal cover $\tilde M$ is a convex set.
The following proposition was proved by Carri\`ere (see \cite{Car}).
We give a brief proof based on our definition of 
$1$-convexity.

\begin{prop}\label{prop:1convex} 
The following statements are equivalent.
\begin{enumerate}
\item[(1)] $M$ is $1$-convex.
\item[(2)] $M$ is convex.
\item[(3)] $M$ is affinely diffeomorphic to 
$\Omega/\Gamma$ where $\Omega$ is a convex domain in 
$\bR^n$ and $\Gamma$ is a group of affine transformations
acting on $\Omega$ freely and properly discontinuously.
\end{enumerate}
\end{prop}
\begin{pf}
We show that (1) implies (2). 
Suppose we are given two segments $s_1$ and $s_2$ in $\tilde M$
with endpoints $x$ and $y$ of $s_1$ and 
$z$ and $t$ of $s_2$. Suppose $t = y$ and 
$s_1$ and $s_2$ is in general position.
We will show that there exists a segment in 
$\tilde M$ with endpoints $x$ and $z$. 
Let us give affine parameters $f_1$ and $f_2: [0, 1] \ra \tilde M$ 
to $s_1$ and $s_2$ 
respectively so that $f_1(0) = y$, $f_2(0) = y$, $f_1(1) = x$, 
and $f_2(1) = z$. For $t \in (0, 1]$, define $T_t$ be 
the triangle in $\tilde M$ with vertices 
$f_1(t), f_2(t),$ and $y$ if it exists. 
By Remark \ref{rem:magnify}, we see that 
$A$ is of form $(0, \tau)$ for $\tau > 0$
or equal $(0, 1]$.

Suppose $A$ of form $(0, \tau)$.
The closure $T'$ of $\bigcup_{t \in A} T_t$ is 
a triangle in $\che M$.
Since the two edges of $T'$ are subset of $s_1$ and $s_2$, 
they do not meet $\ideal{M}$. Since $M$ is $1$-convex, 
this means $T' \subset \tilde M$ and $\tau \in A$.  
This is a contradiction, and $A = (0, 1]$.

This means that $M$ includes a triangle $T$ with 
vertices $x, y, z$, and $x$ and $z$ is connected 
by a segment.

Let $x$ and $y$ be arbitrary points of $\tilde M$.
Then $x$ and $y$ are connected by a broken geodesic. 
The above result tells you that we can find a broken 
geodesic with less number of segments than our initial choice.
An induction shows that $x$ and $y$ are connected by a segment. 

(2) Since $\tilde M$ is a tame set, 
$\dev|\tilde M \ra \dev(\tilde M)$ is an imbedding, 
and hence (3) follows.

That (3) implies (1) follows from the fact that 
$\tilde M$ is affinely diffeomorphic to a convex 
domain in $\bR^n$.
\end{pf}

We will give some examples illustrating $m$-convexity. 
Initial examples are affine manifolds 
whose developing maps are imbeddings.

A convex domain is easily seen to be $i$-convex for each $i$
since the convex hull of any finite set of points of 
the domain is still in the domain.
By the above proposition, every convex affine manifold $M$ is 
$i$-convex for every $i$ since $\tilde M$ is affinely 
diffeomorphic to a convex domain in $\bR^n$.

Let $x_1, x_2$, and $x_3$ denote the standard coordinates 
of $\bR^3$. Remove from $\bR^3$ the set given 
by $x_1 \geq 0, x_2 \geq 0, x_3 \geq 0$. 
Then the remainder is an affine manifold $M$ which is simply 
connected. Clearly, $\tilde M = M$ and $\che M$ can be identified 
with $\bR^3$ removed with the set given by $x_1 > 0, x_2 > 0, x_3 > 0$. 
Then $M$ is not $2$-convex, since the $3$-simplex 
given by $x_1, x_2, x_3 \geq -1$ and $x_1 + x_2 + x_3 \leq 0$ in 
$\che M$ has unique side whose interior intersects 
$\ideal M = \che M - \tilde M$ (see Proposition \ref{prop:mconv}
and Figure \ref{fig:uf3}). 

Remove from $\bR^3$ the set given by $x_1 \geq 0, x_2 \geq 0$.
Then the remaining affine manifold $M$ is simply connected 
and $\tilde M = M$ and $\che M$ can be identified with 
$\bR^3$ removed with the set given by $x_1 > 0, x_2 > 0$. 
Here $\ideal{M}$ equals the union of the set given by 
$x_1 = 0$ and $x_2 \geq 0$ and the set given by 
$x_2 = 0$ and $x_1 \geq 0$.
Then $M$ is $2$-convex since given any $3$-simplex $T$ in $\bR^3$ 
meeting $\ideal{M}$ in the interior of a side $F_1$, $T^o$ or 
the boundary of $F_1$ meets $\ideal{M}$. $M$ is not $1$-convex since 
there exists a $2$-simplex $F$ in $\che M$ with an edge $e$ 
intersecting the line given by $x_1= x_2 = 0$ 
in a point of $e^o$ and $F$ meets $\ideal{M}$ exactly 
at this point only. 
 
A {\em wedge} in $\bR^3$ is a convex set given by 
inequalities $l_1(x) \geq c_1$ and $l_2(x) \geq c_2$ 
for affine functions $l_1$ and $l_2$ and constants $c_1$ and $c_2$,
where $l_1(x) =c_1$ and $l_2 = c_2$ define distinct planes respectively.
If we remove from $\bR^3$ a collection of disjoint wedges, it follows 
easily that the remaining affine manifold is $2$-convex.

A {\em convex cone} in the affine space $\bR^3$ is a convex domain $\Omega$ 
with a distinguished point $p$ such that every point of 
$\Omega$ other than $p$ has an infinite ray in $\Omega$ 
from $p$ passing through it.
A {\em proper convex} cone is a cone that has no complete 
affine line in it.  If we remove from $\bR^3$, a collection of disjoint
properly convex cones, the remaining affine manifold is not $2$-convex. 

Another example is an affine manifold in 
$\bR^3$ given by $x_1^2 - x_2^2 - x_3^2 < -1$. 
This is $2$-convex but not $1$-convex. 
The affine manifold given by $x_1^2 - x_2^2 - x_3^2 > 1$ is 
$1$-convex and $2$-convex.

Lastly, we mention that given a closed real 
projective surface $\Sigma$ of negative Euler characteristic, 
a trivial circle bundle $M$ over $\Sigma$ carries an affine 
structure by Benz\'ecri suspension (see \cite{Ben2}, \cite{Car3}, 
and \cite{Sim}). 
If $\Sigma$ is not convex, then  $M$ is not convex but $2$-convex.
We omit the simple proof of this fact.

We prove Proposition \ref{prop:mapdconv}:
Suppose that $M$ is not $m$-convex, and let $p: \tilde M \ra M$ 
denote the universal covering map.
Then there exists an $(m+1)$-simplex $T$ in $\che M$ such 
that $T \cap \ideal{M} = F_1^o \cap \ideal{M} \ne \emp$ 
for a side $F_1$ of $T$.
Since $T$ is an $(m+1)$-simplex, $\dev| T$ is an imbedding 
onto $\dev(T)$ and $\dev(T)$ is an $(m+1)$-simplex in $\bR^n$. 
Let $f$ be the map $(\dev|\dev(T))^{-1}$ restricted to 
$\dev(T^o) \cup \dev(F_2) \cup \dev(F_3) 
\cdots \cup \dev(F_{m+2})$.
$f$ is an affine immersion to $\tilde M$. 
It is easy to see that $f$ does not extend to all of $\dev(T)$.
Hence, $p\circ f$ is an affine immersion which does not 
extend to $\dev(T)$.

Suppose that $M$ is $m$-convex. Let $\dev: \che M \ra \bR^n$ 
be the developing map. 
Let $f: T^o \cup F_2 \cup \dots \cup F_{m+2}$ be an affine 
immersion into $M$ where $T$ is an $m$-simplex 
with sides $F_1, \dots, F_{m+2}$ in $\bR^n$. 
Let $\tilde f$ be the lift of $f$ to $\tilde M$.
Then $\dev\circ \tilde f$ is also an affine map from 
$T^o \cup F_2 \cup \dots \cup F_{m+2}$ into $\bR^n$.
Since the rank of the map $\dev \circ \tilde f$ is maximal, 
the map extends to a global affine transformation on $\bR^n$, and 
$\dev \circ \tilde f$ is an imbedding.
Hence, $\dev| \tilde f(T^o)$ is an imbedding onto 
an open $m$-simplex $\dev \circ \tilde f(T^o)$ in $\bR^n$.
Since $\tilde f(T^o)$ is a convex subset of $\tilde M$,
the closure $T'$ of $\tilde f(T^o)$ is a tame subset of 
$\che M$ such that $\dev| T'$ is an imbedding onto 
$\clo(\tilde f(T^o))$ by Proposition \ref{prop:tame}. 
Since $T'$ is an $m$-simplex with sides 
$\tilde f(F_2), \dots, \tilde f(F_{m+2})$ and a remaining 
side $F'_1$; $\tilde f(T^o)$ and 
$\tilde f(F_2), \dots \tilde f(F_{m+2})$ 
are subsets of $\tilde M$; and $M$ is $m$-convex,
we have $T' \subset \tilde M$. Clearly, $\tilde f$ extends to 
an affine embedding $f': T \ra T'$, and $p\circ f'$ extends $f$. \qed  

\begin{prop}\label{prop:examples}
\begin{itemize}
\item[(a)] Let $\Omega$ be a convex open domain in $\bR^3$.
Then $\Omega$ removed with the closed set that is the union of 
complete lines (not necessarily finitely many) is $2$-convex. 
\item[(b)] The interior of a closed $3$-dimensional submanifold $M$ 
of $\bR^3$ with smooth boundary without a concave 
point is $2$-convex.
\item[(c)] The interior of a compact affine manifold $M$ 
with trivial holonomy and smooth boundary without 
a concave point is $2$-convex.
\end{itemize}
\end{prop}
\begin{pf}
(a) Let $\Omega'$ be the convex domain $\Omega$ in $\bR^3$ 
removed with the closed set that is the union of
a collection of complete lines. 
Let $T$ be an affine $3$-simplex with sides 
$F_1, F_2, F_3,$ and $F_4$ and an affine immersion 
$f: T^o \cup F_2 \cup F_3 \cup F_4$ into $\Omega'$.  
Since the rank of the affine map $f$ is maximal, 
$f$ is an imbedding and the image equals 
$T^{\prime o} \cup F'_2 \cup F'_3 
\cup F'_4$ for an affine $3$-simplex $T'$ in $\bR^3$ with 
sides $F'_1, F'_2, F'_3, F'_4$. Since the vertices of $T'$ are 
in $\Omega$, $T'$ is a subset of $\Omega$. 

Suppose that $T'$ is not a subset of $\Omega'$.
Since $\Bd \Omega'$ is the union of a sphere $\Bd \Omega$ 
and the union of removed lines, if $T'$ meets $\Bd \Omega'$, then 
$T'$ meets $\Bd \Omega$ or a line $l$.
Since $\Omega$ is convex and vertices of $T'$ are in $\Omega$, 
$T'$ do not meet $\Bd \Omega$. Thus $T'$ meets $l$, one of 
the removed lines; more precisely, $F_1^{\prime o}$ meets $l$. 
$l$ meets $F_1^{\prime o}$ tangentially since 
otherwise $l$ meets $T^{\prime o}$, which is a contradiction.
Since $l$ has endpoints in $\Bd \Omega$, 
$l$ must meet the boundary of the $2$-simplex $F_1$, 
which is compact in the open set $\Omega$. 
But since $\delta F_1$ is a subset of $F_2 \cup F_3 \cup F_4$, 
this is absurd. Therefore, $T' \subset \Omega'$ and 
$f$ extends to an affine embedding $f': T \ra T'$. 
By Proposition \ref{prop:mapdconv}, $\Omega'$ is $2$-convex.

To prove (b), it is sufficient to prove (c).
Let $M$ be a compact affine manifold with smooth boundary and 
trivial holonomy. We assume that $\delta M$ has no concave point. 
Since the holonomy of $M$ is trivial, 
the developing map $\dev: \tilde M \ra \bR^3$ factors
into the composition of the covering map $p: \tilde M \ra M$ and 
$\dev': M \ra \bR^3$. Thus the induced metric $\bdd$ on 
$\tilde M$ agrees with that induced from the metric 
$\bdd'$ on $M$ where $\bdd'$ is induced from $\bR^3$ by 
$\dev'$. Since $M$ is compact, $\bdd'$ is complete, 
and so is $\bdd$ on $\tilde M$.  
Thus, the $\bdd$-completion of the universal cover 
$\tilde M^o$ of $M^o$ can be identified with $\tilde M$,
and $\ideal{M^o}$ with $\delta \tilde M$. 

If $M^o$ is not $2$-convex, then by Proposition \ref{prop:mconv}
there exists a $3$-simplex $T$ in $\tilde M^o$ 
with sides $F_1, F_2, F_3, F_4$ such that $T \cap \delta \tilde M 
= F_1^o \cap \delta \tilde M \ne \emp$.
Let $x$ be a point of $F_1^o \cap \delta \tilde M$.
Let $N$ be an open ball-neighborhood of $x$ in $\tilde M$
such that $p| N$ is a diffeomorphism onto $p(N)$ 
a neighborhood of $p$ in $N$. 
Then the open surface $p(F_1^o \cap N)$ shows 
that $p(x)$ is a concave point of $\delta M$.  
Since this is a contradiction, $M$ is $2$-convex.
\end{pf}

\begin{rem}
If an affine manifold $M$ is $i$-convex\cO then it 
is $j$-convex for each $j,$ $j \geq i$. 
\end{rem}
\begin{pf}
Let $M$ be $i$-convex with $i \leq n-2$. 
We will now show $M$ is $(i+1)$-convex. The rest 
follows by induction.

Let $T$ be an $(i+2)$-simplex in $\che M$ with sides 
$F_1, \dots F_{i+3}$, and assume $F_2, \dots, F_{i+3}$ 
and $T^o$ do not meet $\ideal{M}$. $T$ can be considered 
the union of $(i+1)$-simplices $S_t$ for 
$t \in [0, 1]$ with common side $G_2$ equal to 
$F_2 \cap F_3$. In fact, we let $S_t$ to be the convex hull of 
$tv_2 + (1-t) v_3$ and $G_2$.
Each $S_t$ has sides $G^t_1, \dots G^t_{i+2}$ 
such that $G^t_2 = G_2$ and $G^t_1$ is a subset of $F_1$ 
and for $0 < t < 1$, the interior of $G^t_1$ is a subset 
of $F_1^o$. Since $G^t_2, \dots, G^t_{i+2}$ are respectively 
subsets of $F_2 \cup \cdots \cup F_{i+3}$, it follows 
that $G^t_2, \dots, G^t_{i+2}$ does not meet $\ideal{M}$.  
Since $M$ is $i$-convex, it follows that $S_t$ are all 
subsets of $\tilde M$ and so is $T$. 
Hence $\tilde M$ is $(i+1)$-convex. 
By induction the rest follows. 
\end{pf} 

\section{Tetrahedral convexity and solid-cylinder convexity}

As we said above, from Sections 3 to 7 we assume that 
{\em our affine $3$-manifold is simply connected}.
We have $M = \tilde M$ in this case. 
The purpose of this section is to show that two-convexity 
implies ``tetrahedral two-convexity", Proposition \ref{prop:conveasy},
and further ``solid-cylinder two-convexity", Theorem \ref{thm:convcyl}.  
This theorem will be used to prove the filling of a drum with 
the bottom side removed in the next section.

A point $x$ of a convex subset $A$ of $\bR^n$ 
is said to be {\em exposed} if there exists a supporting 
affine hyperplane $H$ at $x$ such that $H \cap A = \{x\}$
(see Berger \cite[p. 361]{B}). Compare the following proposition
to Proposition \ref{prop:mconv}: Although there is no 
three-simplex with the interior of a side meeting 
the ideal set $\ideals{M}$ and meeting $\ideals{M}$ only there, 
there could still be a three-simplex meeting 
$\ideals{M}$ at a side. 

\begin{prop}\label{prop:conveasy}
Let $M$ be a two-convex closed affine three-manifold
which is simply connected.
Let $T$ be a convex compact three-simplex
in $\che M$, $F_1$ a side of $T$ such that 
$T \cap \ideals{M} = F_1 \cap \ideals{M} \ne \emp,$ and 
$P$ the affine two-space including 
$\dev(F_1)$. Then  
the convex hull $K$ of $\dev(F_1 \cap \ideals{M})$ 
in $P$ has no exposed point belonging to $\dev(F_1^o)$. 
\end{prop}
\begin{pf}
Suppose not. Let $x'$ be an exposed point of $K$
belonging to $\dev(F_1^o)$, and 
$H$ the supporting line in $P$ such that 
$H \cap K = \{x'\}$. It is clear that $x'$ is 
an element of $\dev(F_1 \cap \ideals{M})$. 
Let $s'$ be the segment $\delta H \cap \dev(F_1)$, 
and $s$ the inverse image in $F_1$ of $s'$ 
under $\dev| F_1$; let $x$ be the point of 
$s$ corresponding to $x'$. 
Then $s - \{x\} \subset M$, 
and $x \in \ideals{M}$. 
A component of $\dev(F_1) - s'$ is 
disjoint from $K$. 
Let $\alpha$ be the closure of the corresponding component of 
$\delta F_1 - s$ disjoint from $\ideals{M} \cap F_1$. 
Since $\alpha$ is a compact arc 
disjoint from $\ideals{M}$, there exists a compact three-ball 
neighborhood $\cal U$ of $\alpha$ in $M$. 
By Proposition \ref{prop:extmap}, we may 
assume without loss of generality that 
$\dev| T \cup \cal U$ is an imbedding onto 
a compact subset of $\bR^3$. 
We may choose $\cal U$ so that $\dev(\cal U)$ is 
a compact neighborhood of $\dev(\alpha)$ 
consisting of points of $\bdd$-distance from $\dev(\alpha)$ 
less than or equal to $\eps$ for some small 
positive constant $\eps$.   

Let $v'_1$ be the vertex of $\dev(T)$ opposite to $\dev(F_1)$, 
and $T'$ a convex three-simplex in $\bR^3$ 
with vertex $v'_1$ satisfying the following condition: 
\begin{itemize}
\item The triangle in the unit tangent bundle 
at $v'_1$ determined by $T'$ matches with 
that determined by $\dev(T)$. 
\item $F_1' \cap \dev(F_1) = \dev(s)$
where $F'_1$ is the side of $T'$ opposite to $v'_1$.
\item $\delta T' - F_1^{\prime o} \subset \dev(T \cup {\cal U})$.
\end{itemize}

Let $T'_t$ for $t \in (0, 1]$ be the convex $3$-simplex 
that is the image of $T'$ under the dilatation 
with center $v'_1$ by the magnification factor $t$.
Let $v_1$ be the vertex of $T$ corresponding to $v'_1$, 
and consider the subset $A$ of $(0, 1]$ 
whose element $t$ is such that 
$M$ includes a convex three-simplex $T_t$ such 
that $\delta T_t$ contains $v_1$ and $\dev| T_t$ is 
an imbedding onto $T'_t$.  
By Remark \ref{rem:magnify}, $A$ is of form 
$(0, \tau)$ or $(0, 1]$. 

Suppose that $A$ is of form $(0, \tau)$.  
By Remark \ref{rem:magnify}, there exists a compact convex three-ball 
$T_{\tau}$ in $\che M$ such that 
$\dev| T_{\tau}$ is an imbedding onto $T'_{\tau}$
and $T_{\tau}^o \subset M$.

The third condition  above states that there exists 
a compact subset $D$ of $T \cup {\cal U}$ such that 
$\dev(D)$ equals $\delta T' - F_1^{\prime o}$. 
There exists a compact subset $D_{\tau}$ of $T_{\tau}$,
such that $\dev| D_{\tau}$ is 
an imbedding onto $\delta T'_{\tau} - F_{\tau}^{\prime o}$
where $F_{\tau}^{\prime o}$ is a side 
of $T'_{\tau}$ opposite to $v'_1$. 
We have $\dev(D_{\tau}) \subset \dev(D)$. 
Since $\dev| T \cup {\cal U} \cup T_{\tau}$ is an imbedding 
by Proposition \ref{prop:extmap} (up to choosing $\cal U$ carefully), 
we obtain $D_{\tau} \subset D \subset M$. 
Hence, $T_{\tau} \cap \tilde M_\infty$ is included in 
the interior of the side of $T_{\tau}$ corresponding to 
$F'_{\tau}$. 
By definition of two-convexity, it follows that 
$T_\tau \subset M$ and $\tau \in A$, a contradiction. 

Suppose now $A = (0, 1]$.  
Then $M$ includes a three-simplex $T_1$
such that $\dev| T_1$ is an imbedding onto $T'$.
By Proposition \ref{prop:extmap}, $\dev| T \cup T_1$ is 
an imbedding onto $\dev(T) \cup \dev(T_1)$. 
As $\dev(s)$ is a subset of $T'$, we obtain $s \subset T_1 \subset M$; 
however, this contradicts $x \in \ideals{M}$.  
\end{pf}

\begin{rem} We will meet this type of changing 
the direction or tilting of a side again in Section 8.
\end{rem} 

Let $(\dev, h)$ be a development pair of $M$. 
We have affine coordinate functions 
$x_1, x_2, x_3: \bR^3 \ra \bR$. 
For convenience, we let $x_i$ denote the map 
$x_i \circ \dev$ defined on $M$ for each $i$, $i=1, 2, 3$.  
Let $M_t$ denote $x_1^{-1}(t) \cap M$, which 
is a closed subset of $M$.  

\begin{lem}\label{mt:open} 
For each $t,$ $t\in \bR$, each component of 
$M_t$ is a properly imbedded open surface or an empty set. 
In particular, the union of a collection of 
disjoint simple closed curves in $M_t$ forms the boundary 
of at most one subsurface of $M_t$.  
\end{lem} 
\begin{pf} Since $M$ is without boundary, 
$M_t$ is also without boundary. 
Let $U$ be a component of $M_t$. 
Since $\dev|U$ is an immersion onto 
an open subset of  
an affine hyperplane of $\bR^3$, 
$U$ is not compact. 
Hence, $U$ is an open surface.
The conclusion follows.   
\end{pf}

We now define solid-cylinder maps: 
Given an embedded compact surface $\Omega$ 
in the plane $\bR^2$ with smooth boundary and   
the interval $[a, b]$ of $\bR$ for $a, b \in \bR$, $a < b$,  
think of $\Omega \times [a, b]$ as sitting in 
the Euclidean space $\bR^2 \times \bR$.
A continuous map 
$f: \Omega \times [a, b] \ra \che M$ is called 
a {\em solid-cylinder map\/} if $f$ satisfies the following 
conditions:
\begin{enumerate}
\item[(i)] $f| \Omega \times [a, b)$ is an injective
map into $M$. 
\item[(ii)] $\dev \circ f| \Omega \times [a, b]$ is 
a $C^1$-immersion. 
\item [(iii)] $x_1 \circ \dev \circ f| 
\Omega \times \{t\}$ is a constant map with value $t$ for 
each $t$, $t \in [a, b]$.  
\item [(iv)] $f| \delta \Omega \times [a, b]$ 
is an embedding into $M$.  
\end{enumerate}
In particular, $f| \Omega \times \{t\}$ is an embedding into 
$M_t$ for each $t$, $t \in [a, b)$. 
We will often call $M_t$ and 
$\Omega \times \{t\}$ {\em $t$-level sets}. Thus,  
$f$ preserves levels. 

A priori, $f(\Omega \times \{b\})$ may intersect $\ideals{M}$.
We will show that this does not happen while $M$ is 
$2$-convex.

\begin{lem}\label{cyl:simplex} 
Let $f: \Omega \times [a, b] \ra \che M$ be a solid-cylinder map. 
For every point $x$ of $\Omega^o \times \{b\},$ 
there exists a compact neighborhood 
$V$ such that $f| V$ is an imbedding onto 
a three-simplex $T$ in $\che M$ with 
a side $F_1$ such that $T - F_1 \subset M,$ 
$f(x) \in F_1^o,$ and $x_1 \circ \dev(F_1) =\{b\}$.  
\end{lem}
\begin{pf} 
The conditions (ii) and (iii) show that 
$\dev \circ f$ is an immersion of $\Omega \times [a, b]$ 
into the half space $H$ of $\bR^3$ given by $x_1 \leq b$ while
$\dev \circ f(\Omega \times \{b\})$ is the subset of 
the boundary $\delta H$, which is an affine hyperplane 
of $\bR^3$.   
Thus, $\Omega \times [a, b]$ includes 
a compact neighborhood $V$ of $x$ so that 
$\dev \circ f| V$ is an imbedding onto 
a three-simplex $T'$ with a side $F'_1 \subset \delta H$.  
Since $T'$ is a neighborhood of $\dev \circ f(x)$ in $H$, 
we have $\dev \circ f(x) \in F_1^{\prime o}$. 
Let $T = f(V)$. Then $\dev| T$ is an imbedding 
onto $T'$, and the lemma follows.  
\end{pf}

\begin{thm}\label{thm:convcyl}
Suppose that $M$ is a $2$-convex affine $3$-manifold 
which is simply connected.
If $f$ is a solid-cylinder map\cO then $f$ is a map into $M$. 
\end{thm}
\begin{pf}
By condition (i), $f| \Omega \times \{b\}$ 
is a map into $\che M$. 
Define
$K = (f| \Omega \times \{b\})^{-1}(\ideals{M})$. 
By condition (iv), 
$K$ is a compact subset of $\Omega^o \times \{b\}$.
Suppose that $K$ is not empty.
Let $r$ be the function on $\bR^3$ defined by 
$r = x_2^2 + x_3^2$. 
Let $r' = r \circ \dev \circ f$, 
and choose a maximum point $x$ of 
$r'$ restricted to $K$. By above lemma, 
$\che M$ includes a convex three-simplex $T$ 
with a side $F_1$ such that $f(x) \in F_1^o$ 
and $T - F_1 \subset M$.

Since $F_1 \subset f(\Omega \times \{b\})$, 
we have $\ideals{M} \cap F_1 \subset f(K)$; hence,  
$\ideals{M} \cap F_1 = f(K) \cap F_1$. 
Then $\dev\circ f(x)$ is an exposed point of 
the convex hull of 
$\dev(F_1 \cap f(K)) \supset \dev(F_1 \cap \ideals{M})$
since $x$ is a maximum point of $r'$ in $K$. 
This contradicts Proposition \ref{prop:conveasy}. 
Hence, $K$ is empty. 
\end{pf}

\section{Filling in a drum with the bottom side removed}

The purpose of this section is to 
prove the strongest form of two convexity
that we can always fill in the drum with a bottom side removed.
Let $\delta \Omega$ be the boundary of a compact 
surface $\Omega$, which is 
a union of disjoint simple closed curves.  

\begin{thm}\label{thm:fill}
Let $M$ be a $2$-convex affine $3$-manifold which 
is simply connected.
Let $f: \delta \Omega \times [a, b] \ra M$ 
be a differentiable map satisfying the following conditions {\/\rm :} 
\begin{enumerate} 
\item[(i)] $f| \delta \Omega \times \{t\}$ is 
an embedding into a union of disjoint 
simple closed curves in $M_t$ for each $t$.  
\item[(ii)] $f| \delta \Omega \times [a, b]$ is transversal 
to each level set.
\item[(iii)] $f| \delta \Omega \times \{a\}$ 
extends to an embedding $f': \Omega \times \{a\} \ra M_b$. 
\end{enumerate}
Then $f$ and $f'$ extend to a common solid-cylinder embedding 
${\cal F}$ from $\Omega \times [a, b]$ to $M$.  
\end{thm}

The pair of maps $f$ and $f'$ satisfying (i), (ii), and (iii)
or the union of their images are said to 
be a {\em drum with the top side removed}.  

Actually the way we will use this theorem is that $a$ in (iii) will 
be replaced by $b$. In this case the pair $f$ and $f'$ is said
to be a {\em drum with the bottom side removed}. 
Hence, we use this theorem with the direction 
of ``filling'' reversed.

By the above theorem, there exists a surface $\Omega' \subset M_a$ 
such that $\Omega'$ is diffeomorphic to $\Omega$, 
and $\Omega' \cup \hbox{Im} f \cup \hbox{Im} f'$ 
forms the boundary of ${\cal F}(\Omega \times [a, b])$. 

{\em It is sufficient to prove 
the above theorem when $a = 0$ and $b > 0$}.

The plan to prove this theorem is as follows:
\begin{itemize}
\item We define a flow, defined on $M$ only partially
and not complete, which induces the movement of 
$f(\delta \Omega, t)$ and sends level sets 
to level sets as $t$ increases from $0$ to $b$. 
\item This flow will be used to define a level-preserving
embedding from $\Omega \times [0, \tau)$ to $M$ 
extending $f'$ and $f$ where $\tau$, $0 < \tau < b$, denotes 
the level at which we cannot extend the flow. 
\item Using flow estimates with respect to $\bdd$ 
on $M$, we show that the above map extends 
to a continuous map $F: \Omega \times [0, \tau] \ra \che M$.
\item By showing $\dev\circ F$ is $C^1$, 
we prove that $F$ is a solid-cylinder map.
\item Using Theorem \ref{thm:convcyl}, we show $F$ is 
a map into $M$ and conclude that $\tau = b$. 
\end{itemize}

Recall that $\mu$ is the Euclidean metric on $M$ induced 
by $\dev$, and $M$ has standard coordinates $x_1, x_2$, and $x_3$. 
During the proof of the above theorem, we will state everything
with respect to this Riemannian metric and the coordinates.  

{\em We first need a vector field for our flow.\/} 
The image of $f$ is the union of disjoint 
embedded annuli in $M$. 
Let us call it $\Sigma$. For each point $(x, s)$ of 
$\delta \Omega \times \{s\}$, the derivative 
$\dot{f}(x, s)$ of $f(x, s)$ 
with respect to $s$ yields a vector ${\bf w}$ at $f(x, s)$ 
tangent to $\Sigma$. This defines 
a smooth tangent vector field ${\bf w}$ to $\Sigma$.  
Let ${\bf w}^P$ be the component vector field of ${\bf w}$ 
orthogonal to the $x_1$-direction. Since ${\bf w}^P$ is a smooth 
vector field, we can find a smooth vector field 
${\bf v}^P$ extending ${\bf w}^P$
defined on a neighborhood $\cal U$ of $\Sigma$ and  
orthogonal to the $x_1$-direction. 

Let $\phi$ be 
a smooth function which has compact support in $\cal U$, 
has its range of values lying between $0$ 
and $1$, and is equal to $1$ on $\Sigma$.  
Then $\phi{\bf v}^P$ is regarded as a smooth vector field defined on 
$M$ orthogonal to the $x_1$-direction.    
Since $\phi$ is not zero only in a compact subset, 
the derivatives of $\phi$ of each order are uniformly bounded 
in $M$ by a constant depending only on the order.
Similarly, those of 
${\bf v}^P$ of each order are uniformly bounded 
in the support of $\phi$ by a constant depending 
only on the order. Hence, those of  
$\phi{\bf v}^P$ of each order are uniformly bounded above in $M$ 
by constants depending only on the order.  

Let $e_1$ denote the vector field on $M$ 
in the $x_1$-direction with norm $1$, 
and $\tilde W$ the vector field $\phi{\bf v}^P + e_1 $. 
Then $\tilde W$ extends ${\bf w}$ and is a smooth 
vector field such that the norms of 
itself and its derivatives of order 
less than or equal to two are bounded 
above by a uniform constant $c_{\tilde W}$. 

Let ${\cal O}_0$ denote $f'(\Omega \times \{0\})$,  
and $\Phi_t$ the flow generated by 
$\tilde W$ where $\Phi_0$ equals the identity map 
of $M$. Let us denote by $A$ the 
set of elements $t$ of $[0, b]$ 
such that $\Phi_t(x)$ is defined for every $x \in {\cal O}_0$. 
Since ${\cal O}_0$ is compact, $A$ is open. 

If $\Phi_t(x)$ is defined, then 
$\Phi_{t'}(x)$ is defined for $0 \leq t' \leq t$. 
Thus, $A$ is a connected subset. 

Now, we prove that $A$ is closed.   
Suppose not. Then $A = [0, \tau)$ for some $\tau \in [0, b]$. 
Let $t \in [0, \tau)$. 
Since the $x_1$-component of $\tilde W$ equals $1$ identically, 
and $\tilde W$ extends ${\bf w}$, it follows that 
$\Phi_t({\cal O}_0)$ is a compact surface in $M_t$ with 
boundary equal to the union of 
disjoint smoothly imbedded closed curves 
$\Phi_t(\delta {\cal O}_0)$, which equals 
$f(\delta \Omega \times \{t\})$. 
Let ${\cal O}_t$ denote this surface and
${\cal O}_{[0, \eps]}$ the union of 
${\cal O}_t$ for $t \in [0, \eps]$, 
where $\eps$ is a small positive number 
so that $\eps < \tau$. 

We define an injective immersion 
$F_\tau: \Omega \times [0, \tau) \ra M$ 
by $F_\tau(x, s) = \Phi_s(f'(x, 0))$ 
for $x \in \Omega$ and $s \in [0, \tau)$.  
Then $F_\tau| \Omega \times [0, t]$ for each $t$, 
$0 < t < \tau$, is a solid-cylinder embedding extending $f'$
and $f| \delta \Omega \times [0, t]$ simultaneously.    

From the definition of $F_\tau$, we obtain
for $t, s$, $0 < t + s < \tau$, $t, s > 0$,
and every $x \in \Omega$,   
\[ F_\tau \circ T_s(x, t) = \Phi_s \circ F_\tau(x, t)\]   
where $T_s$ is a translation by a vector $(0, s)$ in $\bR^3$.   
While $t \in [0, \eps]$ for any positive 
number $\eps$ such that $\eps < \tau$, 
$s$ may take any value in $[0, \tau -\eps)$
for the above equation to make sense. 
Thus, for $s \in [0, \tau-\eps)$, 
the following diagram of diffeomorphisms is commutative: 
\begin{equation}\label{comm:equ}
\begin{CD}
\Omega \times [0, \eps] @>{F_\tau}>> {\cal O}_{[0, \eps]}\\
@V{T_s}VV				@V{\Phi_s}VV \\
\Omega \times [s, s+\eps] @>{F_\tau}>> {\cal O}_{[s, s+\eps]}.
\end{CD}
\end{equation}  

Since the differential of $\tilde W$ is bounded with 
respect to $\mu$, Lemma 3 \cite[p. 63]{Abr} implies that   
\[ 
\bdd(\Phi_s(x), \Phi_s(y)) \leq e^{K|s|}\bdd(x, y)
\]
for every $s \in [0, \tau)$, $x, y \in {\cal O}_0$, 
and a constant $K$ independent of $s$.  
Then we have 
\[\bdd(F_\tau(x, s), F_\tau(y, s)) 
\leq e^{K|\tau|}\bdd(f'(x, 0), f'(y, 0)) \] 
for every $x, y \in \Omega$ and $s$, $0 < s <\tau $.  
If the second $\bdd$ below denotes 
the standard Euclidean distance metric 
on $\Omega \times [0, b]$, we have 
for every $x, y \in \Omega$ and 
an independent positive constant $R'$       
\begin{equation}\label{comp:equ}
\bdd(F_\tau(x, s), F_\tau(y, s)) 
\leq R'\bdd((x, s), (y, s)).
\end{equation} 
It is easy to show that the following inequality holds 
if $(y, s)$ is replaced by $(y, s')$ for $s' \in [0, \tau)$
and $R'$ is replaced by another independent positive constant $R$: 
\begin{equation}\label{comp:equ2}
\bdd(F_\tau(x, s), F_\tau(y, s')) 
\leq R\bdd((x, s), (y, s')).
\end{equation} 
Since $\che M$ is complete, $F_\tau$ extends to 
a continuous map $F: \Omega \times [0, \tau] \ra \che M$. 
Therefore, 
\begin{equation}\label{comm2:equ}
\Phi_s(x) = F \circ T_s \circ F_{\tau}^{-1}(x)
\end{equation}
is well-defined with range in $\che M$ for 
$x\in {\cal O}_{[0, \eps]}, s \in [0, \tau-\eps]$.    
Moreover, since $T_s \ra T_{\tau - \eps}$ pointwise as 
$s \ra \tau - \eps$, it follows that  
$\Phi_s| {\cal O}_{[0, \eps]}$ 
converges to $\Phi_{\tau - \eps}$ pointwise 
as $s \ra \tau-\eps$.  

\begin{figure}[h] 
\centerline{\epsfysize=6.75cm
\epsfbox{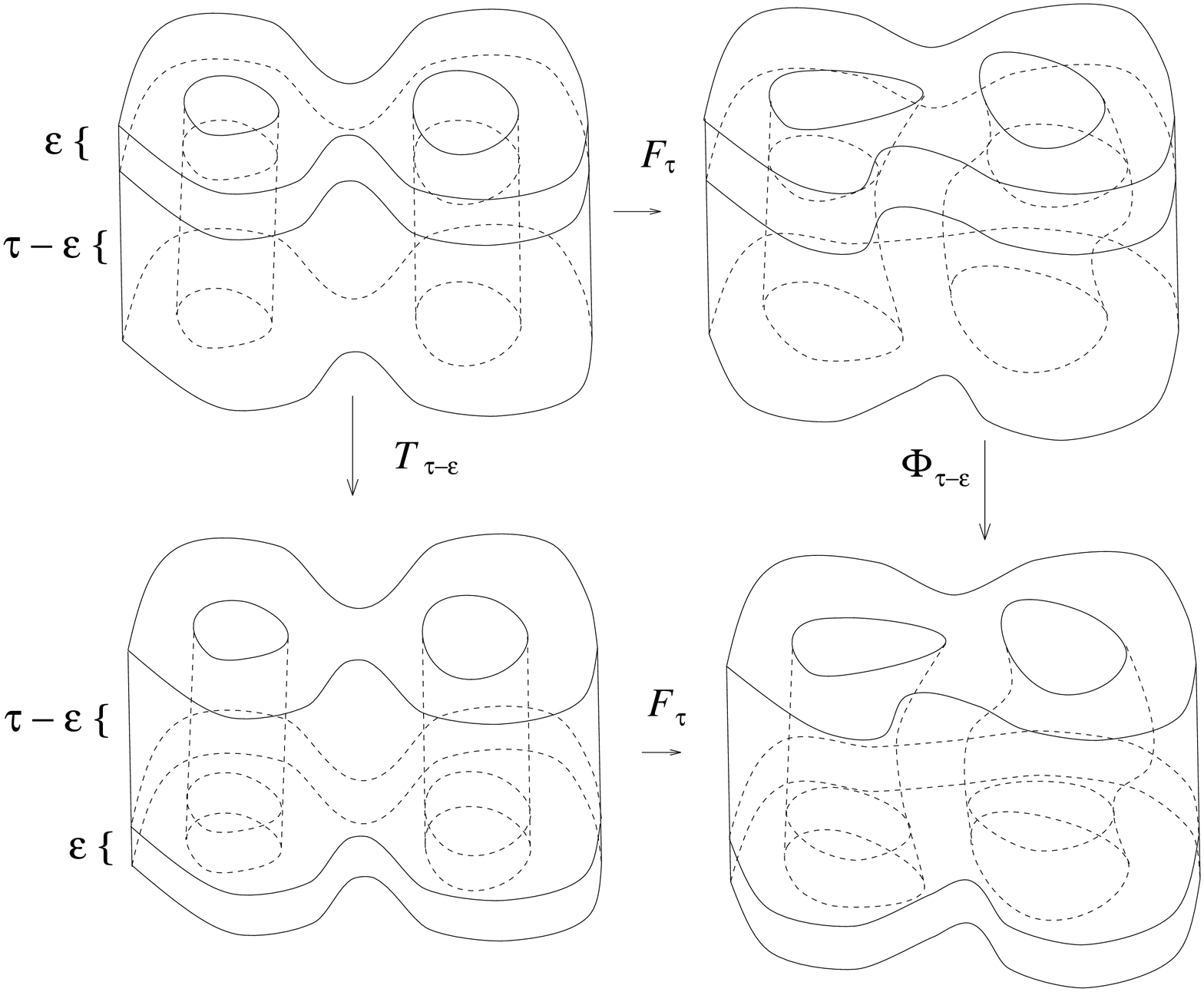}}
\caption{\label{fig:uf4} The flow picture} 
\end{figure} 

{\em We claim that $\dev\circ F$ is a $C^1$-immersion.\/} 
Let $D\Phi_s(x)$ for $x$ and $\Phi_s(x)$ 
in the image of $F_\tau$ denote 
the $3\times 3$-matrix of the differential 
$T_x(M) \ra T_{\Phi_s(x)}(M)$
of $\Phi_s$ at $x$, and $D\tilde W(x)$ 
that of the differential $T_x(M) \ra T_x(M)$
of $\tilde W$ at $x$. 
(Recall the global coordinates $x_1, x_2, x_3$ 
introduced earlier in this section.)
Then 
\begin{equation}\label{floweqn}
\frac{d}{d s} D \Phi_s(x) = D\tilde W(\Phi_s(x))D\Phi_s(x) 
\end{equation}
for $x \in {\cal O}_{[0, \eps]}$ 
and $0 \leq s < \tau - \eps$ (see Lemma 4 of \cite[p. 64]{Abr}).  
Since the norms of $D\tilde W$ and $D^2\tilde W$ are
uniformly bounded above by $C_{\tilde W}$ and $\Phi_0$ is 
the identity map, then  
the proof of Lemma 4 of \cite[p. 65]{Abr} 
shows that the family of $3 \times 3$-matrix-valued functions
$D\Phi_s|{\cal O}_{[0, \eps]}$ for $s \in [0, \tau-\eps)$ 
and that of $9 \times 3$-matrix-valued functions  
$D^2\Phi_s|{\cal O}_{[0, \eps]}$ for $s \in [0, \tau-\eps)$ 
are respectively uniformly bounded.  
Hence, the family of $3 \times 3$-matrix-valued-functions 
$D\Phi_s| {\cal O}_{[0, \eps]}$ 
for $s \in [0, \tau-\eps)$ forms a bounded, equicontinuous 
family of functions. (We benefited from 
a conversation with D.-H. Chae who works in Navier-Stokes
equations).  

Since $\dev$ preserves coordinates $x_1, x_2,$ and $x_3$, 
$D\dev \circ \Phi_s| {\cal O}_{[0, \eps]}$ equals 
$D\Phi_s| {\cal O}_{[0, \eps]}$ as matrix-valued functions,  
and $D^2 \dev \circ \Phi_s| {\cal O}_{[0, \eps]}$ 
equals $D^2 \Phi_s| {\cal O}_{[0, \eps]}$. 
By the Ascoli-Arzel\`a theorem (see \cite[p. 85]{Yosh}), 
choose a subsequence 
$\Phi_{s_j}$ such that $\{s_j\} \ra \tau -\eps$ and 
$D\dev\circ\Phi_{s_j}| {\cal O}_{[0, \eps]}$ 
converges uniformly to a continuous 
$3\times 3$-matrix-valued function $\Upsilon$  
defined on ${\cal O}_{[0, \eps]}$.     

Since
$\dev \circ \Phi_{s_j}| {\cal O}_{[0, \eps]}$ converges 
pointwise to $\dev \circ \Phi_{\tau-\eps}| {\cal O}_{[0, \eps]}$,
Theorem 7.17 \cite[p. 152]{Rudin} 
(a standard fact) implies that 
$\dev\circ\Phi_{\tau - \eps}| {\cal O}_{[0, \eps]}$ is continuously 
differentiable with 
$D\dev\circ\Phi_{\tau-\eps}| {\cal O}_{[0, \eps]} = \Upsilon$. 

We obtain the following equation from equation \ref{floweqn}
by taking determinant of the both sides
\[
\frac{d}{d s} \det D \Phi_s(x) = \det D\tilde W(\Phi_s(x))
\det D\Phi_s(x)\] 
If the initial value $f(0)$ of an ordinary differential 
equation $\frac{d}{ds} f = g f$ for real valued solution $f$ 
is positive, then the solution is positive always,  
as it can be seen from the solution 
\[ f(s) = \exp[\int^s_0 g ds] f(0) \]
Since $\det D\Phi_0 = 1$, $\Upsilon$ has values in nonsingular 
$3 \times 3$-matrices.  

By equation \ref{comm2:equ}, the following 
diagram is commutative: 
\begin{equation}
\begin{CD}
\Omega \times [0, \eps] @>{F_\tau}>> {\cal O}_{[0, \eps]}\\
@V{T_{\tau-\eps}}VV		@V{\Phi_{\tau-\eps}}VV \\
\Omega \times [\tau-\eps, \tau] @>{F}>> \che M. 
\end{CD}
\end{equation}  
Thus,  
\[\dev \circ F| \Omega \times [\tau-\eps, \tau] 
= (\dev \circ \Phi_{\tau-\eps}) \circ 
F_{\tau} \circ T_{\tau -\eps}^{-1}
|\Omega \times [\tau-\eps, \tau]. \]  
Since $D \dev \circ \Phi_{\tau - \eps} = \Upsilon$, 
it follows that 
$\dev \circ F| \Omega \times [\tau-\eps, \tau]$ 
is of class $C^1$, and so is $\dev \circ F$ 
on $\Omega \times [0, \tau]$.    
Moreover, since $\Upsilon$ has values in
nonsingular matrices, 
$\dev\circ F|\Omega \times [\tau-\eps, \tau]$ is an immersion.  
Since $F$ extends $F_\tau$, 
and $F_\tau$ is 
an immersion into $M$, 
$\dev \circ F$ is an immersion, 
and $F: \Omega \times [0, b] \ra \che M$ is a solid-cylinder map.

By Theorem \ref{thm:convcyl}, $F$ is a map into $M$.
Hence, our vector field $\tilde W$ is defined on  
a compact neighborhood of the image of $F$ in $M$. 
Recall that 
$F_\tau: \Omega \times [0, \tau) \ra M$ is given  
by $F_\tau(x, s) = \Phi_s(f'(x, 0))$.  
One defines $\Phi_\tau(f'(x, 0))$ by 
letting it equal to $F(x, \tau)$. 
Then $\Phi_\tau$ is clearly a flow generated by 
$\tilde W$, and $\tau$ belongs to our set $A$. 
This is a contradiction. Therefore, 
$A$ is closed and $A = [0, b]$. 

Define ${\cal F}: \Omega \times [0, b] \ra M$ 
by ${\cal F}(x, s) = \Phi_s(f'(x, 0))$ for $s \in [0, b]$.  
Then $\cal F$ satisfies the conclusion of Theorem 
\ref{thm:fill}.  

\section{Filling a Disk with One Critical Point and an Annulus with 
Two Critical Points}

Let $\Omega$ be a compact surface.  
Suppose that $F: \Omega \ra M$ is a smooth embedding 
with the property that $F(\delta \Omega)$ 
is a subset of $M_t$ for some $t \in \bR$, 
$x_1 \circ F$ is a Morse function
with $k$ critical values, and 
$t$ not a critical value.  
Then $F$ is said to be an {\em amenable embedding\/} of level $t$
with $k$ critical values.  
We say that the image of $F$ is an {\em amenable surface\/}, 
and if $x_1\circ F$ increases in the inner direction on 
$\delta \Omega$, then $F$ is a {\em positive\/} amenable embedding; 
otherwise, $F$ is a {\em negative\/} amenable embedding.

The purpose of this section is to fill in an amenable 
disk with one critical point and an amenable annulus with 
two critical points. From the following Proposition \ref{prop:Euler}
these have the smallest number of critical points 
in their respective homeomorphism classes. 
The filling of the disk follows from the 
drum filling that we obtained above.
The proof of the filling of annulus is as follows: 
The idea is to cut off the surface into two parts
at the level $k$ just above the level of the saddle point, 
where the higher one has the unique maximum point 
and the lower one has the saddle point. The higher  
surface is homeomorphic to a disk and we can use the filling of 
the disk. The lower surface is homeomorphic to 
a pair of pants with two boundary components at 
a lower level and one at the level just above the level of 
the saddle. The lower surface is cut and pasted near the saddle into 
the union of two annuli by finding two disks in 
a small neighborhood of the saddle in $M$. 
Since the higher surface is filled,  
$M_k$ include a disk $D$ with which the higher surface bound 
a $3$-ball. $D$ is then altered according to the surgery so that 
$D$ together with the two annuli form 
a drum with the bottom side removed.  The top side may 
consist of the union of two disjoint disks or an annulus.  
We then fill in the drum and take the union with the ball 
filling the higher disk from which we add or subtract to fit 
in the original surface. Depending on whether the top 
side is the union of two disjoint disks or an annulus, 
the result is a $3$-ball or a solid torus.
(The both processes are described in Figure \ref{fig:movie0}.)

The level sets of $x_1 \circ F$ form a singular 
foliation on $\Omega$ with components of $\delta \Omega$ 
as leaves. Let us call the foliation 
the {\em level-set foliation} on $\Omega$ for $F$.  
For such foliations indices of singularities are well-defined. 
(See Casson-Bleiler \cite[p. 71]{Casson}.) 

\begin{prop}\label{prop:Euler}
The sum of indices of all singularities of 
the level-set foliation on $\Omega$ for $F$ equals the Euler 
characteristic of $\Omega.$ 
\end{prop}  
\begin{pf} This follows since the boundary components 
are leaves of the foliation. See \cite[p. 72]{Casson}. 
\end{pf}

Each point of $M$ has a neighborhood $U$ such that 
$\dev| U$ is a chart. Hence, $U$ includes a convex 
open ball neighborhood of $x$. 

\begin{prop}\label{prop:diskfill1}
Let $\disk$ be a compact disk with smooth boundary, and
$F: \disk \ra M$ an amenable embedding 
of level $t$ with one critical point.  
Then there exists an embedding $G: \disk \ra M_t$ 
homotopic to $F$ relative to $\delta \disk$
with $F(\disk^o) \cap G(\disk^o) = \emp$. 
Moreover\cO the union of the imbedded disks 
$F(\disk)$ and $G(\disk)$ 
is the boundary of a compact three-ball in $M$.  
\end{prop} 
\begin{pf}
Assume without loss of generality that $F$ is positive. 
Let $\Sigma = F(\disk)$, and let $f_1$ 
denote $x_1| \Sigma$.   
Let $k$ be the unique critical value of $f_1$
and $z$ the critical point. Then $z$ is 
the maximum point of $f_1$.  
Choose a convex open ball neighborhood $B$ of $z$. 
Then there exists a real number $k'$ near $k$, $k' < k$, so that 
the level set $f^{-1}(k')$ at $k'$ 
is a simple closed curve in $\Sigma \cap B$, say $\alpha$, 
which equals the boundary of a disk $\Sigma_\alpha$ in $\Sigma^o$. 
Since $M_{k'} \cap B$ is a convex open disk, 
$\alpha$ is the boundary of a disk $D^2_{k'}$ in $M_{k'}$.   
Now, the union of $D^2_{k'}$ and $\Sigma_\alpha$
is the boundary of a three-ball in $B$. 
Let us call the ball $B_1$. 

By Theorem \ref{thm:fill}, it follows that 
there exists a disk $D^2_t$ in the level $t$ 
so that the union of $D^2_{k'}$, $D^2_t$, and 
the embedded annulus $f_1^{-1}([k', t])$ 
is the boundary of a three-ball. Let us call it $B_2$. 
Then $B_1 \cup B_2$ is the three-ball with boundary 
$\Sigma \cup D^2_t$. The proposition follows easily now.  
\end{pf}    

\begin{figure}[h] 
\centerline{\epsfysize=4.25cm
\epsfbox{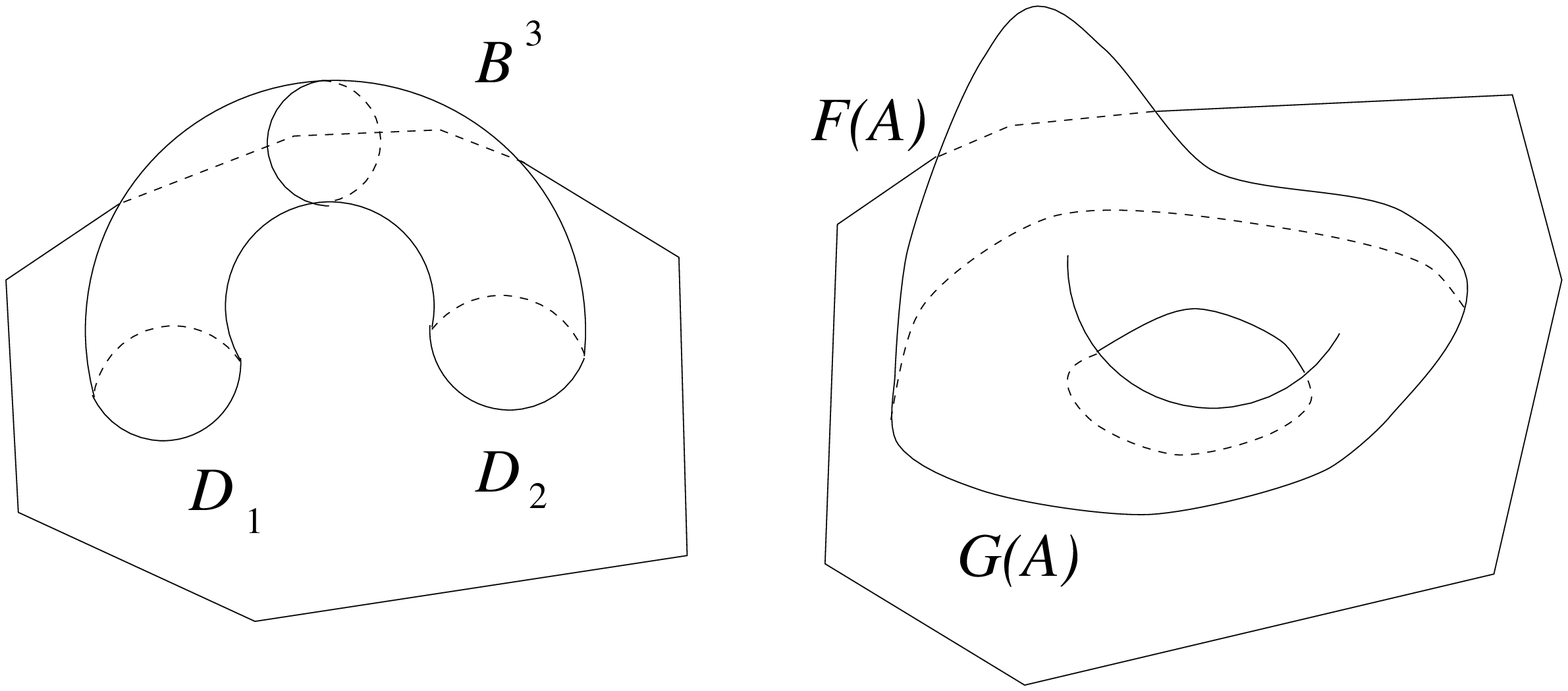}}
\caption{\label{fig:uf5} Examples of two types of annuli and 
fillings by a $3$-ball and a solid torus respectively.}
\end{figure} 

\begin{prop}\label{prop:annfill}
Let $A$ be a compact annulus\cO and
$F: A \ra M$ a positive amenable embedding of 
level $t$ with two critical points. 
Then one of the following holds{\/\rm :\/} 
\begin{itemize}
\item[($\alpha$)] $F(\delta A)$ is the 
boundary of the union of two disjoint disks 
$D_1$ and $D_2$ in $M_t$ 
so that $D_1 \cup D_2 \cup F(A)$ is the boundary 
of a compact subset of $M$ homeomorphic to a three-ball. 
\item[($\beta$)] There exists an embedding $G: A \ra M_t$ 
homotopic to $F$ relative to $\delta A$,    
so that $G(A) \cup F(A)$ equals the boundary 
of a compact subset in $M$, homeomorphic to a solid torus. 
\end{itemize} 
\end{prop}  
\begin{pf}
The critical points of $A$ consist of one index $1$ critical
point and one index $-1$ critical point since $\chi(A) = 0$ 
by Proposition \ref{prop:Euler}. 

Let $\Sigma$ be the image $F(A)$, and $f_1 = x_1| \Sigma$.   
Let $l$ and $k$, $l > k$, be 
critical values, and $w$ and $z$ the corresponding 
critical points respectively. Then $w$ is a local maximum point and 
$z$ a saddle point. Let $\eps$ be 
a small positive number so that $k + \eps < l$. 
($\eps$ will be chosen more precisely below.)

Let $\Sigma_0 = f_1^{-1}([k + \eps, l])$.  
Then the higher surface $\Sigma_0$ is a positive amenable disk 
with one critical point. By Proposition \ref{prop:diskfill1},
$M_{k+\eps}$ includes a disk $D_{k+\eps}$,  
the union of which with $\Sigma_0$ is the 
boundary of a compact three-ball in $M$. 
Let $B_1$ denote the three-ball.

Obviously, $t < k < l$. 
Let $\Sigma_1 = f_1^{-1}([t, k + \eps])$, 
which is homeomorphic to a sphere removed 
with three disjoint open disks, i.e., a pair of pants. 
Let $\delta_1 \Sigma_1$ denote  
the unique boundary components of $\Sigma_1$ 
at the level $k + \eps$, 
and $\delta_2 \Sigma_1$ the union of 
the two boundary components in the level $t$. 
Then $D_{k+\eps}$ satisfies 
$\delta D_{k+\eps} = \delta \Sigma_0 = \delta_1 \Sigma_1$. 

Near the saddle point of $\Sigma$, we will do 
the smooth surgery for the lower 
surface $\Sigma_1$ as indicated by Figure \ref{fig:surg}. 

\begin{figure}[h]  
\centerline{\epsfysize=5.25cm
\epsfbox{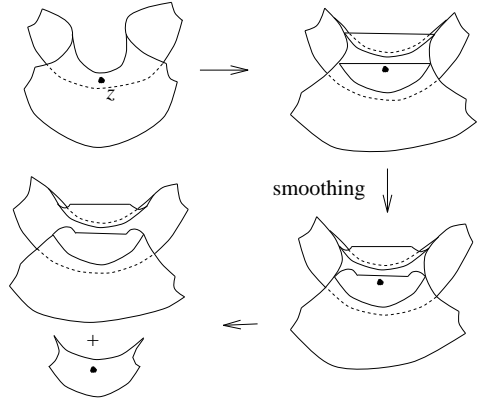}}
\caption{\label{fig:surg} The local smooth surgery process} 
\end{figure} 

{\em We will make a chart near $z$ so that $\Sigma$ is 
mapped to a canonically shaped saddle}.
(i) Let $B$ be an open ball neighborhood of $z$ in $M$. 
Choose a coordinate system $u^1, u^2$ in a neighborhood 
$V$ of $z$ in $\Sigma \cap B$ so that 
$u^1(z) = u^2(z) = 0$ 
and $f_1 = k + (u^1)^2 - (u^2)^2$ hold in $V$
(see \cite{Mil}). 
Extend $u^1$ and $u^2$ to smooth functions on $B$, and  
let $F: B \ra \bR^3$ be given by 
$F(x) = (x_1(x), u^1(x), u^2(x))$. 
Then $D_{z}F$ is nonsingular. 
Hence, there exists a neighborhood 
$U$ of $z$ in $B$ so that 
$F| U$ is a diffeomorphism onto a neighborhood of 
$(k, 0, 0)$ in $\bR^3$. We may assume without 
loss of generality that $\Sigma \cap U$ is 
a connected subsurface.   

For a positive number $\delta$, we let
\[
B_{\delta}(k, 0, 0) = \{(x_1, u^1, u^2) \in \bR^3|
|x_1 - k| < \delta, (u^1)^2 < \delta, (u^2)^2 < \delta \}. 
\]
Let $\Sigma' = F(\Sigma \cap U)$. 
We choose $\delta'$, $\delta' > 0$,  
so that $B_{\delta'}(k, 0, 0) \subset F(U)$
and $B_{\delta'}(k, 0, 0) \cap \Sigma'$ 
equals 
\[ \{(x_1, u^1, u^2) \in B_{\delta'}(k, 0, 0)| 
x_1 = k + (u^1)^2 - (u^2)^2 \}.\]
That is, $\Sigma'' = B_{\delta'}(k, 0, 0) \cap \Sigma'$ 
is realized as a graph of a function. 

(ii) Let $\cal F$ denote the foliation on $M$ 
whose leaves are components of level sets of $x_1$. 
{\em We choose imbeddings of two $\cal F$-transverse disks in 
$B_{\delta'}(k, 0, 0)$ that does a surgery on $\Sigma_1$, 
and obtain a union of disjoint
$\cal F$-transverse annuli with smooth boundary.\/} 
We may choose our cutoff 
number $\eps$ to be less than $3\delta'$. 
Let $\delta = 3\eps$ and define domains
\begin{eqnarray}
D &=&\{ (x_1, u^1)| k + (u^1)^2 - \delta < x_1 < k+\delta, 
(u^1)^2 < \delta\}\nonumber\\
D^u &=&\{(x_1, u^1) \in D| x_1 \geq k + (u^1)^2 - \delta/2\}\\
D^d &=&\{(x_1, u^1) \in D| x_1 \leq k + (u^1)^2 - \delta/2\}.
\nonumber 
\end{eqnarray}

\begin{figure}[h] 
\centerline{\epsfysize=4.5cm
\epsfbox{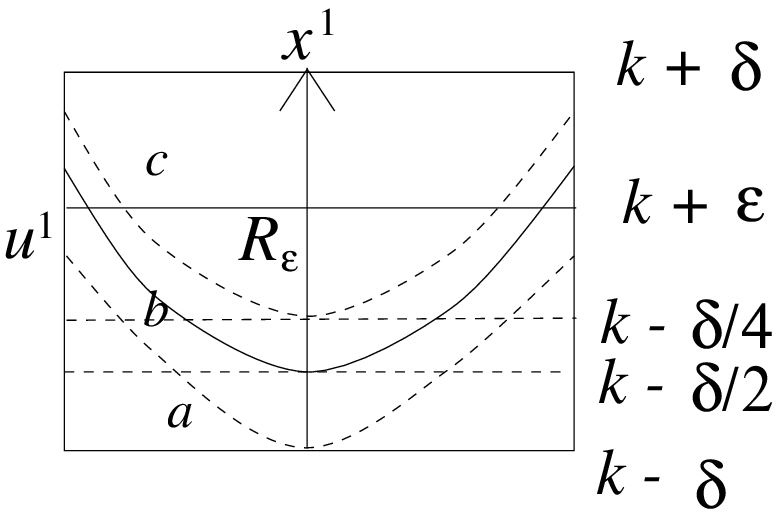}}
\caption{\label{fig:uf7} $D, D^u, D^d, N(D^d),$ and $R^\eps$.
$a$ is given by $x_1 = (u^1)^2 + k -\delta$, 
$b$ by $x_1 = (u^1)^2 + k - \delta/2$, and 
$c$ by $x_1 = (u^1)^2 + k - \delta/4$.
} 
\end{figure} 

We define a continuous map $g: D \ra \bR$ by 
\begin{equation}
g(x_1, u^1) = \left\{\begin{array}{cc}
	     \sqrt{\delta/2} \quad, & (x_1, u^1) \in D^u \nonumber\\ 
		&	\nonumber \\
	     \sqrt{k + (u^1)^2 - x_1} \quad, & (x_1, u^1) \in D^d.
	     \end{array}\right.
\end{equation}
We smooth this function. Let $N(D^d)$ be the 
set of points $(x_1, u^1)$ of $D$ with $x_1 < k + (u^1)^2 - \delta/4$
and $\Phi^d: D \ra [0, 1]$ a function with bounded 
derivatives such that   
\[ 
\Phi^d(x_1, u^1) =  \left\{ \begin{array}{cc}
                 1, & (x_1, u^1) \in D^d \\
                 0, & (x_1, u^1) \in D - N(D^d).
                 \end{array} \right.
\]
Since 
$\sqrt{k + (u^1)^2 - x_1}$ is defined on $N(D^d)$,  we define 
\[
g^s(x_1, u^1) = 
(1-\Phi^d(x_1, u^1))\sqrt{\delta/2} + 
\Phi^d(x_1, u^1) \sqrt{k + (u^1)^2 - x_1}.
\]   
Then $g^s$ is a positive-valued smooth function such that 
\begin{eqnarray}
g^s(x_1, u^1) &=& \sqrt{k + (u^1)^2 - x_1} 
\quad\mbox{for}\quad (x_1, u^1) \in D^d \nonumber,\\
g^s(x_1, u^1) &=& \sqrt{\delta/2} 
\quad\mbox{for}\quad (x_1, u^1) \in D - N(D^d),\\  
g^s(x_1, u^1) &\geq& \sqrt{k + (u^1)^2 - x_1} 
\quad\mbox{for}\quad (x_1, u^1) \in N(D^d) - D^d. \nonumber    
\end{eqnarray}

Since $\delta = 3\eps$, the line in $D$ 
given by $x_1 = k + \eps$ intersects  
the arc $x_1 = k + (u^1)^2 - \delta/2$ in $D$ at two points. 
Let $R^\eps$ be the region in $D$ 
given by $x_1 \leq k + \eps$ and 
$x_1 \geq k + (u^1)^2 - \delta/2$.   
Let $G^+: D \ra B_{\delta}(k, 0, 0)$ 
be given by $G^+(x_1, u^1) = (x_1, u^1, g^s(x_1, u^1))$ 
and $G^-: D \ra B_{\delta}(k, 0, 0)$ 
by $G^-(x_1, u^1) = (x_1, u^1, -g^s(x_1, u^1))$
for $(x_1, u^1) \in D$.  
Then $G^+$ and $G^-$ are embeddings into $B_{\delta}(k, 0, 0)$. 
Let $\cal E$ be the foliation on $B_\delta(k, 0, 0)$ 
with leaves that are level sets under $x_1$. 
Let $G^+(R^\eps)$ be denoted by $D^+$, and $G^-(R^\eps)$ by $D^-$. 
Then $D^+ \cup D^-$ are disjoint imbedded disks
transverse to $\cal E$. 

\begin{figure}[h] 
\centerline{\epsfysize=3cm
\epsfbox{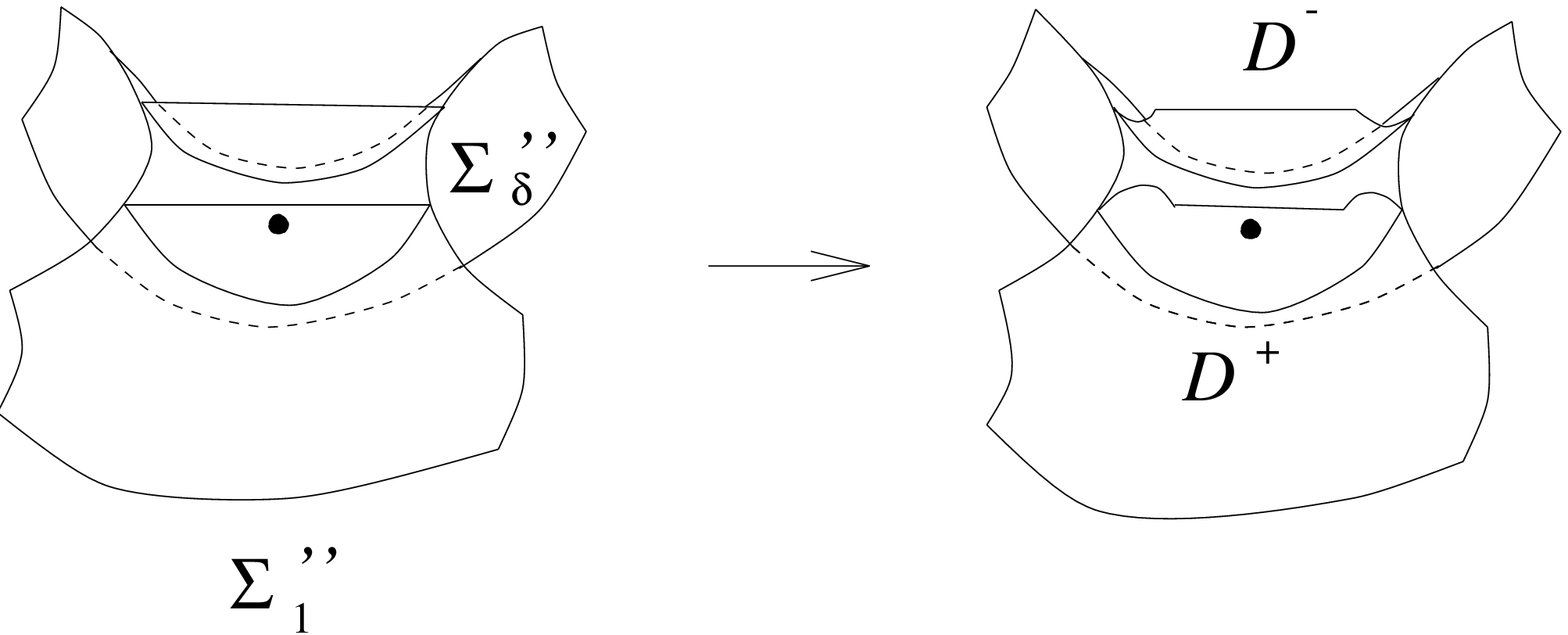}}
\caption{\label{fig:uf8}} 
\end{figure} 

Let us recall 
\[\Sigma'' = \{(x_1, u^1, u^2) \in B_{\delta'}(k, 0, 0) | 
x_1 = k + (u^1)^2 + (u^2)^2 \}; \]
let 
\begin{eqnarray*}
\Sigma''_1 &=& 
\{(x_1, u^1, u^2) \in \Sigma''| x_1 \leq k+\eps\},\\ 
\Sigma''_\delta &=& \{(x_1, u^1, u^2) \in \Sigma''_1|  
-\sqrt{\delta/2} \leq u^2 \leq \sqrt{\delta/2}\}, 
\end{eqnarray*}
so that $(\Sigma''_1 - \Sigma''_\delta) \cup D^+ \cup D^-$ 
is a smooth manifold properly imbedded in $B_{\delta}(k, 0, 0)$ 
diffeomorphic to the union of two disjoint disks
transverse to $\cal E$, a result of a surgery, 
since it is the union of 
\[G^+(\{(x_1, u^1) \in D| x_1 \leq k+\eps\}) 
\hbox{ and } G^-(\{(x_1, u^1) \in D| x_1 \leq k+\eps\}),\] 
two manifolds diffeomorphic to disks.  

Let $P$ be the plane in $B_\delta(k, 0, 0)$ given by $x_1 = k+\eps$. 
Then $P \cap \Sigma''_1$ is the union of two 
arcs on $P$ given by $k+\eps = k + (u^1)^2 - (u^2)^2$.   
Recall that $R^\eps$ is the region in $D$ 
given by $x_1 \leq k + \eps$ and 
$x_1 \geq k + (u^1)^2 - \delta/2$.   
Let $I_\eps$ be the segment 
consisting of points of $R^\eps$ whose 
$x_1$-values equal to $k+\eps$. Then the union of 
$P \cap \Sigma''_\delta$, $G^+(I_\eps)$, and $G^-(I_\eps)$ 
is the boundary of a closed disk $E'$ in $P$,  
and the union of $E', \Sigma''_\delta, D^+$, 
and $D^-$ is the boundary of 
a compact three-ball $\cal B$ in $B_\delta(k, 0, 0)$.  

\begin{figure}[h] 
\centerline{\epsfysize=6cm
\epsfbox{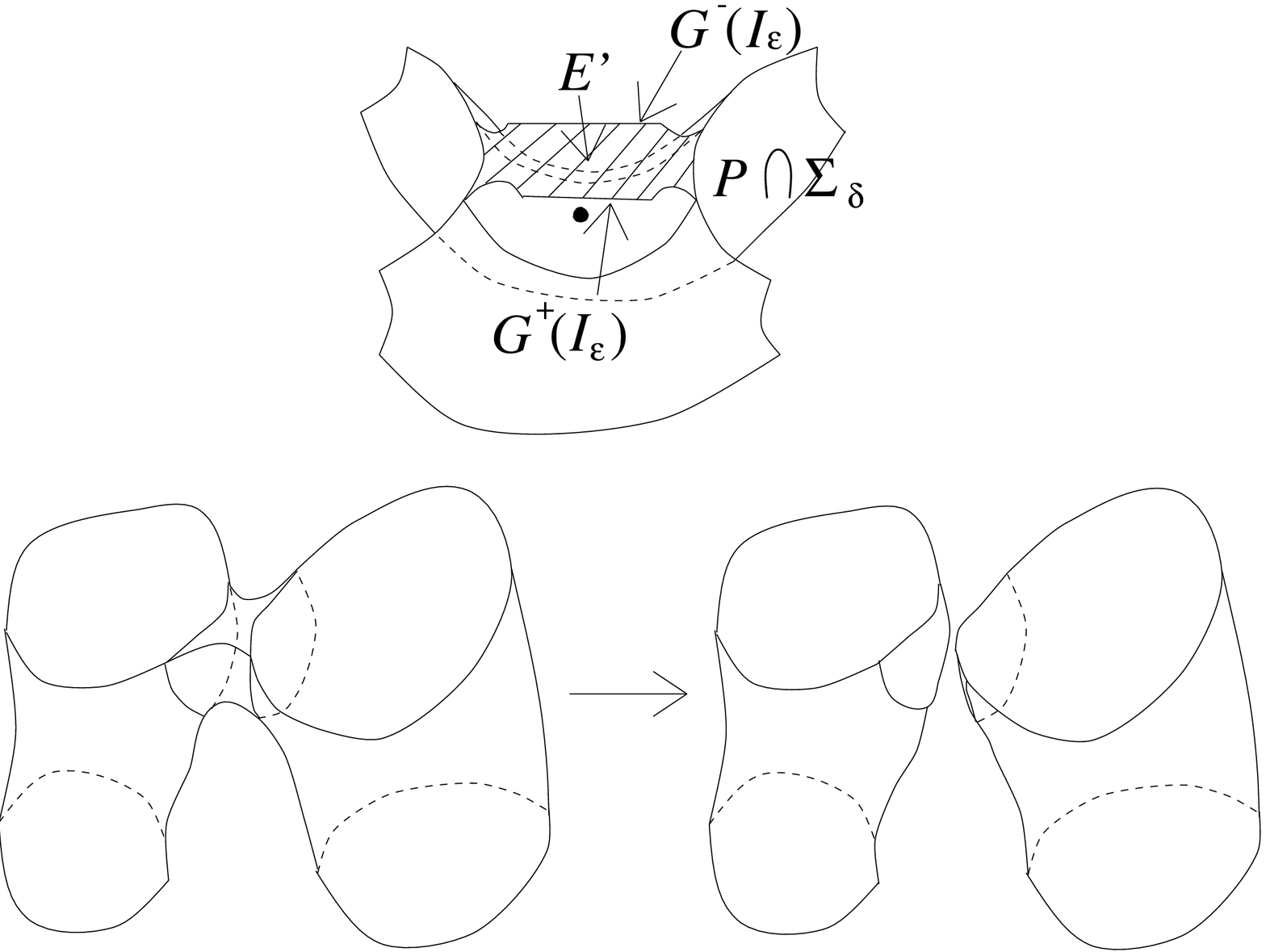}}
\caption{\label{fig:uf9}} 
\end{figure} 

(iii) {\em We do the smooth surgery now}. 
Let ${\cal D}^+$, ${\cal D}^-$, $E$, and $\Sigma_\delta$ 
be $(F| U)^{-1}(D^+)$, $(F| U)^{-1}(D^-)$, 
$(F| U)^{-1}(E')$, and $(F| U)^{-1}(\Sigma''_\delta)$
respectively. 
Recall that $\Sigma_1 = f_1^{-1}([t, k+\eps])$. Let 
\[
\Sigma_1^s = 
(\Sigma_1 - \Sigma_\delta) \cup {\cal D}^+ \cup {\cal D}^-,
\] 
which is the result of a surgery and is diffeomorphic to 
the union of two disjoint annuli. 
Since $(F| U)^{-1}$ maps leaves of 
$\cal E$ into leaves of $\cal F$, $\Sigma_1^s$ is 
transverse to $\cal F$. 

(iv) {\em We now try to find the top face for the above annuli}.
Recall that $M_{k+\eps}$ includes the disk $D_{k+\eps}$ 
with $\delta D_{k+\eps} = \delta \Sigma_0 = \delta_1 \Sigma_1$ 
from the beginning of the proof.  
$\delta D_{k+\eps} \cap E$ is a union of two compact arcs
$\alpha_1$ and $\alpha_2$ in $\delta_1 \Sigma_1$
where $\alpha_1 \cup \alpha_2 = 
\Sigma_\delta \cap \delta_1 \Sigma_1$.   
Since $E - \alpha_1 -\alpha_2$ is included 
in one of two components of 
$M_{k+\eps} - \delta_1\Sigma_1$,  
we have the following two cases: 
\begin{itemize}
\item[(a)] $E - \alpha_1 - \alpha_2 \subset D^o_{k+\eps}$,
\item[(b)] $(E -\alpha_1 -\alpha_2) \cap D_{k+\eps} = \emp$.  
\end{itemize}

\begin{figure}[h] 
\centerline{\epsfysize=3cm
\epsfbox{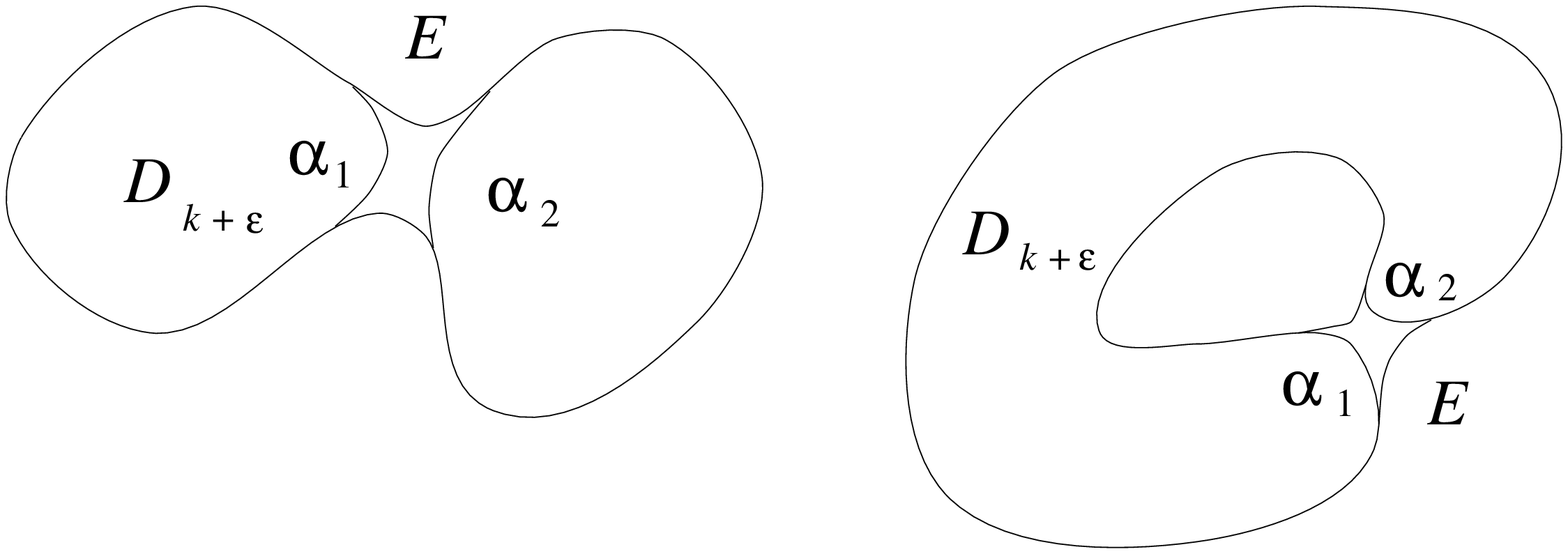}}
\caption{\label{fig:uf10} $D_{k+\eps}$ and $E$.}
\end{figure} 

(a) In this case 
$D_{k+\eps} - (E^o \cup \alpha_1^o \cup \alpha_2^o)$ 
is the union of two disjoint compact disks, 
which is denoted by $\cal O$. 
Then the boundary of $\cal O$ equals 
\begin{equation}\label{eqn:surg}
(\delta D_{k+\eps} - \Sigma''_\delta)    
\cup (F|U)^{-1}(G^+(I_\eps)) 
\cup (F|U)^{-1}(G^-(I_\eps)),
\end{equation}
which is the union of boundary components of $\Sigma_1^s$ 
at the level $k+\eps$. 

We see that $\cal O$ forms the top of the drum 
${\cal O} \cup \Sigma^s_1$ with the bottom side removed.
Since $\Sigma^s_1$ is transverse to $\cal F$, 
Theorem \ref{thm:fill} shows that  
there exists a surface in $M_t$ 
which is the union of two disjoint compact disks   
and whose boundary equals $\delta_2 \Sigma_1$.  
Call the disks $D_1$ and $D_2$. 
Then ${\cal O} \cup \Sigma_1^s \cup D_1 \cup D_2$ is 
the boundary of a compact subset $N_2$ of $M$, 
which is the union of two disjoint subsets homeomorphic 
to three-balls. 
Let 
\[B_2 = N_2 \cup (F| U)^{-1}({\cal B}).\] 
Then $B_2$ is homeomorphic to a three-ball.  

Recall that the union of $\Sigma_0 = f_1^{-1}([k + \eps, l])$,  
where $l$ is the maximal value of $f_1$ on $\Sigma$,
and the disk $D_{k+\eps}$ in $M_{k + \eps}$ 
equals the boundary of a compact three-ball $B_1$ in $M$. 
Since $\delta B_2$ equals 
$D_{k+\eps} \cup \Sigma_1 \cup D_1 \cup D_2$,  
the balls $B_1$ and $B_2$ meet exactly at $D_{k+\eps}$.  
Consequently, $B_1 \cup B_2$ is homeomorphic to 
a three-ball. Moreover, the boundary 
of $B_1 \cup B_2$ equals $\Sigma \cup D_1 \cup D_2$.  
Hence, we obtain $(\alpha)$. 

(b) In this case, $D_{k+\eps} \cup E$ is a compact 
annulus in the level $k+\eps$ with smooth boundary,
as described in equation \ref{eqn:surg}.    
As in (a), Theorem \ref{thm:fill} shows that 
there exists a compact annulus $\Omega_2$
with boundary $\delta_2 \Sigma_1$ in the level $t$ 
such that the union of $D_{k+\eps} \cup E$,  
$\Omega_2$, and $\Sigma_1^s$ is the boundary of 
a compact three-manifold $N_2$ homeomorphic to 
$\Omega_2 \times [t, k+\eps]$, a solid torus.        

Let \[B_2 = (N_2 - (F| U)^{-1}({\cal B})) \cup \Sigma_\delta.\] 
While $B_2$ is a compact three-manifold still 
homeomorphic to a solid torus,
the boundary of $B_2$ equals the union of $D_{k+\eps}$, 
$\Omega_2$, and $\Sigma_1$. 
We can show as in (a) that $B_1 \cup B_2$ is  
homeomorphic to a solid torus, and 
the boundary of $B_1 \cup B_2$ 
equals $F(A) \cup \Omega_2$. 
Hence, there exists an embedding $G: A \ra \Omega_2$ 
homotopic to $F$ relative to $\delta A$.   
We obtained $(\beta)$. 
\end{pf}

\begin{figure}[h]
\centerline{\epsfysize=11cm
\epsfbox{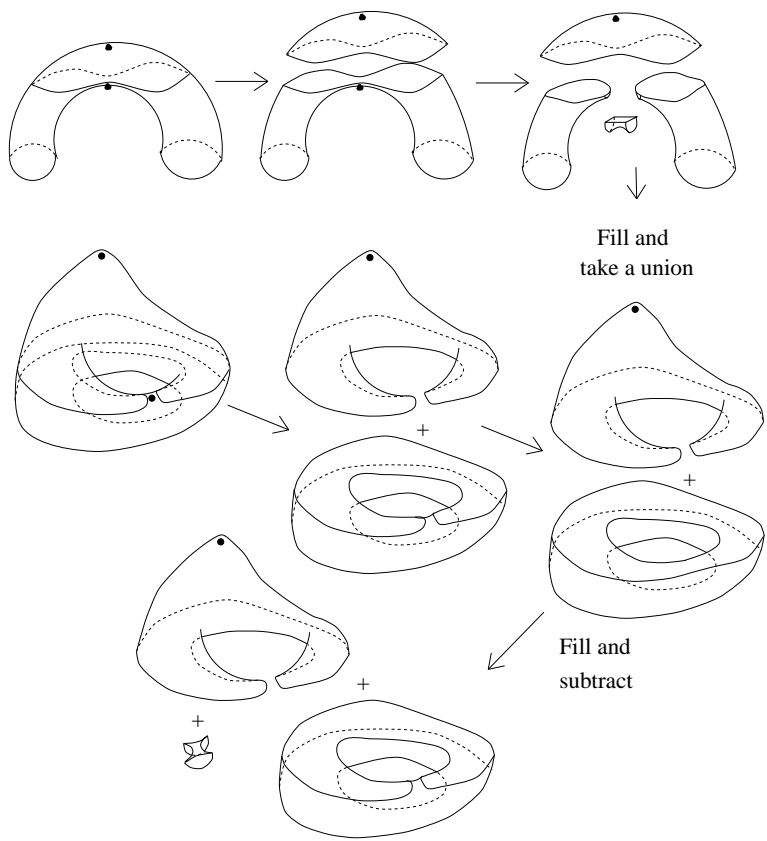}}
\caption{\label{fig:movie0} Our filling processes} 
\end{figure}

\section{Replacing disks and controlling the relative 
maximum points}
\label{sec:control}

In this section, we show that if a portion of 
an amenable embedding $F$ is an amenable embedding 
homotopic to an embedding of a disk into a level set 
then we can replace this portion by an amenable 
disk with one critical point and still obtain
an amenable embedding (see Proposition \ref{prop:smoothing})
provided the interior of the imbedded disk is disjoint from 
the image of $F$. We modify this result further a bit in order 
to control the location and the size of the second 
relative maximum point using Lemma \ref{lem:diskmap} 
(see (iv) and (v) of Proposition \ref{prop:smoothing} below). 
We even produce gradient-like vector fields for these functions.

Let $F$ be an amenable embedding of 
a surface $\Omega$. Suppose that 
$\Omega_1$ is a subsurface of $\Omega^o$ 
diffeomorphic to a disk
so that $F| \Omega_1$ is a positive amenable 
embedding, of level $t$, 
homotopic relative to $\delta\Omega_1$ 
to an embedding $G: \Omega_1\ra W$ 
for a subsurface $W$ in $M_t$. 
Let us fix orientations of $\Omega$ 
and $\bR^2$ so that $(x_2\circ G, x_3\circ G)$ is 
an orientation preserving map into $\bR^2$. 
We let $f = x_1 \circ F$.

\begin{prop}\label{prop:smoothing}
Suppose that $W$ includes no components of 
$F(f^{-1}(t))$ other than $F(\delta \Omega_1)$.
Let $\delta$ be an arbitrary positive number. 
Then there exists a smooth amenable embedding 
$F': \Omega \ra M$ with following properties~{\/\rm :\/}  
\begin{itemize}
\item[(i)] $F'$ is homotopic to $F$ relative to  
$\Omega - \Omega_1^o$.
\item[(ii)] The maximum point $y$ of $x_1\circ F'| \Omega_1$ 
is unique with maximum value less than $t + \delta$. 
\item[(iii)] The number of critical points of $x_1\circ F'$ 
is less than or equal to that of $x_1\circ F$. 
\end{itemize}
\end{prop}
\begin{pf} 
Assume without loss of generality that $t >1$. 
Since $W$ is a compact subset of $M_t$, 
a compact neighborhood $N$ in $M$ includes $W$. 
Let $N_t(W)$ be a compact-disk neighborhood of 
$W$ in the interior of $M_t \cap N$. 
Hence, there exists a $\delta'$, $0 < \delta' <1/2, \delta$,
so that $p \in M$ obtained by starting from 
a point $x \in N_t(W)$ and going in the positive 
or negative $x_1$-direction 
by a distance $s$ less than or equal to $\delta'$ still lies in
$N$. Let this point $p$ be denoted by $x \pm se_1$
where the sign depends on the direction.   
The set of all such points is denoted 
by $N(W)$, which is a neighborhood of $W$
homeomorphic to $N_t(W) \times [t-\delta', t+\delta']$.   

Let $N_{t'}(W)$ denote $M_{t'} \cap N(W)$ for 
$t' \in [t-\delta', t+\delta']$. 
We also require that the interval $[t -\delta', t+\delta']$ 
contains no critical value of $f$. 
Since we assumed that $W$ includes no components of 
$F(f^{-1}(t))$ other than $F(\delta \Omega_1)$, we 
can choose $\delta'$ sufficiently small 
so that $F^{-1}(N_{t'}(W))$ is a simple closed curve 
$\alpha_{t'}$, which is a component of the $t'$-level set
of $x_1\circ F$.   
Thus, $F^{-1}(N(W))$ is a compact annulus 
that is a neighborhood of $\delta \Omega_1$
and is foliated by simple closed curves $\alpha_{t'}$. 
In particular, whenever $t', t'' \in [t-\delta', t+\delta']$,
the union of two level curves 
$\alpha_{t'}$ and $\alpha_{t''}$, $t' \ne t''$, 
is the boundary of a compact annulus $A_{t', t''}$
included in $F^{-1}(N(W))$. 

\begin{figure}[h] 
\centerline{\epsfysize=3.5cm
\epsfbox{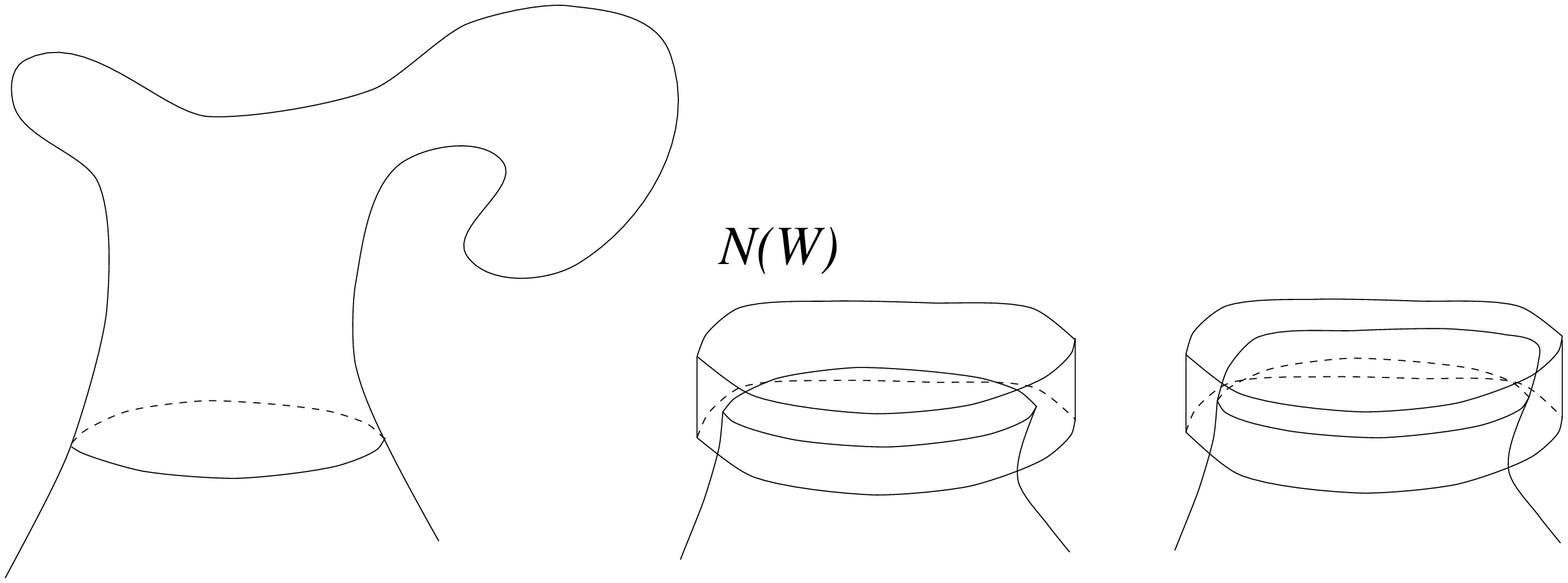}}
\caption{\label{fig:uf12} Smoothing.} 
\end{figure} 

{\em We will now construct an embedding 
from $\Omega_1$ into $N(W)$.\/}
Assume without loss of generality that $t > 0$.
There is a diffeomorphism 
\begin{equation}\label{eqn:Phi}
\Phi: A_{t, t+\delta'/2} \ra E
\end{equation}
where $E$ is the annulus in $\bR^2$ whose boundary is 
the union of  
two circles of respective radii $t$ and $t-\delta'/2$ with 
centers at the origin $O$ 
mapping $\alpha_{t'}$ to a circle of radius 
$t - (t' -t)$ with center $O$. 
Clearly $\Phi$ extends to a diffeomorphism 
of $\Omega_1$ with the disk $D$ in $\bR^2$ of radius $t$ with 
center $O$. 

We choose a smooth function $f':D \ra [0, 3\delta/4]$
as follows: 
\begin{equation}\label{eqn:fprime}
f'(x, y) = t - \sqrt{x^2 + y^2}  
\quad \mbox{for} \quad t-\delta'/2 \leq \sqrt{x^2 + y^2} \leq  t
\end{equation}
and constant along circles with center $O$ 
with unique critical point $O$ taking the maximum value $3\delta'/4$.
Then $f' \circ \Phi$ has constant value $t' -t$ on 
$\alpha_{t'}$ for $t' \in [t, t + \delta'/2]$, and  
$\Phi^{-1}(O)$ is the maximum point of $f'\circ \Phi$. 
Consequently, $\Phi^{-1}(O)$ does not belong to $A_{t, t+\eps}$
for $\eps < \delta'/2$. 

Let us define a map $G': \Omega_1 \ra N(W)$ 
by $G'(x) = G(x) + (f' \circ \Phi(x))e_1$. 
Since this can be seen as a graph map, $G'$ 
is an embedding. Moreover, 
\begin{equation}\label{equ:x1com}
x_1\circ G' = x_1\circ F \quad\mbox{on}\quad A_{t, t+\delta'/2} 
\end{equation}
since $x_1\circ G'$ has constant value $t'$ on $\alpha_{t'}$.   
Therefore, $G'|\alpha_{t'}$ is an embedding onto 
a simple closed curve in $M_{t'}$ for 
$t' \in [t, t+\delta'/2]$, while 
so is $F|\alpha_{t'}$.  

Let $\eps$ be a positive number so that $2 \eps$ 
is less than $\bdd(\Omega, \delta N(W))$.
Since $G'|\delta \Omega_1 = F|\delta \Omega_1$, 
it follows that for every $\eps$, $\eps > 0$, 
there exists $\gamma$, $0 < \gamma < \delta'/2$, such that  
\begin{eqnarray}\label{eqn:small}
\bdd(F(x), F(\delta \Omega_1)) &< & \eps, \nonumber \\
\bdd(G'(x), F(\delta \Omega_1)) &< & \eps, \\
\bdd(G'(x), F(x)) & < & \eps
\quad\mbox{if}\quad x \in A_{t, t+\gamma}. \nonumber
\end{eqnarray}  
Each point of $F(A_{t, t+\gamma}) \cup G'(A_{t, t+\gamma})$ 
is contained in an open $\eps$-$\bdd$-ball included in $N(W)$. 
Let $\eps_1 < \gamma$. Since the balls in the cover are convex,  
we obtain 
\[hF(x) + (1-h)G'(x) \in N(W) \hbox{ if } x \in A_{t, t+\eps_1}\] 
for every real number $h$, $0 \leq h \leq 1$.

{\em We smooth $G'$ and $F$ near $\delta \Omega_1$.}
Let $\eps$ be a real number with 
$0 < \eps < \eps_1, \delta'$, and
$\phi_\eps: \Omega \ra [0, 1]$ a smooth function  
such that $\phi_\eps$ is $1$ on 
$(\Omega - \Omega_1) \cup A_{t, t+\eps/4}$ 
and $0$ on the compact disk $\Omega^o_1 - A_{t, t+\eps}^o$
and is constant along $\alpha_{t'}$ for 
$t' \in [t, t+\eps]$.      
We define $F_\eps: \Omega \ra N(W)$ by 
\[
F_\eps(x) = \phi_\eps(x)F(x) + (1-\phi_\eps(x))G'(x). 
\] 
Then $F_\eps$ equals $F$ on $\Omega - \Omega_1^o$
and equals $G'$ on $\Omega^o_1 - A_{t, t+\eps}^o$.  
In particular, equation \ref{equ:x1com} shows that 
\begin{equation}\label{equ:x1com2}
x_1\circ G' = x_1\circ F_\eps \quad\mbox{on}\quad \Omega_1. 
\end{equation}
Moreover, $F_\eps$ is homotopic with $F$ relative to  
$\Omega - \Omega_1^o$ since $F_\eps$ depends 
continuously on $\eps$ for $\eps > 0$ 
and, as $\eps \ra 0$, $F_\eps|\Omega_1$ converges 
to $G'|\Omega_1$ which is homotopic to 
$F|\Omega_1$ relative to $\delta \Omega_1$. 

{\em We now show that $F_\eps$ is an embedding}.
Let us orient the curves $\alpha_{t'}$ continuously 
for $t \leq t' \leq t+\delta$.
Since the derivative of $F$ along $\alpha_t$ is 
identical with that of $G'$,  
a continuity argument shows that  
there exists $\eps_2$ so that for $\eps \leq \eps_2$, 
the derivative of $F$ on $A_{t, t+\eps}$ 
along $\alpha_{t'}$ has positive $\mu$-inner
product with that of $G'$. 
Thus, the derivative of $F_\eps$ along $\alpha'_t$ 
is a nonzero vector field. Since the derivatives of $F$ and $G'$ along 
$\alpha_{t'}$ have zero $x_1$-components, 
that of $F_\eps$ has zero $x_1$-component.  

Since for $0 < \eps < \delta'$, at each point of $A_{t, t+\eps}$, 
the $x_1$-component of the derivative of $G'$ 
in the inner direction orthogonal to 
$\alpha_{t'}$ is positive,  
equation \ref{equ:x1com2} shows the positivity of  
that of $F_\eps$. 
Therefore, $F_\eps$ is an immersion at 
points of $A_{t, t+\eps}$ for $\eps < \eps_1, \eps_2, \delta'$.  
Since $F$ and $G'$ are immersions on $\Omega - \Omega_1$ 
and on $\Omega_1 - A_{t, t+\eps}$ respectively, and 
$F_\eps$ equals $F$ or $G'$ on the respective sets,  
$F_\eps$ is an immersion on all of 
$\Omega$ for $\eps < \eps_1, \eps_2, \delta'$.  

Since we have $F_\eps|\alpha_t = F|\alpha_t= G|\alpha_t$, 
Theorem 2.1.4 \cite[p. 37]{Hirsch} shows that 
there exists a positive constant $\eps_3$ 
so that $F_\eps|\alpha_{t'}$ is an embedding 
into $M_{t'}$ for every $\eps$ and $t'$ with $0<\eps <\eps_3$ and
$t'\in [t, t+\eps]$. 
This implies that $F_\eps$ is the desired embedding
for $\eps < \eps_i, \delta'$, $i=1, 2, 3$.  
\end{pf}

\begin{lem}\label{lem:diskmap}
Let $\disk$ be a compact disk in the plane $\bR^2,$
$z$ a point of $\delta \disk,$ and $W$ 
a compact disk with smooth boundary in an open surface.   
Let $U$ be a compact-disk neighborhood of $z$ 
in $\disk$ where $U \cap \delta \disk$ is 
diffeomorphic to a closed interval. 
Suppose that $f: U \cup \delta \disk$ 
is an imbedding into $W$ so that 
$f(\delta \disk) \subset \delta W$
and $f$ extends to a smooth map $f'$ in 
a neighborhood of $U \cup \delta \disk.$  
Then there exists a diffeomorphism 
$g:\disk \ra W$ extending $f.$ 
\end{lem} 
\begin{pf} 
We claim that $f$ extends to an embedding $f''$ of
an open neighborhood $U'$ of $U\cup \delta \disk$ 
in $\disk$. We may assume without loss of 
generality that $\disk$ and $W$ equal the standard disk 
in $\bR^2$. One can easily extend $f| \delta \disk$ to  
an immersion $f_1$ from a neighborhood $N_1$ of 
$\delta \disk$ in $\disk$ to a neighborhood of $\delta W$ 
in $W$. Let $N$ be an open neighborhood of $U$ in $\disk$
so that $f'|N: N \ra W$ is a well-defined immersion
extending $f| U \cup (N \cap \delta {\cal D})$. 
We can do this since the ranks of the derivatives of $f$ 
at points of $U$ equal $2$ and hence those of $f$ at 
points of $\delta {\cal D}$ sufficiently near points of $U$ 
are also two.  
Let $\psi: \disk \ra [0, 1]$ be a smooth function such that 
$\psi$ equals $0$ on $U$ and $1$ 
on $\disk - N$. Define a map $f'': N_1 \cup N \ra W$ 
by $f''(x) = \psi(x)f_1(x) + (1 - \psi(x))f'(x)$, so that      
$f''$ equals $f$ on $U \cup \delta \disk$. 
Since $(f'|N) | N \cap \delta \disk$ equals 
$f| N \cap \delta \disk$,  
the derivative of $f_1$ in the angular direction 
equals that of $f'|N$ on $N \cap \delta \disk$.   

This and the fact that 
$f_1$ and $f'|N$ are immersions show that  
$\disk$ includes a sufficiently thin neighborhood $N_2$
of $\delta \disk$ such that 
$f''$ restricted to $N_2$ is an immersion
as an explicit matrix calculation in polar 
coordinates will verify: 
Give the polar coordinates on $\cal D$ and $W$ so 
that $\theta$ increases along the boundary orientation of 
$\cal D$ and $r$ decreases in the radial inward 
direction of the standard disk $\cal D$. 
Assume that $f_1$ and $f'$ preserves orientation.
Then under the coordinate system $(r, \theta)$, 
\[ D f_1 = \left( \begin{array}{cc} \alpha & 0 \\ \beta & b
\end{array} \right), \quad\quad  
D f' = \left( \begin{array}{cc} \alpha' & 0 \\ \beta' & b' 
\end{array} \right) \] 
on $\delta \cal D$ for $\alpha, b, \alpha', b' > 0$  
since $\delta {\cal D} \mapsto \delta W$ 
and the inward vector of $\cal D$ 
is mapped to inward vector of $W$ for both of 
$f_1$ and $f'$.
We have 
\[ Df'' = \left( \begin{array}{cc} \psi_r(f_{11} - f'_1) +
\psi \alpha + (1-\psi) \alpha' & 0 \\
* & \psi_\theta(f_{12} - f'_2) + \psi b + (1-\psi) b' 
\end{array} \right) \] 
where $f_{11}$ and $f'_1$ denote the $r$-component of 
$f_1$ and $f'$ respectively and $f_{12}$ and $f'_2$ 
denote the $\theta$-components of $f_1$ and $f'$ 
respectively. Since $f_1 = f'$ on $\delta {\cal D} \cap N$, 
it follows that $D f''$ is not zero on $\delta \cal D$.
Therefore, $f''$ is an immersion on a neighborhood of $\delta \cal D$.

Since $f''| U \cup \delta \disk$ is injective, $N \cup N_2$ 
includes a neighborhood $U'$ of $U \cup \delta \disk$ in 
$\disk$ so that $f''|U'$ is an embedding.  

In $U'$, there is a neighborhood $A$ of $U \cup \delta \disk$, 
diffeomorphic to a compact annulus, with 
one boundary component equal to $\delta \disk$ 
and the other in $\disk^o$.  
Then $f''|A$ is a diffeomorphism of $A$ into 
an annulus $A'$ in $W$ whose one boundary 
component is $\delta W$. 
Since $\disk^o - A^o$ is a compact disk with 
smooth boundary and so is $W^o - f(A')^o$, 
there is a diffeomorphism $f'''$ between these disks
extending the map $f''$ restricted to the component 
of $\delta A$ in $\disk^o$.  
One can now smooth $f''$ and $f'''$ to 
a diffeomorphism $\disk \ra \disk$ 
(see Theorem 8.1.9 \cite[p. 182]{Hirsch}).
\end{pf} 

We will need two definitions: 
A smooth vector field $\zeta$ on $\Omega$ is 
an {\em $f$-gradient-like vector field\/} 
for a real-valued Morse function $f$ if 
\begin{itemize}
\item $\zeta(f) > 0$ in the complement of critical points of $f$, and 
\item given each critical point $p$ of $f$, there are coordinates 
$(u^1, u^2)$ in a neighborhood $U$ of $p$ so that 
$f= f(p) +(-1)^p(u^1)^2 +(-1)^q(u^2)^2$, and $\zeta$ has coordinates 
$((-1)^p u^1, (-1)^q u^2)$ throughout $U$, where $p$ and $q$ 
are integers equal to $1$ or $2$.  
\end{itemize}
(See Milnor \cite[p. 20]{Mil}.)

A {\em cubical neighborhood} of a point 
$x$ in $M$ is an open convex ball $B$
with compact closure $\clo(B)$ in $M$
so that $\dev| B$ is an embedding onto the set 
\[\{(x_1, x_2, x_3) \in \bR^3| |x_i - a_i| < b 
\quad \mbox{for every} i = 1, 2, 3\}\] 
for some real numbers $a_1, a_2, a_3,$ and $b$, $b > 0$. 
A {\em bottom side\/} of $B$ 
is the side of $B$ corresponding to $x_1 = a_1 - b$.
$B$ has natural coordinates $x_1\circ \dev, 
x_2 \circ \dev,$ and $x_3 \circ \dev$, 
which were denoted by $x_1, x_2$, and $x_3$.  

{\em Proposition \ref{prop:smoothing} can be extended as follows\/}:  

Let $B$ be a cubical neighborhood of 
$z \in \delta \Omega_1$. 
Suppose that there exists an $x_1\circ F$-gradient-like 
vector field $\zeta$ on $\Omega - \Omega_1^o$ 
so that the map $(x_2\circ F, x_3\circ F)$ restricted on 
a compact-disk neighborhood $V$ of $z$ in $\Omega_1$, 
$V \subset F^{-1}(B)$, is 
an orientation-preserving diffeomorphism onto 
an open subset of $\bR^2$ with the given orientation. 
(Recall that we fix an orientation of $\Omega$ 
and $\bR^2$ so that $(x_2\circ G, x_3\circ G): \Omega_1 \ra \bR^2$ is 
orientation-preserving.) 
Then we can choose $F'$ so that the following 
statements hold in addition to (i), (ii), and (iii) 
in Proposition \ref{prop:smoothing}: 
\begin{itemize}
\item[(iv)] There exists an $x_1\circ F'$-gradient-like vector 
field $\zeta'$ on $\Omega$ extending $\zeta$ smoothly  
so that the closure $\Gamma$ of 
the $\zeta'$-trajectory from $z$ to the relative maximum point $y$
of $x_1 \circ F'$ is a smooth arc and 
lies in the interior of $V$ in $\Omega_1$. 
(We have $F'(y) \in B$.)
\item[(v)] 
$(x_2\circ F', x_3\circ F')|V = (x_2\circ F, x_3\circ F)|V$.  
\end{itemize}

\begin{pf}
To prove (iv) and (v), we will change $F_\eps$ and $G$ 
in the proof of Proposition \ref{prop:smoothing}. 
Since $(x_2 \circ F, x_3 \circ F)| V$ maps into $W$
by the orientation condition, 
Lemma \ref{lem:diskmap} shows that   
there is a diffeomorphism $G_1: \Omega_1 \ra W$ so that 
\begin{eqnarray}\label{eqn:fixedmap}
(x_2\circ G_1, x_3\circ G_1)| V &=& (x_2\circ F, x_3\circ F)| V \\
G_1|\delta \Omega_1 &=& F|\delta \Omega_1.
\end{eqnarray}

Recall $f'$ and $\Phi$ from equations \ref{eqn:fprime} 
and \ref{eqn:Phi} respectively.
Since $\Phi$ extends to a neighborhood of $\Omega_1$, 
it follows that there exists an $f'$-gradient-like 
vector field $\zeta_1$ on $D$ so 
that $\zeta_1$ pulls back to an $f'\circ \Phi$-gradient-like 
vector field on $\Omega_1$ extending $\zeta$. 

Let $V' = \Phi(V)$, and $\Gamma_1$ the closure of 
the $\zeta_1$-trajectory from $\Phi(z)$ to $O$. 
There exists $\eps_4$, $\eps_4 > 0$, 
so that $\Gamma_1$ intersected with 
a circle of radius $t-\eps$ with center $O$ for 
$0 \leq \eps \leq \eps_4$ is a unique point 
of the interior of $V'$ in $D$.  
Let $\Psi: D \ra D$ 
be a diffeomorphism fixing all points of distance $t'$
from $O$, $t-\eps_4 \leq t' \leq t$, 
and so that $\Psi(\Gamma_1)$ is a subset of 
the interior of $V'$ in $D$.   
We define $G'_1:\Omega_1 \ra N(W)$ 
by $G'_1(x) = G_1(x) + f'\circ \Psi^{-1} \circ \Phi(x) e_1$. 

Let $\zeta_2$ denote the vector field 
defined on $\Omega_1$ pulled-back from $\zeta_1$ 
by $\Psi^{-1}\circ \Phi$. 
Then $\zeta_2$ extends $\zeta$; $\zeta_2$ 
is $f'\circ\Psi^{-1}\circ \Phi$-gradient-like vector field. 
We let $\zeta'$ 
denote the extended vector field. 
Let $\Gamma$ be $\Phi^{-1}\circ \Psi(\Gamma_1)$, which 
is the closure of the $\zeta'$-trajectory 
from $z$ to $\Phi^{-1}\circ\Psi(0)$, the maximum point 
of $f'\circ \Psi^{-1}\circ \Phi$ on $\Omega_1$,
and $\Gamma$ is a subset of the interior of 
$V$ in $\Omega_1^o$. 

\begin{figure}[h] 
\centerline{\epsfysize=3.5cm
\epsfbox{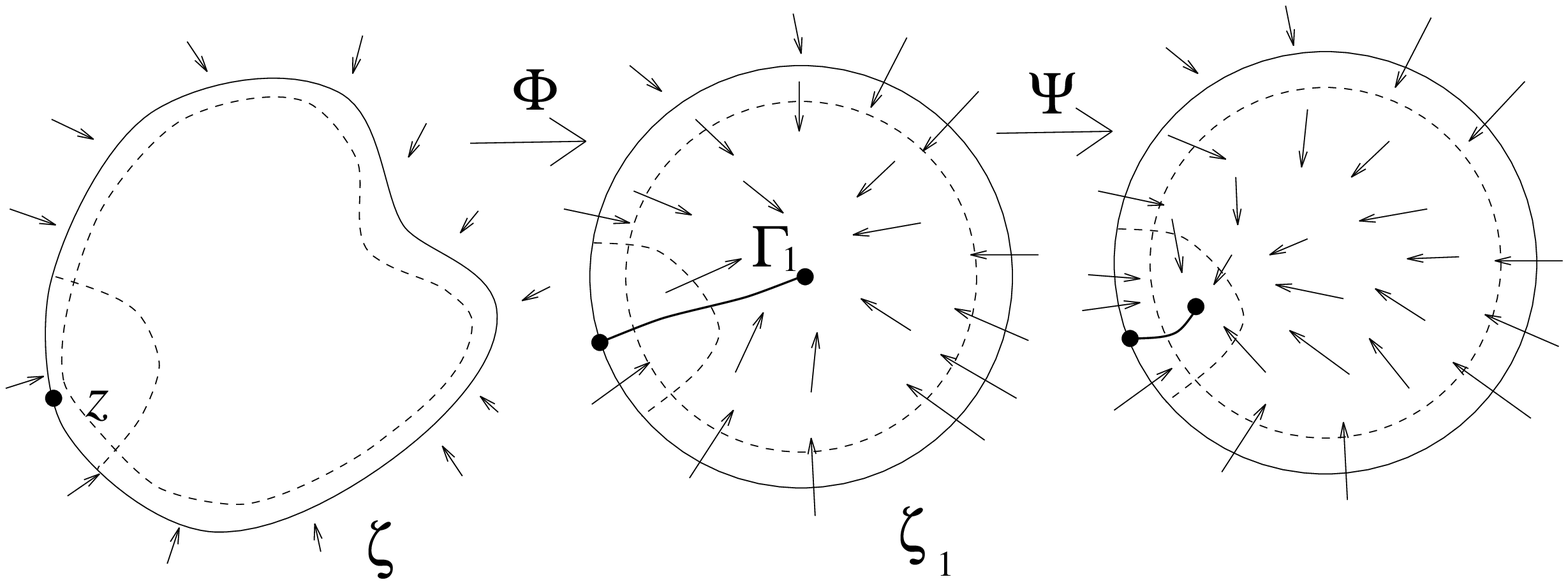}}
\caption{\label{fig:uf13}} 
\end{figure} 

We take the numbers $\eps_i$ for $i=1, 2, 3$ 
from the proof of Proposition \ref{prop:smoothing}. 
Let $\eps$ satisfy $0 < \eps < \eps_i, \delta'$, 
$i=1, 2, 3, 4$. 
We define $F^1_\eps: \Omega_1 \ra N(W)$ by 
\[
F^1_\eps(x) = \phi_\eps(x)F(x) + (1-\phi_\eps(x))G'_1(x). 
\] 
Then $F^1_\eps(x)$ is homotopic to $F$ relative to  
$\Omega - \Omega_1^o$ as in the proof of the part (i) of the lemma.   
As above, $F^1_\eps$ is an embedding satisfying (ii) and (iii), and  
\[
x_1\circ G'_1 = x_1\circ F^1_\eps \quad\mbox{on}\quad \Omega_1,
\] 
exactly as in equation \ref{equ:x1com2}.
Therefore, $\zeta'$ is a gradient-like-vector field 
for $x_1\circ F^1_\eps$ as well. Since $\Gamma$ 
is the closure of the $\zeta'$-trajectory from 
$z$ to the maximum point $\Phi^{-1}\circ \Psi(0)$
of $x_1\circ F^1_\eps$ in $\Omega_1$, 
(iv) is satisfied.    
Now, for $x \in V$, equation \ref{eqn:fixedmap} shows that  
\begin{eqnarray}
(x_2\circ F^1_\eps(x), x_3\circ F^1_\eps(x)) 
&=& \phi_\eps(x)(x_2\circ F(x), x_3\circ F(x)) \nonumber\\ 
&+& (1-\phi_\eps(x))(x_2\circ G_1(x), x_3\circ G_1(x))\\
&=& (x_2\circ F(x), x_3\circ F(x)); \nonumber
\end{eqnarray}
hence, (v) is satisfied. 
\end{pf}

We will need the following corollary later.
\begin{cor}\label{cor:inner}
Let $N(\Omega_1)$ be the neighborhood 
$\Omega_1 \cup A_{t-\delta', t}$ of $\Omega_1$ 
in $\Omega^o.$ Assume that $W$ includes no component of $F(f^{-1}(t)),$ 
other than $F(\delta \Omega_1).$ 
Then there exists an amenable embedding $F'': \Omega \ra M$
with the following properties {\/\rm :}
\begin{itemize}
\item[(i)] $F''$ is homotopic to $F$ relative to 
$\Omega - N(\Omega_1)^o.$
\item[(ii)] The maximum point $y$ of
$x_1\circ F''|N(\Omega_1)$ is unique
with maximum value less than $t$.
\item[(iii)] The number of critical points of $x_1\circ F'$
is less than or equal to that of $x_1\circ F.$
\item[(iv)] The number of components of 
$(x_1\circ F'')^{-1}(t)$ is less than 
that of \break $(x_1\circ F)^{-1}(t)$ by one or more.  
\end{itemize}
\end{cor}         
\begin{pf}
By the above choice of $\delta',$ 
$F| N(\Omega_1)$ is a positive amenable 
embedding of level $t-\delta'$ and  
homotopic to an imbedding into $N(W)$ relative 
to $N(\Omega_1) - \Omega_1^o.$  
Since $N(W)$ is diffeomorphic to 
a three-ball $N_t(W) \times [t-\delta', t+\delta'],$ 
$F|N(\Omega_1)$ is homotopic 
relative to $\delta N(\Omega_1)$ to 
an embedding into $N_{t-\delta'}(W).$  
Now, Proposition \ref{prop:smoothing} implies 
the conclusions of the corollary. 
\end{pf}

\section{Homotopy of a Disk with Three Critical Points}
\label{sec:hdthree}
We prove that an amenable embedding of a disk with 
three critical points is homotopic to 
an embedding into a level set of $M$. 
The proof involves Morse cancellation theory 
in this geometric setting.
Such a disk can either have two relative maximum points 
and a saddle point or have one maximum and 
a relative minimum point and a saddle point.  
In the first case, by using Proposition \ref{prop:smoothing}, 
we try to push the one relative maximum very near 
the saddle and then we cancel them by constructing 
the explicit homotopy in a small convex open set 
called a cubical neighborhood. To do this, we use 
a flow line of a gradient-like vector field 
in the neighborhood connecting 
the relative maximum to the saddle point 
that we constructed at the end of Section \ref{sec:control}
(see Milnor \cite{Mil}).  
In the neighborhood of the union of arcs in the flow lines
containing the saddle point and relative maximum and 
a point above the relative maximum, we find 
an expression of $f$ as an integral and modify the integral
expression. In the second case, 
we use the $3$-ball and the solid torus obtained 
in Proposition \ref{prop:annfill} to obtain homotopies.
The both processes are captured in Figures \ref{fig:mv1} 
and \ref{fig:mv2} respectively.

\begin{prop}\label{prop:diskfill2}
Let $\disk$ be a disk\/{\rm ;}
let $F: \disk \ra M$ be a positive amenable embedding 
of level $t$ with three critical values. 
Then there exists an embedding $G: \disk \ra M_t$ 
homotopic to $F$ relative to $\delta \disk$  
\end{prop}
\begin{pf}
Let $f = x_1\circ F$. 
Since the Euler characteristic of $\disk$ is $1$, 
$f$ has two relative extreme points and a saddle point
by Proposition \ref{prop:Euler}.
Let $v_1, v_2,$ and $v_3$ 
denote critical values, and $z_1, z_2,$ and $z_3$ 
the corresponding critical points in $\disk$ respectively.   
By relabeling if necessary, 
we can assume without loss of generality that 
there are only following cases: 
\begin{itemize}
\item[(a)] $z_1$ is a maximum point, $z_2$ a local maximum point, 
and $z_3$ a saddle point where $v_1 > v_2 > v_3 > t$.  
\item[(b)] $z_1$ is a maximum point, $z_2$ a local minimum point,  
and $z_3$ a saddle point where $v_1 > v_3 > v_2$ 
($v_1, v_3 > t$).
\end{itemize} 

\begin{figure}[h]
\centerline{\epsfysize=6cm
\epsfbox{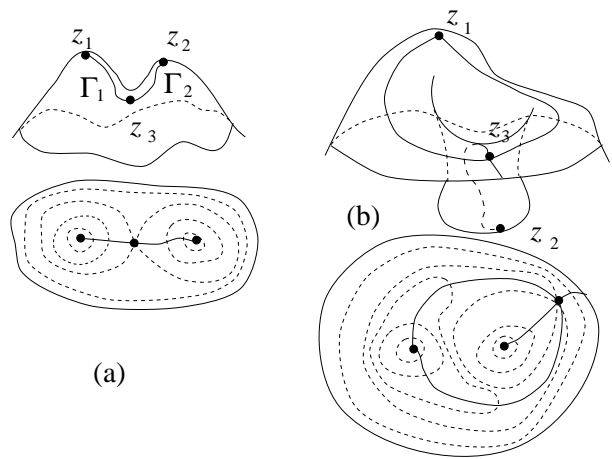}}
\caption{\label{fig:uf14} Examples of (a) and (b).} 
\end{figure} 

(a) Let $B$ be a small cubical neighborhood of $F(z_3)$ in $M$
so that $F^{-1}(B)$ is an open disk neighborhood of $z_3$. 
We can write 
$F(x) = (f(x), x_2\circ F(x), x_3\circ F(x))$ for $x \in F^{-1}(B)$ 
under the natural coordinates of $B$.  
Let ${\cal V} = F^{-1}(B)$ whose closure 
$\clo(\cal V)$ is a compact disk in $\disk^o$.
By choosing $B$ sufficiently small, we may
assume without loss of generality that 
the map defined by \[x \mapsto (x_2\circ F(x), x_3\circ F(x))\] 
is a diffeomorphism of $\clo({\cal V})$ 
onto a closed rectangular disk in $\bR^2$, corresponding to 
the bottom side of $\clo(B)$. (See the proof of 
Proposition \ref{prop:annfill}.)   

Choose $\eps$, $\eps > 0$, so that $v_3 + \eps < v_2$ 
and for the components $\alpha_1$ and $\alpha_2$ of 
$f^{-1}(v_3 + \eps)$, 
which are simple closed curves,  
both $F(\alpha_1)$ and $F(\alpha_2)$ pass $B$. 
Since $v_1, v_2 > v_3$, it follows from 
the flow argument that 
the curves $\alpha_1$ and $\alpha_2$ 
are the respective boundaries of disjoint disks
$\disk_1$ and $\disk_2$ in $\disk^o$.  
We may assume without loss of 
generality that $\disk_1^o$ contains $z_1$,
and $\disk_2^o$ contains $z_2$. 

There exists a coordinate system $u^1, u^2$ 
in a neighborhood $U$ of $z_3$ 
so that 
\[f= v_3 -(u^1)^2 + (u^2)^2
\]
holds in $U$ (see the diagram 5 of Milnor \cite[p. 15]{Mil2}).  
For a gradient-like vector field $\zeta$ for $f$,   
we are able to choose $\eps$ sufficiently small 
and points $r$ and $z$ so that the following statements hold: 
\begin{itemize}
\item 
Let $r$ be the point on $\disk_1^o \cap \cal V$ 
on the $\zeta$-trajectory from 
$z_3$ to $z_1$ such that $f(r) > v_3 + \eps + \delta$ 
for some number $\delta$, $0 < \delta < \eps$.   
The closure $\Gamma_1$ of the $\zeta$-trajectory 
from $z_3$ to $r$ is a subset of $\cal V$.  
($\Gamma_1$ is a smooth simple arc with endpoints $r$ and $z_3$.) 
\item Let $z$ be the point on $\delta \disk_2 \cap \cal V$
on the $\zeta$-trajectory from $z_3$ to $z_2$. 
The closure $\Gamma_2$ of the $\zeta$-trajectory 
from $z_3$ to $z$ is a simple arc in $\cal V$. 
\end{itemize} 

{\em We will now modify $F$\/}.
By Proposition \ref{prop:diskfill1}, $F(\alpha_1)$ and
$F(\alpha_2)$ are respective boundaries of 
disks $E_1$ and $E_2$ in $M_{v_3 + \eps}$. 
Since the union of $E_2$ 
and $F(\disk_2)$ is the boundary of 
a three-ball, if $E_1$ is included in $E_2$, then  
$F(\disk_1)$ has to lie in the three-ball as 
$F$ is an imbedding. This contradicts the assumption that 
$f(z_1)$ is greater than $f(z_2)$. 
Therefore, $E_2$ does not include $E_1$, and 
$E_2$ is disjoint from $F(\alpha_1)$.   

Recall that $(x_2 \circ F, x_3\circ F)$ is 
an imbedding from ${\cal V}$ to the subset of 
$\bR^2$ corresponding to the bottom side of the cubical neighborhood.
(See Figure \ref{fig:orient}.) 
We naturally choose an orientation on $\disk$ so 
that $(x_2 \circ F, x_3\circ F)| \cal V$ is orientation-preserving.
Let $\alpha_2$ be given a boundary orientation from $\disk_2$, 
and so the image of ${\cal V} - \disk_2$ lies to the right of 
\[ (x_2\circ F(\alpha_2 \cap {\cal V}), x_3\circ F(\alpha_2\cap {\cal V}))
\hbox{ in } \bR^2. \] 

Let $G: \disk_2 \ra E_2$ be the imbedding homotopic to $F| \disk_2$. 
We claim that $(x_1 \circ G, x_2\circ G)$ 
on ${\cal V} \cap \disk_2$ is also orientation preserving. 

Suppose not. Then $E_2$ lies to the right of 
$F(\alpha_2)$ in $M_{v_3+\eps}$ since $(x_2 \circ G, x_3\circ G)$ 
is orientation-reversing immersion. 
Hence, $E_2$ is a disk in $M_{v_3 + \eps}$ 
that is the closure of the right component of 
$M_{v_3 + \eps} - F(\alpha_2)$. 
Since $B \cap M_{v_3 + \eps} \cap F(\disk)$ 
is the union of two components $F(\alpha_1 \cap \cal V)$ and 
$F(\alpha_2 \cap \cal V)$ and nothing else, it follows that 
$E_2 \cap B$ is a component of $B \cap M_{v_3 + \eps}$ 
removed with these two arcs, to the right of 
$F(\alpha_2 \cap \cal V)$. 
Since $F(\alpha_1 \cap \cal V)$ lies to the right 
of $F(\alpha_2 \cap \cal V)$, the closure of 
$E_2$ contains some points of $F(\alpha_1)$. This is 
absurd since the boundary of $E_2$ equals $F(\alpha_2)$ 
and $E_2^o$ does not meet $F(\alpha_1)$.

\begin{figure}[h] 
\centerline{\epsfysize=2.75cm
\epsfbox{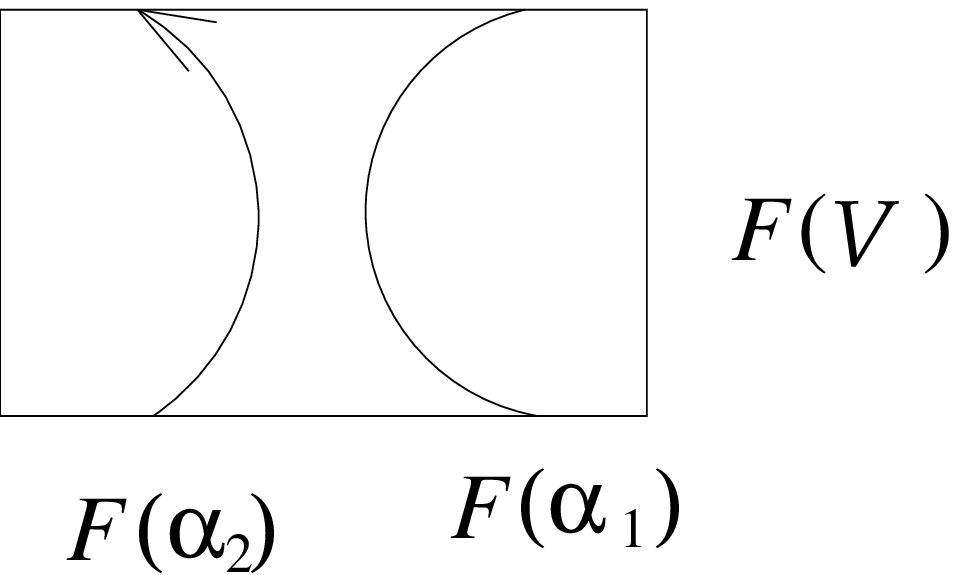}}
\caption{\label{fig:orient}} 
\end{figure} 

By the above claim, Proposition \ref{prop:smoothing}, and 
the subsequent result, 
there exists an amenable embedding 
$F': \disk \ra M$ with the following properties {\/\rm :\/}  
\begin{itemize}
\item[(i)] $F'$ is homotopic to $F$ relative to  
$\disk - \disk_2^o$.
\item[(ii)] The maximum point $y$ of $x_1\circ F'| \disk_2$ 
is unique with maximum value less than $t + \eps + \delta$. 
\item[(iii)] The number of critical points of $x_1\circ F'$ 
equals that of $x_1\circ F$. 
\item[(iv)] There exists a gradient-like vector 
field $\zeta'$ on $\disk$ extending 
$\zeta| \disk - \disk_2^o$ smoothly
so that the closure of 
the $\zeta'$-trajectory $\Gamma_3$ from $z$ to $y$ 
is a smooth arc and lies in ${\cal V} \cap \disk_2$. 
\item[(v)] $(x_2\circ F', x_3\circ F') = (x_2\circ F, x_3\circ F)$
holds on $\clo({\cal V}) \cap \disk_2$.  
\end{itemize}

To obtain (v), we may have to take a slightly larger cubical neighborhood $B$. 
In fact, since $F' = F$ on $\disk - \disk_2^o$, 
we have $(x_2\circ F', x_3\circ F') = (x_2\circ F, x_3\circ F)$ on $\cal V$. 
We let $f' = x_1 \circ F'$. Then   
$F'(x) = (f'(x), x_2\circ F(x), x_3\circ F(x))$ for 
every $x \in \cal V$ in the natural coordinates of $B$.  

Since $\Gamma_1 \cap \Gamma_2 = \{z_3\}$,   
$\Gamma_2 \cap\Gamma_3 = \{z\}$, and  
$\Gamma_1 \cap \Gamma_3 = \emp$ hold,  
the union $\Gamma = \Gamma_1 \cup \Gamma_2\cup \Gamma_3$ 
is a smoothly imbedded compact arc containing $z_3$ in the interior 
and endpoints $r$ and $y$.
Clearly, $F'(\Gamma)$ is a subset of $B$. 
There exists a positive constant $\delta_1$ so that 
the minimum $\bdd$-distance from $\delta \clo(B)$ 
and $F'(\Gamma)$ is greater than $2\delta_1$.   
Let $B'$ be a cubical neighborhood 
of $F'(\Gamma)$ in $B$ so that the $\bdd$-distance 
from each point of $\delta \clo(B')$ to $\delta \clo(B)$
equals $2\delta_1$. 
For $\delta_1$, the equation \ref{eqn:small} 
shows that there exists $\eps_5$, $\eps_5 >0$, such that  
\begin{equation}
\bdd(F_\eps(x), F(x)) \leq \delta_1 
\quad\mbox{if}\quad x \in A_{t, t+\eps_5}.
\end{equation}  
We may assume that $F' = F_\eps$ for $\eps < \eps_5$. 
Therefore, we have obtained: 

\begin{lem}\label{lem:inject}
$F'(x)$ does not belong to $B'$ if 
$x \in \disk - \cal V$. 
\end{lem}

\begin{figure}[h] 
\centerline{\epsfysize=6.5cm
\epsfbox{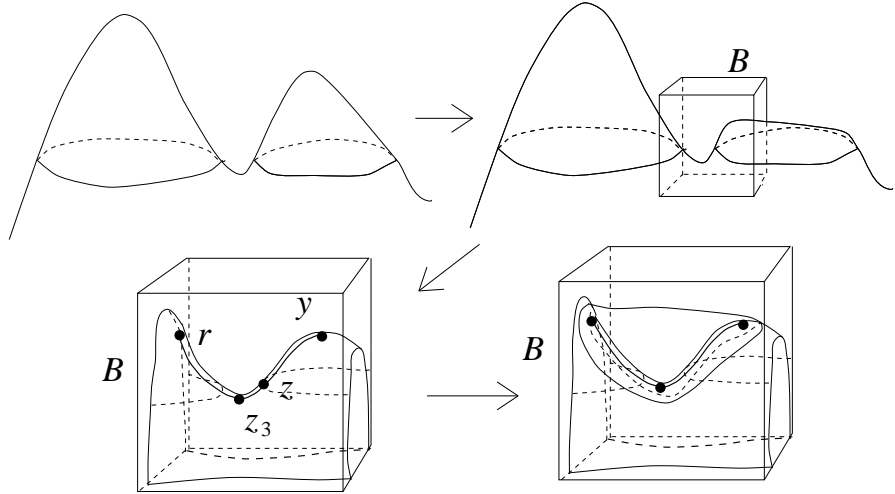}}
\caption{\label{fig:mv1} The cancellation steps in the first case} 
\end{figure} 

By a slight extension of the arguments in Section 5 of 
Milnor \cite{Mil}, 
we prove Lemma \ref{lem:coordfun} in the appendix that 
there exists an open neighborhood $U$ of $\Gamma$ in $F^{-1}(B')$ 
with a coordinate chart $u = (u^1, u^2)$ so that 
the $\zeta'$-trajectory extending $\Gamma$ is mapped 
into the $u^1$-axis, $u(r) = (-1, 0)$, 
$u(z_3) = (0, 0)$, $u(y) = (1, 0)$, and for every $x$, $x \in U$,    
\[f'(x) = f'(r) + 2\int^{u^1(x)}_{-1} v(t)dt - (u^2(x))^2\]
for a real-valued function $v$ defined on the 
real line $\bR$ that is positive for $0 < t <1$,  
and zero at $0$ and $1$ and negative elsewhere. 

The property of $v$ also can be obtained from the equation:
\[
2\int^1_0 v(t) dt = f'(y) - f'(z_3) > 0 
\quad\mbox{and}\quad 
2\int^0_{-1} v(t) dt = f'(r) - f'(z_3) < 0
\] 
since $f'$ is decreasing on the interval in $\Gamma$ 
from $r$ to $z$ and increasing on that in $\Gamma$ 
from $z$ to $y$. Thus, $v$ is negative on the interval $(-1, 0)$ and 
positive on $(0, 1)$ and negative on $(1, 1+\mu]$
for some small $\mu, \mu > 0$. 

Now we choose a negative-valued function $w$ defined on $\bR$
agreeing with $v$ in $\bR - (-1, 1 + \mu)$ satisfying the following 
conditions: 
\begin{itemize}
\item[(i)] $(x, 0)$ belongs to $u(U)$ for $x \in [-1 - 2\mu, 1 + 2\mu]$, 
\item[(ii)] $2\int^{1 + \mu}_{-1} w(t)dt = f'(y') - f'(r) < 0$ 
where $y' = u^{-1}(1 + \mu, 0)$, $f'(y') > f'(z)$, and 
\item[(iii)] $\int_{-1}^x w(t) dt \geq \int_{-1}^x v(t) dt$ 
for $x \in \bR$. 
\end{itemize}

\begin{figure}[h] 
\centerline{\epsfysize=5cm
\epsfbox{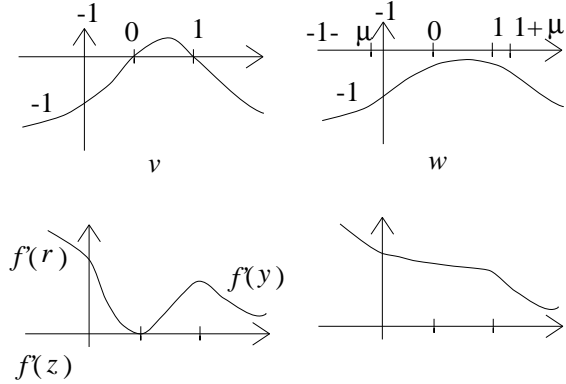}}
\caption{\label{fig:uf17} Graphs of functions $v$ and $w$ and 
their primitive functions $f'(z)$.} 
\end{figure} 

One can find a nonnegative smooth function $s: \bR^2 \ra [0, 1]$  
that equals $1$ on the segment in the $u^1$-axis with 
endpoints $(-1 - 2\mu, 0)$ and $(1 + 2\mu, 0)$, 
has a support in a compact subset of $u(U)$,  
and whenever $s(u^1, u^2) \ne 0$, 
we have $\partial s/\partial u^2(u^1, u^2) < 0$ for $u^2 > 0$, 
$\partial s/\partial u^2 (u^1, u^2) > 0$ for $u^2 < 0$.  
Define a real-valued function $g$ on $U$ by   
\begin{eqnarray}
g(x) &=& f(r) - (u^2(x))^2 \nonumber \\ 
&+&2(1-s(u^1(x), u^2(x)))\int^{u^1(x)}_{-1} v(t) dt 
+ 2s(u^1(x), u^2(x))\int^{u^1(x)}_{-1} w(t)dt. 
\nonumber 
\end{eqnarray}
Then $g(x) = f(x)$ for $x \in U - K$, where $K$ is  
the support of $s \circ u$. 
By (iii) and the condition on $\partial s/\partial u^2$,  
$\partial g/\partial u^2$ is nonzero whenever $u^2 \ne 0$.  
If $u^2 = 0$, $\partial g/\partial u^1$ is nonzero since $w$ is 
nonzero everywhere. 
Thus, $g$ has no critical points in $U$.  
Since the map given by $x \ra (x_2\circ F(x), x_3\circ F(x))$ 
for $x \in \cal V$ is 
a homeomorphism of $\cal V$ onto an open rectangle on $\bR^2$, 
the map given by $x \ra (g(x), x_2\circ F(x), x_3\circ F(x))$ 
for $x \in U$ and $x \ra (f'(x), x_2\circ F(x), x_3\circ F(x))$ for 
$x \in {\cal V} - U$ is an embedding into 
a subset of $B$ with the natural coordinates.  
Denote this map by $\hat g: {\cal V} \ra B$. Define 
a map $F'': \disk \ra M$ by $F''(x) = \hat g(x)$ if 
$x \in \cal V$ and $F''(x) = F'(x)$ if $x \in \disk - \cal V$. 
Since $\hat g(U) \subset B'$, 
Lemma \ref{lem:inject} and the fact that 
$F'| \disk - U$ and $\hat g: {\cal V} \ra B$ are injective 
imply that $F''$ is an embedding.   
Clearly, $F''$ is an amenable embedding 
homotopic to $F$, and $x_1 \circ F''$ has 
a unique critical point, which is a maximum point.   
Hence, the conclusion follows 
from Proposition \ref{prop:diskfill1}. 

(b) Let $\eps$ be a small positive number so that 
$v_3 - \eps > v_2$.    
The level set $f^{-1}(v_3 - \eps)$ is the union of two 
disjoint simple closed curves $\alpha_1$ and $\alpha_2$. 
Since $z_1$ is a local maximum point and $z_3$ a saddle point, 
Proposition \ref{prop:Euler} implies that the union of 
$\alpha_1$ and $\alpha_2$ is 
the boundary of an annulus $A$ in $\disk^o$ with $z_1, z_3 \in A^o$.  
We may assume without loss of generality that 
$\alpha_1$ is the boundary of a disk $\disk_1$
in $\disk^o$ such that $\alpha_2 \subset \disk_1^o$. 
Let $\disk_2$ be the closed disk in $\disk^o$ with boundary 
$\alpha_2$. Then $\disk_2 \subset \disk_1^o$, 
$z_2 \in \disk_2^o$. 
By Proposition \ref{prop:annfill}, one of the following item 
is true:
\begin{itemize}
\item[($\alpha$)] $F(\delta A)$ is the boundary of  
the union of two disjoint disks $D_1$ and $D_2$ in $M_{v_3 - \eps}$, 
so that $D_1 \cup D_2 \cup F(A)$ is the boundary 
of a compact subset $B_1$ of $M$ homeomorphic to a three-ball. 
(We may assume without loss of generality that 
$\delta D_1 = F(\alpha_1)$ and $\delta D_2 = F(\alpha_2)$.)
\item[($\beta$)] There exists an embedding $J: A \ra M_{v_3 - \eps}$ 
homotopic to $F|A$ relative to $\delta A$,    
so that $J(A) \cup F(A)$ equals the boundary 
of a compact subset in $M$, homeomorphic to a solid torus. 
\end{itemize} 

\begin{figure}[h] 
\centerline{\epsfysize=9.5cm
\epsfbox{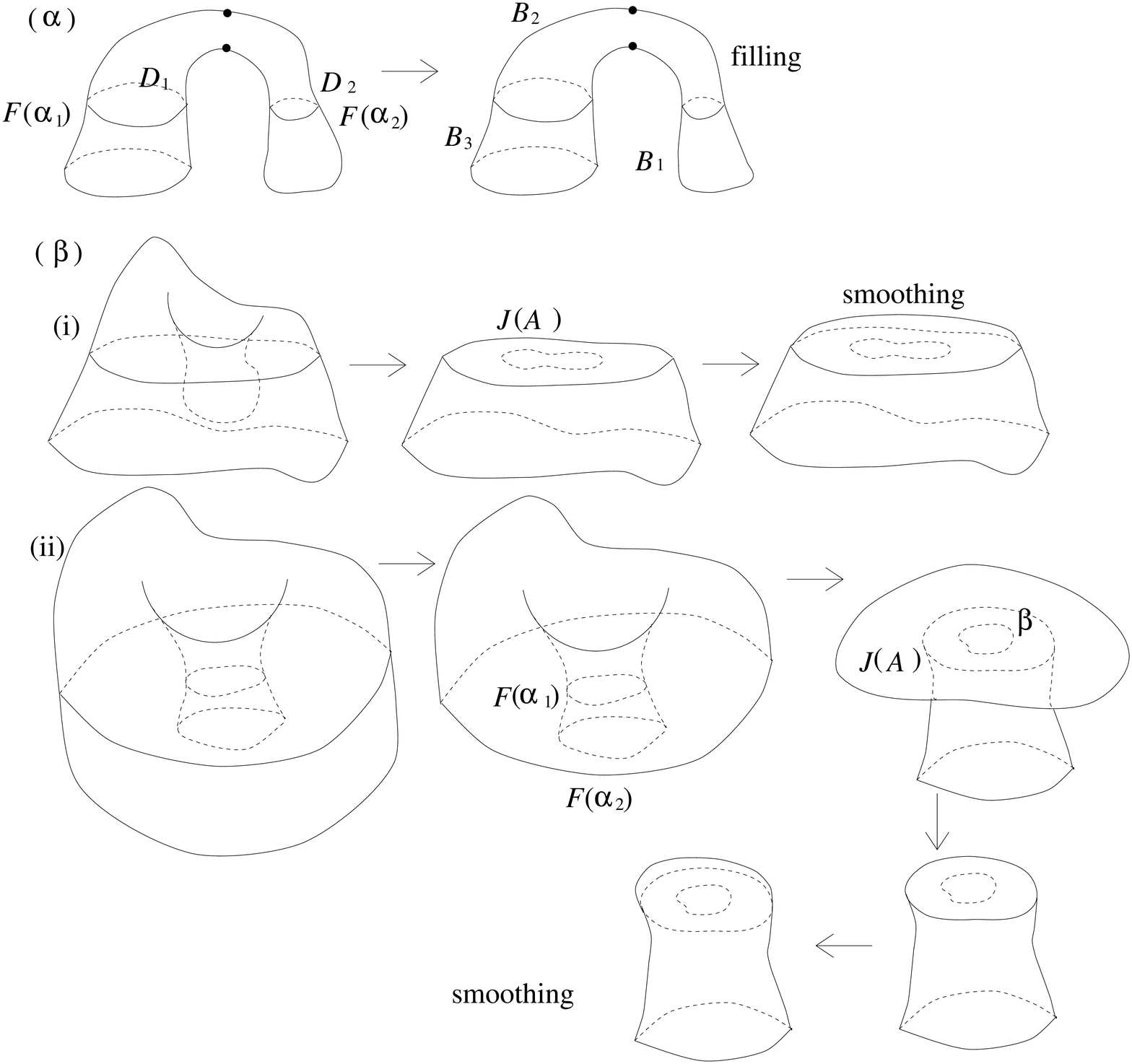}}
\caption{\label{fig:mv2} The cancellation steps in the second case
($\alpha$) and ($\beta$).} 
\end{figure} 

($\alpha$) The map $F| \disk_2$ is a negative 
amenable embedding of level $v_3 -\eps$ with one 
critical point. By reversing the $x_1$-direction 
and using Proposition \ref{prop:diskfill1}, we obtain
a disk $D'$ in $M_{v_3 - \eps}$ with boundary $F(\alpha_2)$ 
such that $D' \cup F(\disk_2)$ is a boundary of 
a three-ball $B_2$. 
By Lemma \ref{mt:open}, $D' = D_2$.  
Hence, $B_1\cap B_2 = D_2$ and $B_1 \cup B_2$ is a compact 
subset of $M$ homeomorphic to a three-ball.  
The component $A_1$ of $f^{-1}([t, v_3 -\eps])$ including 
$\delta \disk$ is homeomorphic to an annulus, 
and the boundary component of $A_1$ in the level $v_3-\eps$ 
equals $\alpha_1$, where $F(\alpha_1) = \delta D_1$.    
By Theorem \ref{thm:fill}, there exists 
a disk $D_t$ in $M_t$ so that the union of $F(A_1)$, $D_t$, and  
$D_1$ is the boundary of a three-ball $B_3$. 
Then $B_1 \cup B_2 \cup B_3$ is homeomorphic to a three-ball 
and the boundary of $B_1 \cup B_2 \cup B_3$ equals $D_t \cup F(\disk)$.  
Hence, $F$ is homotopic relative to $\delta \disk$ 
to an embedding $G:\disk \ra D_t$. 

($\beta$) As in ($\alpha$), 
$F|\disk_2$ is homotopic relative to $\delta \disk_2$ 
to an embedding $j':\disk_2 \ra M_{v_3-\eps}$.   
Clearly, $F(\alpha_1)$ either is disjoint from $j'(\disk_2)$ or 
a subset of $j'(\disk_2)^o$. 

(i) Suppose that $F(\alpha_1)$ is disjoint from 
$j'(\disk_2)$. Then the union of $F(\alpha_1)$ and $F(\alpha_2)$ 
is the boundary of the annulus $J(A)$, 
such that $J(A) \cup j'(\disk_2)$ 
is an imbedded disk in $M_{v_3-\eps}$.  
Let $H_1:\disk \ra M$ be defined by 
$H_1(x) = F(x)$ if $x \in \disk - \disk_1^o$, 
$H_1(x) = J(x)$ if $x\in A$, and 
$H_1(x) = j'(x)$ if $x \in \disk_2$. 
Then $H_1|\disk_1$ is an imbedding onto the disk 
$J(A) \cup j'(\disk_2)$. Theorem 8.1.8 
in \cite{Hirsch} shows that $H_1|\disk_1$ is 
homotopic relative to $\delta \disk_1$ 
to an embedding $H_2$ from $\disk_1$ 
onto $J(A) \cup j'(\disk_2)$ in $M_{v_3 -\eps}$.      

Applying Proposition \ref{prop:smoothing} to $F| \disk_1$ 
and $H_2$, we obtain a smooth amenable embedding 
with a unique critical point, which is a maximum point. 
Proposition \ref{prop:diskfill1} completes the proof.  

(ii) Suppose that $F(\alpha_1) \subset j'(\disk_2)^o$.  
Then by Lemma \ref{mt:open}, 
$J(A)$ is an annulus in $j'(\disk_2)$ with boundary 
$F(\alpha_1) \cup F(\alpha_2)$.  
Define a map $H_3: \disk \ra M$ by 
$H_3(x) = F(x)$ if $x \in \disk - \disk_1^o$, 
and $H_3(x) = J(x)$ if $x \in A$, and 
$H_3(x) = j'(x)$ if $x \in \disk_2$. 
Then $H_3$ is homotopic to $F$ relative to $\disk - \disk_1^o$. 
Let us choose a smooth closed curve $\beta$ in the open disk 
$j'(\disk_2) - J(A)$, so that the union of $\beta$ and $F(\alpha_1)$ 
is the boundary of an annulus, which we call $A_1$.
Then $H_3$ is homotopic relative to $\disk - \disk_1^o$
to an imbedding $H_4: \disk \ra M$ defined by 
$H_4(x) = F(x)$ if $x \in \disk - \disk_1^o$, 
$H_4|A$ is an embedding onto $A_1$  
so that $H_4(\alpha_1) = F(\alpha_1)$, 
and $H_4(\alpha_2) = \beta$, and $H_4| \disk_2$ 
is an embedding onto the disk that is the closure 
of $j'(\disk_2) - A_1 - J(A)$. Now, $H_4|\disk_1$ 
is an imbedding onto the closure of 
$j'(\disk_2) - J(A)$, a disk. Again by Theorem 8.1.8 
in \cite[p. 181]{Hirsch}, $H_4|\disk_1$ is 
homotopic relative to $\delta \disk_1$ 
to an embedding $H_5$ from $\disk_1$ 
onto the closure of $j'(\disk_2) - J(A)$ in $M_{v_3 -\eps}$.      

Similarly to (i), the proof is completed by Proposition 
\ref{prop:smoothing} and Proposition \ref{prop:diskfill1}.
\end{pf}

\section{The incompressible level surfaces}

We prove that every amenable disk with boundary in 
a level set is homotopic to an embedding into the level set. 
Our basic idea is to find a highest saddle point and cut 
the embedding just below the level of the saddle point. 
Then what is above is the union of components 
that are either disks with three 
critical points, annuli with two critical points, 
and a disk with one critical point. 
There is exactly one component that is not a disk with 
one critical point. We reduce the number of 
the critical points on it by methods we 
described in Sections 5 and 6. 

Since every embedded disk with boundary in 
a level set can be deformed to an amenable disk, 
this shows that each component of a level set is 
incompressible in $M$.  
Finally, we show that each component of a level set is 
diffeomorphic to $\bR^2$ using the incompressible 
surface theory. By a result of Palmeira \cite{Pal}, $M$ is 
diffeomorphic to $\bR^3$. This completes the proof 
of Theorem \ref{thm:main}. 

\begin{figure}[h] 
\centerline{\epsfysize=2.5cm
\epsfbox{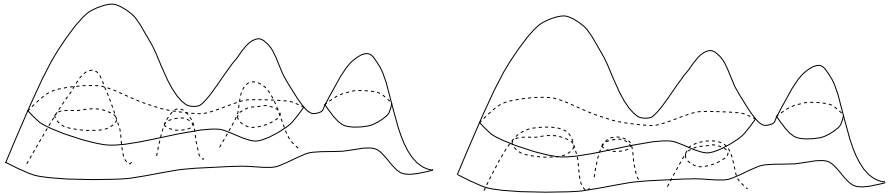}}
\caption{\label{fig:uf20} Induction Moves}
\end{figure} 

By Proposition \ref{prop:diskfill1},  
the following {\em induction hypothesis\/} is 
satisfied for $k=1$: 

\begin{hypo}
Let $\disk$ be a standard disk and 
$F: \disk \ra M$ a positive amenable embedding
of level $t$ with $k$ critical values. 
Then there exists an embedding $G$ of $\cal D$ into $M_t$ 
homotopic to $F$ relative to $\delta \disk$.  
\end{hypo}

Now assume that the induction hypothesis holds for 
$k \leq i - 1$ for an integer $i$, $i > 1$. 
We will show that the induction hypothesis holds for $k=i$. 
Let $F: \disk \ra M$ be a positive amenable embedding
of level $t$ with $i$ critical values. Let $f = x_1 \circ F$,  
$v$ the maximum of the critical values of index $1$, and $z$ the 
saddle point corresponding to $v$.   
We choose $\eps$, $\eps > 0$, so that $[v-\eps, v]$  
contains no other critical value. Letting $v' = v-\eps$, 
we see that the components of $f^{-1}(v')$ are simple closed curves  
$\alpha_1, \dots, \alpha_n$
and the components of $f^{-1}([v', \infty))$  
smooth subsurfaces $A_1, \dots, A_m$ of $\disk^o$ with boundaries 
each of whose components equals $\alpha_i$ for some $i$.  
We may assume without loss of generality that 
exactly one component $A_1$ contains 
the unique saddle point $z$ in the interior. 
Since $F| A_i$ is positive amenable, $A_i$ has at least 
one local maximum point.  
Proposition \ref{prop:Euler} show that $\chi(A_1) \geq 0$. 
Hence, $A_1$ is diffeomorphic to (i) a disk, or (ii) an annulus.
Since there is no saddle point in $A_j$ for $j > 1$, 
we have $\chi(A_j) > 0$. Thus, each $A_j$ for $j = 2, \dots, m$ 
is diffeomorphic to a compact disk
and contains a unique critical point, which is 
a local maximum point. 

(i) In this case, every $A_i$ is diffeomorphic 
to a disk. Since each $A_i$ has a unique boundary component, 
we have $m=n$ and may assume without
loss of generality that $\alpha_i$ 
is the boundary of $A_i$ in $\disk$
for each $i$, $i=1, 2, \dots, n$. 
Since $\chi(A_1) = 1$ and $A_1$ has one saddle point,
$A_1$ has three critical points, 
one of which is a saddle point and the other 
two are local extreme points.  
One of them is a local maximum point, and 
if the other one is a local minimum point, then 
the level of the highest saddle $v$ is higher than
that of the local minimum point, and since 
$[v-\eps, v)$ contains no critical value, $v-\eps$ 
is still higher than the level of the local minimum 
point. This is a contradiction since $A_1$ has no points 
below $v-\eps$. Hence $A_1$ has two local maximum points 
and a saddle point. 

Since by Propositions \ref{prop:diskfill1} and \ref{prop:diskfill2}, 
each $F(\alpha_i)$ for $i=1, \dots, n$, 
is the boundary of a compact disk in $M_{v'}$, 
the closure of precisely one component of 
$M_{v'} - F(\alpha_i)$ is homeomorphic to a disk
by Lemma \ref{mt:open}.   
Therefore, we can define an innermost closed curve
among $F(\alpha_1), \dots, F(\alpha_n)$ as follows: 
$F(\alpha_j)$ is an {\em innermost closed curve\/}
({\em with respect to $M_{v'}$}\/) if 
the disk with boundary $F(\alpha_j)$
in $M_{v'}$ includes no 
$F(\alpha_i)$ for $i \ne j$. 
We consider an innermost $F(\alpha_j)$, $j \ne 1$, included 
in the disk in $M_{v'}$ with boundary $F(\alpha_1)$.  
By Corollary \ref{cor:inner}, we obtain an amenable embedding 
$F': \disk \ra M$ homotopic to $F$ relative to $\disk - N(A_j)^o$ for 
a compact-disk neighborhood $N(A_j)$ of $A_j$ with equal number of 
singularities, and the number of components of 
$f^{\prime -1}(v')$ where $f' = x_1 \circ F'$ is 
less than that of $f^{-1}(v')$ by one. 
Hence, by an induction, 
we obtain a positive amenable embedding $F'': \disk \ra M$ 
homotopic to $F$ relative to 
\[
\disk -\bigcup_{j=2}^n N(A_j)^o
\]     
so that $F''(\alpha_1)$ is now an innermost component
and the number of critical points of $F''$ 
equals that of $F$.   
(By our choice in Corollary \ref{cor:inner}, 
we may assume that $N(A_j)$ is disjoint from $A_1$ and 
$N(A_k)$ for $k \ne j$, $k =2, \dots, n$.)  
Then by Proposition \ref{prop:diskfill2}, 
$F''|A_1$ is homotopic to an embedding 
onto a disk in $M_{v'}$ relative to $\delta A_1$. 
By Proposition \ref{prop:smoothing}, we obtain an amenable embedding 
$F''':\disk \ra M$ homotopic to $F''$ 
relative to $\disk - A_1^o$  
and $x_1 \circ F'''$ has two less number of critical points than 
$x_1 \circ F''$. Then by the induction hypothesis, 
$F$ is homotopic to an embedding $G: \disk \ra M_t$ 
relative to $\delta \disk$.  

(ii) We have $m=n-1$, and may assume without 
loss of generality that $\alpha_1$ and 
$\alpha_n$ are the boundary components of the annulus $A_1$
and that $\alpha_i$ is the boundary component of 
a compact disk $A_i$ for each $i$, $i = 2, \dots, n-1$. 
We may further assume without loss of generality that  
$\alpha_1$ is the boundary of a disk $\disk_1$ in $\disk^o$  
such that $\alpha_n \subset \disk_1^o$. 
Let $\disk_2$ be the disk with boundary 
$\alpha_n$ in $\disk$; we have $\disk_2 \subset \disk_1^o$.  
Then $F| \disk_1$ is a positive 
amenable embedding of level $v'$, and    
the map $F| \disk_2$ is a negative amenable embedding 
of level $v'$ with the number of critical points 
less than that of $F$ by at least two. 
Reversing the direction 
of $x_1$ and applying the induction hypothesis
show that $F|\disk_2$ is homotopic to 
an embedding $H_2: \disk_2 \ra M_{v'}$ relative to 
$\delta \disk_2$. Applying Proposition \ref{prop:smoothing}
to $H_2$, 
we obtain a positive amenable embedding $H_1: \disk_1 \ra M$ 
homotopic to $F| \disk_1$ relative to 
$\disk_1 - \disk_2^o$ so that $x_1\circ H_1$ has 
only three critical points, which are a local maximum point, 
a saddle point, and a local minimum point respectively.   

By Proposition \ref{prop:diskfill2}, we obtain an embedding 
$H: \disk_1 \ra M_{v'}$ homotopic to $F|\disk_1$ 
relative to $\delta \disk_1$. Let $W$ be the image 
disk $H(\disk_1)$, where $\delta W = F(\alpha_1)$.    
Let us define a map $F': \disk \ra M$ 
defined by $F'(x) = F(x)$ if $x \in \disk - \disk_1^o$ 
and $F'(x) = H(x)$ if $x \in \disk_1$. 
Now $F'$ may not be an imbedding but is 
homotopic to $F$ relative to $\disk - \disk_1^o$.  
For each $i$, $i=2, \dots, n-1$, 
$A_i$ either is a subset of $\disk_2^o$ or is disjoint from 
$\disk_1$, and $F(\alpha_i)$  
either is included in $W^o$ or is disjoint from $W$
since $F$ is an imbedding. 
Let 
\[J = \{j| F(\alpha_j) \subset W^o, 
A_j \not\subset \disk_2^o, 2 \leq j \leq n-1\}.\] 
Since $\alpha_j$, $j \in J$, is the boundary of a disk $A_j$
with one critical point, Proposition \ref{prop:diskfill1} shows that  
$F'| A_j$ is homotopic relative to $\delta A_j$ to   
an embedding of $A_j$ into $W^o$. 
Suppose that $F(\alpha_j)$, $j \in J$, 
is an innermost one among $F(\alpha_k)$ for $k \in J$ in $M_{v'}$.   
Corollary \ref{cor:inner} implies that 
$F'$ is homotopic relative to $\disk - N(A_j)^o$   
to a map $F'': \disk \ra M$ so that 
$x_1\circ F''(N(A_j)) < v'$.
By an induction eliminating such innermost 
$\alpha_j$, $j \in J$, we obtain an imbedding $L: \disk \ra M$ 
so that $L(x) \in W$ if and only if $x \in \disk_1$,   
$x_1\circ L(N(A_j)) < v'$ for $j \in J$,   
and $L$ is homotopic to $F$ relative to
\[\disk - \disk_1^o - \bigcup_{j \in J} N(A_j)^o,\] 
where $N(A_j)$ is obtained as above for each $j$,
and $L$ is smooth except on $\delta \disk_1$.  
(The neighborhoods $N(A_j)$, $j=2, \dots, n-1$, 
are assumed to be mutually disjoint.)

By Proposition \ref{prop:smoothing}, 
we obtain an amenable embedding 
$L':\disk \ra M$ homotopic to $L$ 
relative to $\disk - \disk_1^o$, and have 
a unique critical point in $\disk_1^o$, which
is a local maximum point. Thus we reduced the number of 
critical points by at least two.    
The induction hypothesis implies that there 
exists an embedding $G: \disk \ra M_t$ 
homotopic to $L'$ relative to $\delta \disk$ 
and, hence, to $F$. 
This completes the induction argument.   

We proved:
\begin{thm}\label{thm:amena} 
Let $F: \disk \ra M$ be a positive amenable embedding
of level $t$. 
Then there exists an embedding $G: \disk \ra M_t$ 
homotopic to $F$ relative to $\delta \disk$.  
\end{thm}

We will follow Chapter 6 of \cite{Hemp}. 
(Peter Scott and a number of topologists 
told me that the following lemma is true.)
Let $M$ be an orientable smooth open three-manifold. 
A {\em surface} is a properly imbedded connected 
two-manifold---possibly noncompact.   
We say that a surface $F$ is incompressible in $M$ 
if none of the following conditions is satisfied. 
\begin{itemize}
\item $F$ is a two-sphere that is the boundary of a homotopy 
three-cell in $M$, or 
\item there is a two-cell $D$ in $M$ such 
that $D\cap F = \delta D$, 
and $\delta D$ not contractible in $F$. 
\end{itemize} 

\begin{lem}\label{lem:incompressible}
Let $S$ be a two-sided surface in $M,$ 
and let $i_*:\pi_1(S) \ra \pi_1(M)$ be 
the homomorphism induced by the inclusion map. 
Then if the kernel of $i_*$ is not trivial\cO 
then there is an embedded two-cell $D$ 
in $M,$ transversal to $S$, such that
$D \cap S = \delta D,$ 
and $\delta D$ is not contractible in $S$. 
\end{lem} 

\begin{cor}\label{cor:incompressible}
If $S$ is a two-sided incompressible surface in 
$M,$ then $i_*$ is injective. 
\end{cor} 

Now, let $M$ be a simply connected 
open two-convex affine three-manifold. 
Recall that $x_1$ is a function on $M$ defined 
by $x_1 = x_1\circ \dev$ for a coordinate 
function $x_1$ on $\bR^3$. We defined $M_t$ to be 
the inverse image $x_1^{-1}(t)$ in $M$ for a real number $t$. 
Let $F$ be a component of $M_t$, which 
is a properly imbedded submanifold of $M$. 
By Lemma \ref{mt:open}, $F$ is not homeomorphic to a sphere. 
We claim that $F$ is incompressible. 
Suppose that there is an embedded two-cell $D$ in $M$, 
transversal to $F$, such that $D\cap F = \delta D$ holds. 
It is easy to see that $D$ can be perturbed to be an imbedded disk 
so that $x_1|D$ is a Morse function.
Theorem \ref{thm:amena} shows that $\delta D$ bounds a disk in $F$.  

Since $M$ is simply connected, 
$F$ is simply connected by Corollary \ref{cor:incompressible}. 
Hence $F$ is diffeomorphic to $\bR^2$. 
We conclude that $M$ is foliated by leaves that 
are diffeomorphic to $\bR^2$. The leaves are closed in $M$
since $M_t$ is closed in $M$.
By Corollary 3 of Palmeira \cite{Pal}, $M$ is 
diffeomorphic to $\bR^3$. This completes 
the proof of Theorem \ref{thm:main}.   

\section{The shrinkable dimension and $d$-convexity} 
\label{sec:dconv}

We first give our definition of shrinkable 
dimension and list Lie groups where this dimension is 
calculated. The shrinkable dimension of 
an affine group $\Gamma$ is the number $d$ such that 
for any element of the linear part of group $L(\Gamma)$
there is a direction of at least $(n-d)$-dimension 
of the standard $n$-ball which does not shrink to the origin. 
(See Figure 1.) The dimension is calculated for Lie groups  
using the representation restricted to the maximal tori.

The purpose of this section is to prove 
Theorem \ref{thm:dconv}. 
$M$ is a closed affine $n$-manifold
with a development map $\dev:\tilde M \ra \bR^n$ and 
the holonomy homomorphism $h: \pi_1(M) \ra \Aff(\bR^n)$.
We assume that $h(\pi_1(M))$ having shrinkable 
dimension less than equal to $d$. If $M$ is not 
$d$-convex, then $\che M$ includes a $(d+1)$-simplex  
such that $T \cap \ideal{M} = F_1^o \cap \ideal{M} \ne \emp$ 
for a side $F_1$ of $T$ by Proposition \ref{prop:mconv}. 

\begin{figure}[h]
\centerline{\epsfysize=9cm
\epsfbox{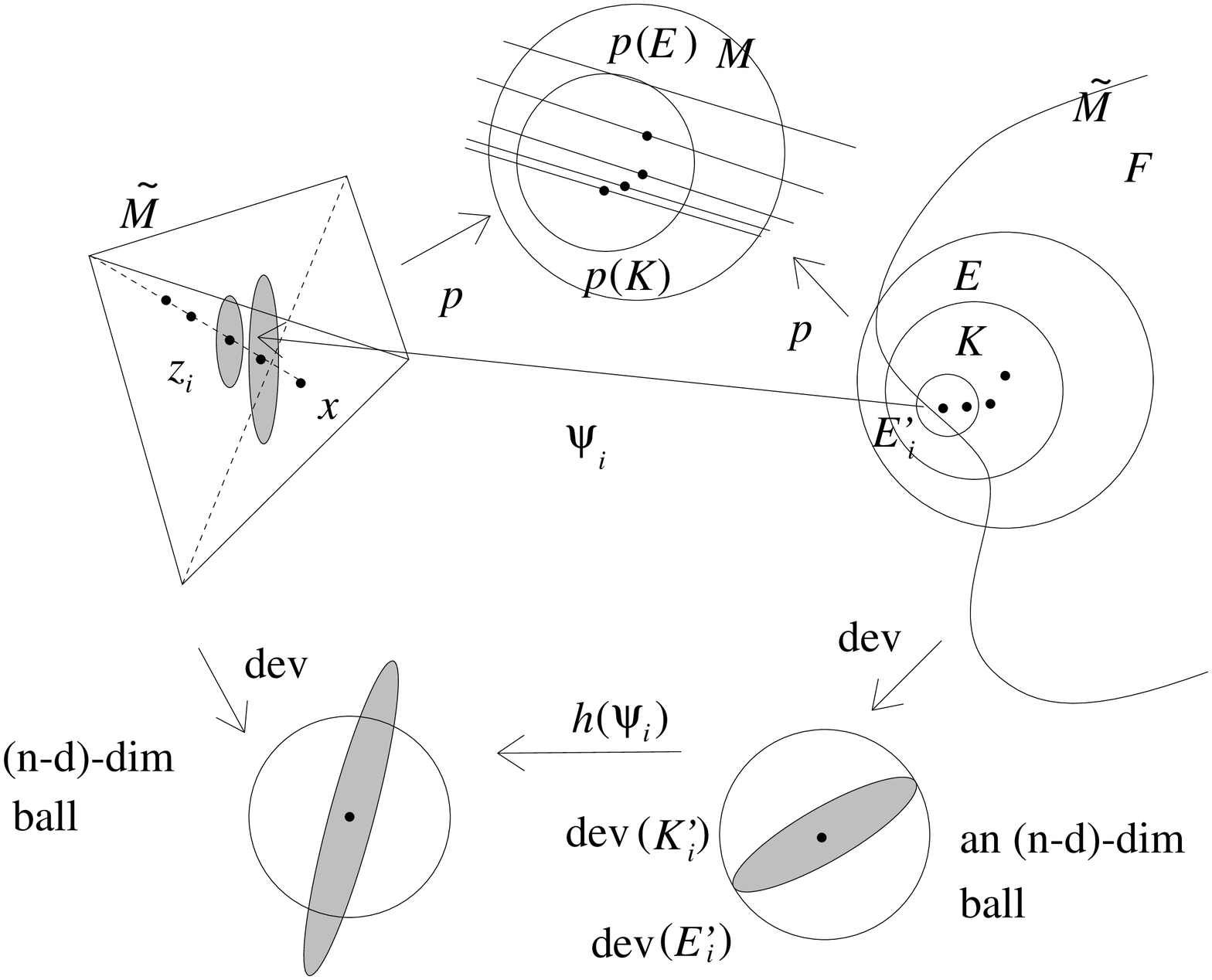}}
\caption{\label{fig:sdarg}} 
\end{figure}

We prove that if the diameter of $F_1^o \cap \ideal{M}$ 
is small, then we use the method of Carri\`ere in \cite{Car}. 
That is, we look at a ray starting from  
the vertex $v_1$ opposite $F_1$ and ending at a point 
of $F_1^o \cap \ideal{M}$. We look at a sequence of 
balls that this ray projected to $M$ passes through and using 
the fact that the shrinkable dimension is 
less than or equal to $d$, we show that there exists 
a sequence of ellipsoids on the ray that does not
degenerate in $n-d$ dimensions and hence meet a compact 
subset of $T \cap \tilde M$ infinitely often. 
Since these ellipsoids are 
subsets of a fixed compact subset of $\tilde M$, 
we obtain a contradiction to proper-discontinuity 
of the deck transformation group action of $\tilde M$.
(See Figure \ref{fig:sdarg}.) 

If the diameter of $F_1^o \cap \ideal{M}$ is not so small, 
then we change $T$ by changing the normal direction of 
$F_1$ in the confine of $\che M$. This can be done, 
and for a generic choice of the normal direction, 
we can show that the diameter can become arbitrarily 
small. Our major effort is directed toward controlling 
the size of $F_1^o \cap \tilde M_\infty$ (see Figure \ref{fig:bending}).  
This will complete the proof of Theorem \ref{thm:dconv}.

Recall the standard Euclidean metric $\bdd$ on $\bR^n$. 
We give a definition of shrinkable dimension, which 
is stronger than that of Carri\`ere \cite{Car}
but probably just as useful.  
Let $\bdd'$ be the pull-back metric $f^*\bdd$ 
for a linear automorphism $f: \bR^n \ra \bR^n$; 
let $\Gamma$ be a subgroup of $\GL(n, \bR)$ 
acting on $\bR^n$ in a standard manner.  
For $\eps > 0$, 
let $E_\eps$ be the standard $n$-$\bdd'$-ball of 
radius $\eps$ given by $\{v \in \bR^n| \bdd'(O, v) \leq \eps\}$ 
where $O$ is the origin.
For $i \geq 1$, 
an {\em $i$-$\bdd'$-ball $E^i_\eps$ of radius $\eps$\/} 
is the intersection of $E_\eps$ with 
a vector subspace of dimension $i$ of $\bR^n$.  
We say that $\Gamma$ has {\em shrinkable dimension\/} 
less than or equal to $d$ for $0 \leq d \leq n-1$,  
or say $\hbox{sd} \Gamma \leq d$, if 
given $\eps$, $\eps > 0$, 
there exists a positive constant $\eta$
depending only on $\eps$
such that for every $\vth \in \Gamma$ 
there exists an $(n-d)$-$\bdd'$-ball $E_{\eps, \vth}$ of 
radius $\eps$ such that 
\begin{equation}\label{eq:disc1}
\bdd'(\delta \vth(E_{\eps, \vth}), O) \geq \eta 
\end{equation}
where $\delta \vth(E_{\eps, \vth})$ is the boundary of 
$\vth(E_{\eps, \vth})$. (See Figure 1.) 
Since $\vth$ is linear, the condition is independent of $\eps$. 
That is, if the condition is satisfied for 
a certain $\eps$, then the condition is satisfied for 
every positive $\eps$. 
Moreover, the condition is independent of which 
Euclidean metric we use; that is, if $\bdd'$ is 
given by the pull-back metric $g^*\bdd$ 
for another linear automorphism $g: \bR^n \ra \bR^n$, 
then the above condition is true for a different choice of $\delta$.  
This follows since $\bdd'$ and $\bdd$ are 
quasi-isometric metrics on $\bR^n$. 
Thus, we need to verify equation \ref{eq:disc1} only for 
a certain metric. 


We give some examples: Let $G$ be 
a connected reductive subgroup of $\GL(n, \bR)$. 
Then $G$ admits an internal Cartan decomposition
\[G = KTK\] 
where $K \subset G$ is a maximal compact subgroup 
and $T$ a real maximal torus (see \cite[p. 249]{Helg}).
Suppose that we can conjugate simultaneously so that  
$K \subset O(n)$, and $T \subset D(n)$,
where $O(n)$ is the orthogonal subgroup and $D(n)$ 
the group of diagonal matrices.   
Thus, with respect to a fixed basis $\bold f$ of $\bR^n$ 
the matrix $M(g)$ of every element $g$ of $G$ can be written 
$M(g) = R_1(g)D(g)R_2(g)$ where $R_1(g), R_2(g) \in O(n)$ 
and $D(g)$ is a diagonal matrix with positive diagonal 
element. We may further assume that the diagonal 
entries of $D(g)$ can be written in the decreasing 
order, i.e., $a_{11} \geq a_{22} \geq \cdots \geq a_{nn}$. 
For each element $g$ of $G$, let $n(g)$ be 
the number of diagonal elements of $D(g)$ less than $1$. 
(We remark that $n(g)$ or $n(g^{-1}) < n/2$ if $a_{ii} = 1$ for 
a certain number of $i$s.) 

If $n(g) \leq d$ for every $g \in G$, then 
we claim that the shrinkable dimension of $G$ is 
less than or equal to $d$: 
Let $\bdd'$ be the metric pulled-back from $\bdd$ 
by the linear map sending the basis $\bold f$ to 
the standard basis of $\bR^n$.  
We consider the $(n-d)$-$\bdd'$-ball $E^{n-d}_{\eps}$ of 
radius $\eps$ that is the intersection of $E_\eps$ with 
the vector subspace given by $x'_{n-d+1}=0, \dots, x'_n=0$
where $x'_1, \dots, x'_n$ are coordinates associated with $\bold f$. 
For each $g \in G$, 
let $E_{\eps, g}$ be $R_2(g)^{-1}(E^{n-d}_\eps)$, 
an $(n-d)$-$\bdd'$-ball of radius $\eps$.  
Then $g(E_{\eps, g}) = R_1(g)D(g)(E^{n-d}_\eps)$. 
Since $a_{11}, \dots, a_{n-d, n-d} \geq 1$, 
it follows that 
\[\bdd'(\delta M(g)(E_{\eps, g}), O) = 
\bdd'(\delta R_1(g)D(g)(E^{n-d}_\eps), O) \geq \eps.\] 

Let $\SLt(n, \bR)$ be the group of $n\times n$-matrices 
of determinant $1$, and ${\rm Sp}(2n, \bR)$ the group of 
symplectic linear maps $\bR^{2n} \ra \bR^{2n}$ 
with respect to the standard symplectic form. 
Let $O(p, q)$ denote the subgroup of $\GL(n, \bR)$ 
of linear maps preserving the quadratic form 
given by 
\[x_1^2 + \cdots + x_p^2 - x^2_{p+1} - \cdots 
- x^2_{p+q}, p+q = n.\] 
From the dimension of the real maximal torus 
and their standard representation (see \cite{Car} 
and \cite[p. 310]{Oni}), 
we can deduce  
\begin{equation}\label{eqn:collapse}
\hbox{sd}\SLt(n, \bR) \leq n-1, \quad 
\hbox{sd} O(2, q) \leq 2 \quad (2 \leq q), \quad
\hbox{sd} {\rm Sp}(2n, \bR)\leq n.  
\end{equation} 

There is a homomorphism $L$ from the group $\Aff(\bR^n)$ of
affine transformations of $\bR^n$ to $\GL(n, \bR)$  
given by sending $g$ to its linear part $L(g)$. 
A subgroup $G$ of $\Aff(\bR^n)$ of $\bR^n$ 
has {\em shrinkable dimension 
less than or equal to $d$} if its linear part 
$L(G)$ has shrinkable dimension less than or 
equal to $d$.  

We will now give the proof of Theorem \ref{thm:dconv}.
Let $M$ be a closed affine $n$-manifold 
and $(\dev, h)$ a development pair.   
Then by assumption the linear part 
of the holonomy group $h(\pi_1(M))$ has 
shrinkable dimension 
less than or equal to $d$. 
Suppose that $M$ is not $d$-convex.
We will derive a contradiction. 

Let $p:\tilde M \ra M$ be the covering map, and 
$F$ a fundamental domain of $\tilde M$. Since $\tilde M$ 
is an open manifold, each point $x$ of $F$ is  
contained in the interior of 
a convex compact ball $D$ in $\tilde M$ 
such that $p| D$ is an imbedding onto a closed ball in $M$. 
(We can do this by choosing $D$ within a chart.)
Since $\clo(F)$ is compact, using the Lebesgue number, 
we may assume that for each point $x$ of $F$, $\tilde M$ includes 
a compact $\bdd$-ball $B$ of radius $2\eps$ with center $x$ 
so that $p| B$ is an imbedding onto a ball in $M$, 
and $\dev| B$ is an imbedding onto a $\bdd$-ball of 
radius $2\eps$ with center $\dev(x)$ (see Proposition \ref{prop:extconv}).  
We say that $B$ is a {\em tiny ball of radius $2\eps$ and 
center $x$\/}.  

Let $\delta$ be a positive number for this $\eps$ 
satisfying equation \ref{eq:disc1} for 
the linear part of $h(\pi_1(M))$. 
Since $M$ is not $d$-convex, there exists
a $(d+1)$-simplex $T$ in $\che M$ and
$F_1$ the side of $T$ so that 
\[T \cap \ideal{M} = F_1^o \cap \ideal{M} \ne \emp\] 
by Proposition \ref{prop:mconv}.
By definition of simplices, $\dev| T$ is 
an imbedding onto a $(d+1)$-simplex $\dev(T)$ in $\bR^n$. 
Let $v_1$ be the vertex of $\dev(T)$ opposite to $\dev(F_1)$.

{\em We first study the case where
the $\bdd$-diameter of $\ideal{M}^T$ in $T$ is 
less than or equal to $\delta/2$ following Carri\`ere's argument
\cite{Car}.\/} 
Let $x \in F_1 \cap \ideal{M}$ and $\alpha$  
the line starting from $v_1$ and ending at $x$. 
We have $\clo(\alpha) - \alpha = \{x\}$, 
where $\clo(\alpha)$ is a segment in $T$ with endpoints 
$x$ and $v_1$.   
Let $M$ have a complete metric denoted by $\bar d$. 
Then $\tilde M$ has an induced complete metric, 
which is also denoted by $\bar d$.  
Under the covering map $p:\tilde M \ra M$, 
$p|\alpha$ is a $\bar d$-semi-infinite 
arc in $M$. Then since $M$ is compact, 
$p|\alpha$ passes arbitrarily close to a point $z$ in $M$. 

Let $\tilde z$ be a point of $\tilde M$ 
corresponding to $z$ in the fundamental domain $F$. 
$\tilde M$ includes a tiny ball $E$ with 
$\bdd$-radius $2\eps$ and center $\tilde z$.    

Since $p|\alpha$ enters $p(K)$ for the tiny ball 
$K$ in $E$ of $\bdd$-radius $\eps/2$ with center $\tilde z$  
infinitely often, we may choose a monotone sequence of points 
$\{z_i\}, z_i \in \alpha$, converging to $x$, 
so that $p(z_i) \in p(K)$ for each $i$.
There exists a sequence $\{z'_i\}$ of points of $E$   
such that $p(z'_i) = p(z_i)$ and, hence,  
$z_i = \psi_i(z_i')$ for a deck transformation 
$\psi_i$ for each $i$. 
For each $i$, $p(E)$ includes a compact ball $E_i$ with
$p(z_i) \in E_i$ so that $E_i = p(E'_i)$ for 
a tiny ball $E'_i$ in $E$ of $\bdd$-radius $\eps$ with center $z'_i$.
(See Figure \ref{fig:sdarg} from now on.)
Since the holonomy group $h(\pi_1(M))$ has 
shrinkable dimension less than or equal to $d$, 
it follows that for each $\psi_i$, $E'_i$ includes 
a compact set $K'_i$ such that 
$\dev(K'_i)$ is a compact $(n-d)$-$\bdd$-ball of radius $\eps$ 
with center $\dev(z'_i)$ and 
\[\bdd(h(\psi_i)(\dev(z'_i)), h(\psi_i)(\delta \dev(K'_i))) 
\geq \delta\] 
for every $i$; that is, we have for each $i$
\[\bdd(\dev(z_i), \delta \dev(\psi_i(K'_i))) 
\geq \delta. \]
Since $\dev| \psi_i(K'_i)$ is a $\bdd$-isometry, 
$K'_i$ is an $(n-d)$-$\bdd$-ball satisfying 
\begin{equation}\label{eqn:deld}
\bdd(z_i, \delta \psi_i(K'_i)) \geq \delta. 
\end{equation}

The set $\dev(T)$ is a subset of a unique $(d+1)$-dimensional 
affine subspace $P^{d+1}$ in $\bR^n$. Let $S$ be the component of 
$\dev^{-1}(P^{d+1}) \cap \tilde M$ including $T^o$. 
Let $l: \bR^n \ra \bR^{n-d-1}$ be the affine function whose zero 
set equals $P^{d+1}$. Since we have $\psi_i(K'_i) \cap S =
(l\circ \dev| \psi_i(K'_i))^{-1}(O)$, the set $\psi_i(K'_i) \cap S$ 
equals a convex ball of dimension $\geq 1 = (n-d) - (n -d -1)$. 
Therefore, $\psi_i(K'_i) \cap S$ includes a segment $s_i$ 
so that $z_i \in s_i^o$ and $\bdd(z_i, \delta s_i) \geq \delta$ for 
each $i$ by equation \ref{eqn:deld}. 

Let $T_{3\delta/4}(x)$ be the subset of 
points of $T$ of $\bdd$-distance $3\delta/4$ from $x$.  
Since the $\bdd$-diameter of $F_1 \cap \ideal{M}$ 
is less than or equal to $\delta/2$, 
$T_{3\delta/4}(x)$ is a compact subset of $\che M$
not meeting $\ideal{M}$; hence, $T_{3\delta/4}(x)$ 
is a compact subset of $\tilde M$.  
Let $l_i = s_i \cap T$. Then 
$l_i$ is a segment in $T$ so that $z_i \in l_i^o$. 
An elementary geometric argument shows that
$l_i$ meets $T_{3\delta/4}(x)$ 
if $z_i$ is in the $\delta/8$-$\bdd$-neighborhood of $x$ in $T$ 
(compare with Carri\`ere \cite[Fig. 4]{Car}). 
Since $l_i \subset \psi_i(E)$ for 
the above compact subset $E$ of $\tilde M$, 
$\psi_i(E) \cap T_{3\delta/4}(x)$ is not empty for 
infinitely many $i$. Since $T_{3\delta/4}(x)$ is 
a compact subset of $\tilde M$, this contradicts 
the fact that the deck transformation group acts 
on $\tilde M$ properly discontinuously.  

\begin{figure}[h] 
\centerline{\epsfysize=5.5cm
\epsfbox{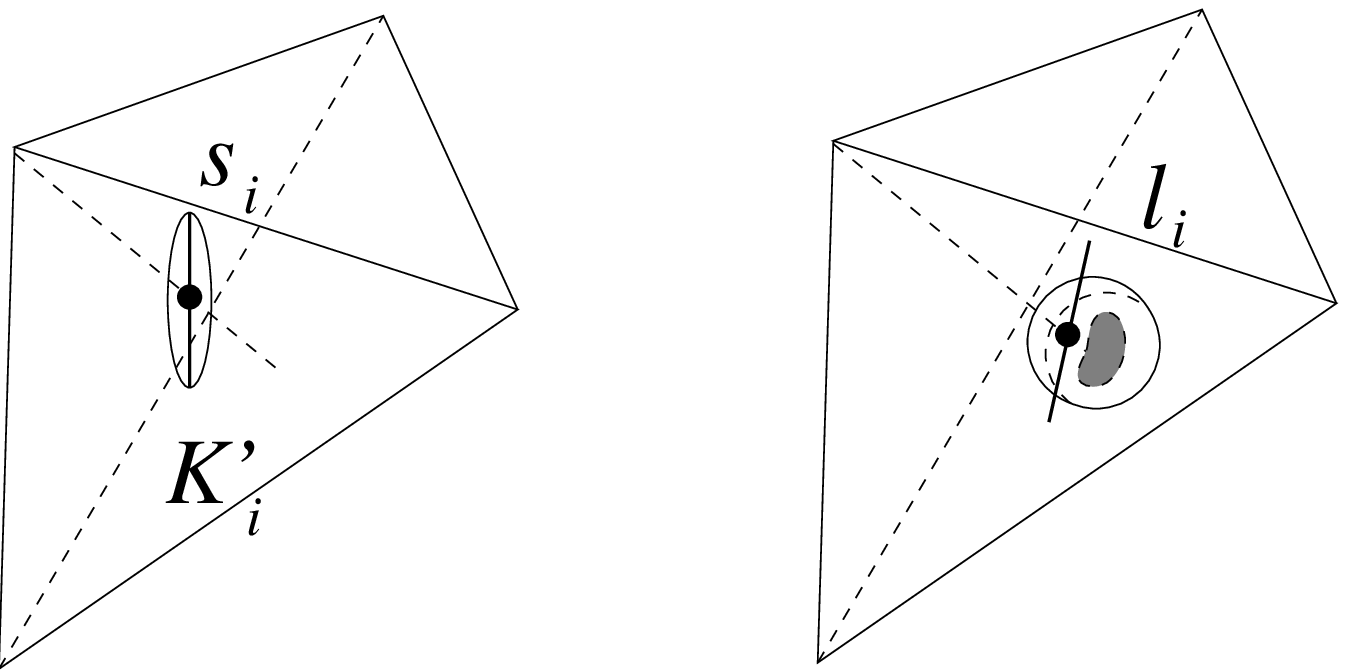}}
\caption{\label{fig:uf22}} 
\end{figure} 

{\em Suppose now that the $\bdd$-diameter of 
$F_1 \cap \ideal{M}$ is greater than $\delta/2$.\/} 
By tilting the face $F_1$ of $T$, we will now produce a new 
$(d+1)$-simplex $T'$ with \[T'\cap \ideal{M} 
\subset F_1^{\prime o} \cap \ideal{M} \ne \emp\] 
where $F_1^\prime$ is a side of $T'$
and the $\bdd$-diameter of $F_1' \cap \ideal{M}$ less than 
or equal to $\delta/4$. This will complete the proof 
of Theorem \ref{thm:dconv}.   

\begin{figure}[h] 
\centerline{\epsfysize=5cm
\epsfbox{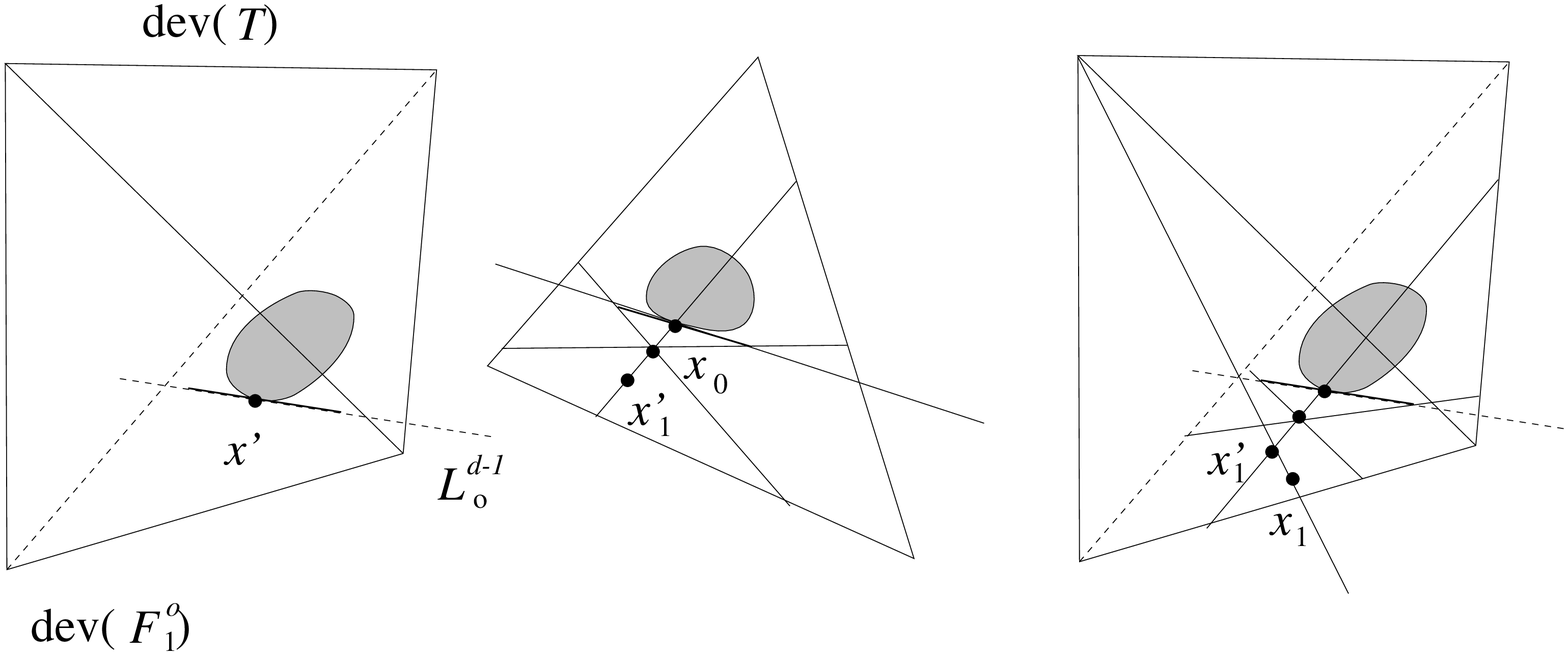}}
\caption{\label{fig:labels}} 
\end{figure}

\begin{figure}[h] 
\centerline{\epsfysize=9cm
\epsfbox{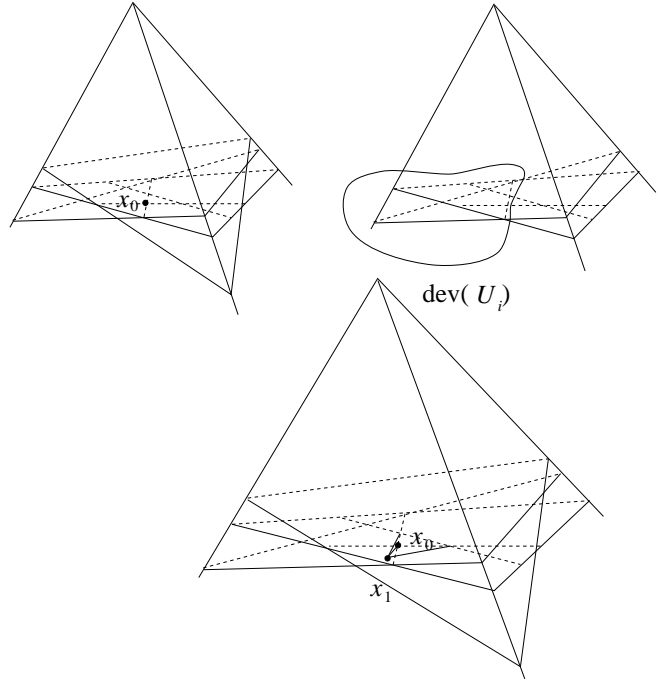}}
\caption{\label{fig:bending} Bending picture.}
\end{figure}

Let us consider the convex hull $K'$ of 
$\dev(F_1 \cap \ideal{M})$ in $\bR^n$. 
Since $\dev(F_1^o)$ is a subset of $P^d$ for a $d$-dimensional affine 
subspace $P^d$ of $P^{d+1}$, $K'$ is a subset of $P^d$.
For the purpose of choosing $T'$ and $F_1'$, 
we will now choose following objects (see Figure \ref{fig:labels}):
\begin{itemize}
\item[$L_0^{d-1}$:] Since $K'$ is convex, 
there exists an exposed point $x'$ of $K'$ with respect to $P^d$; 
that is, $x'$ is the unique intersection point of 
$K'$ with an affine hyperplane $L_0^{d-1}$ in $P^d$. 
Then we have $x' \in \dev(F_1 \cap \ideal{M})$ and  
$L_0^{d-1} \cap K' = \{x'\}$. 
\item[${\cal L}_0^{d-1}$:] We choose a $(d-1)$-simplex ${\cal L}_0^{d-1}$ 
in $L_0^d \cap \dev(F_1^o)$ containing $x'$ in the interior and 
included in the $\bdd$-ball of radius $\leq \delta/8$ with center $x'$.     
\item[$l^{d-2}_i$:] Let $l^{d-2}_i$ for $i=1, \dots, d$ 
denote the sides of ${\cal L}_0^{d-1}$. 
\end{itemize}
Let $P^d_0$ be the component of $P^d-L_0^{d-1}$ 
disjoint from $K'$. 
Choose a segment $s$ of $\bdd$-length $\leq \delta/8$ 
in the intersection of the closure of $P^d_0$ and $\dev(F_1^o)$,
transverse to $L_0^{d-1}$, and starting from $x'$. 
\begin{itemize}
\item[$L^{d-1}_i$:] We choose a point $x_0$ in $s^o$ sufficiently
close to $x'$ so that the affine hyperplane $L^{d-1}_i$ in $P^d$ including 
$x_0$ and $l^{d-2}_i$ is disjoint from $K'$ for every $i$, 
$i=1, \dots, d$.   
\item[${\cal L}^{d-1}_i$:] Let $L^{d-1}_i \cap \dev(F_1)$ 
be denoted by ${\cal L}^{d-1}_i$ for each $i$, $i=1, \dots, d$. 
\item[${\cal H}^d_i$:] For each $i$, $i=0, \dots, d$, 
let ${\cal H}^d_i$ denote the intersection of $\dev(F_1)$ 
with the closure of the component of $P^d-L^{d-1}_i$ disjoint 
from $K'$. Then ${\cal H}^d_i \cap K' = \emp$ for $i=1, \dots, d$, 
and ${\cal H}^d_0 \cap K' =\{x'\}$.
\end{itemize}

If we let $\tilde M_P = \dev^{-1}(P^{d+1}) \cap \tilde M$, then  
$\tilde M_P$ includes a compact 
$(p+1)$-ball neighborhood $U_i$ 
of the inverse image of ${\cal H}^d_i$ under $\dev| T$ 
for each $i$, $i=1, \dots, d$.
By Proposition \ref{prop:extmap} and choosing $U_i$ appropriately, 
we may assume without loss of generality that $\dev|T \cup U_i$ is an 
imbedding onto a compact subset of $P^{d+1}$.      

Let $P^2$ be a two-dimensional affine 
subspace including the segment $s$ and $v_1$.
Then $P^2 \cap \dev(F_1)$ is a segment $t$
passing through the interior of $\bigcap_{i=1}^d {\cal H}^d_i$
and $x_0$. Choose a point $x'_1$ in 
$t \cap \bigcap_{i=1}^d {\cal H}^{d, o}_i$, 
and a point $x_1$ in the exterior of $\dev(T)$
on the semi-infinite line starting 
from $v_1$ passing through $x'_1$.
We choose $x'_1$ and $x_1$ so that 
their $\bdd$-distances from $x'$ are less than
$\delta/4$. Then $x_0 = \lambda v_1 + \mu x_1 + \nu x'$ 
where $\lambda, \mu, \nu$ are the barycentric coordinates
and $\lambda + \mu + \nu = 1$. 
We choose $x_1$ sufficiently close to $x_0$ 
so that $\mu$ is close to $1$ and, hence, we obtain
\begin{equation}\label{eqn:mupos}
\mu > 0. 
\end{equation}

For each $i$, $i=1, 2, \dots, d$,
$P^{d+1}$ includes a $(d+1)$-simplex $T^{d+1}_i$  
with the following properties: 
\begin{itemize}
\item The $d$-simplex in the unit tangent sphere  
at $v_1$ determined by $T^{d+1}_i$ matches with 
that determined by $\dev(T)$. 
\item $F^d_{1, i} \cap \dev(F_1) = {\cal L}^{d-1}_i$
where $F^d_{1, i}$ is the side of $T^{d+1}_i$ 
opposite to $v_1$.
\item $x_0, x_1 \in F^d_{1, i}$.
\end{itemize}
Then it is easy to see that
$T^{d+1}_i \supset {\cal H}^d_i$, and 
$T^{d+1}_i \cap K' = \emp$. Hence, $x'$ is 
not an element of $T^{d+1}_i$. Since   
$F^d_{1, i} \cap {\cal L}^{d-1}_0 \supset l^{d-2}_i$
and we tilted $F^d_{1, i}$ away from $x'$,   
we have $T^{d+1}_i \cap ({\cal L}^{d-1}_0)^o = \emp$. 

Since $L^{d-1}_1, \dots, L^{d-1}_d$ are in general position in $P^{d+1}$, 
$F^d_{1, 1}, \dots, F^d_{1, d}$ are in general position in $P^{d+1}$. 
Hence, $\bigcap_{i=1}^d F^d_{1, i}$ is a segment 
containing $x_0$ and $x_1$.   

Recall $\tilde M_P = \dev^{-1}(P^{d+1}) \cap \tilde M$.  
We may choose $x_1$ close enough to $x'_1$ so that 
$T^{d+1}_i \subset \dev(T \cup U_i)$ for 
every $i$, $i=1, \dots, d$. 
Let $\hat T^{d+1}_i$ be the $(d+1)$-simplex 
that is the inverse image of $T^{d+1}_i$ 
under the imbedding $\dev|T \cup U_i$, where 
$\dev| \hat T^{d+1}$ is an imbedding onto $T^{d+1}_i$.
Since $T \cup U_i$ is a subset of $\tilde M$, 
so is $\hat T^{d+1}_i$. 
Letting $T^* = \bigcup_{i=1}^d T^{d+1}_i$, and 
$\hat T^* = \bigcup_{i=1}^d \hat T^{d+1}_i$,
we obtain by Proposition \ref{prop:extmap} and an induction 
that $\dev| T \cup \hat T^*$ is an imbedding onto 
$\dev(T) \cup T^*$ in $P^{d+1}$.  

For a simplex $R$ in $\bR^n$ and points $w_1, \dots, w_m$, 
let $[R, w_1, \dots, w_m]$ denote the affine simplex spanned by 
the vertices of $R$ and $w_1, \dots, w_n$ if 
it is nondegenerate.  
Since $l^{d-2}_i \subset F^d_{1, i}$ and $x_0, x_1 \in F^d_{1, i}$, 
the $d$-simplex $[l^{d-2}_i, x_0, x_1]$ is a subset of $F^d_{1, i}$. 
The $(d+1)$-simplex $[{\cal L}^{d-1}_0, x_0, x_1]$ 
is included in the $\bdd$-ball 
of radius $\delta/4$ with center $x'$ by our choice of 
${\cal L}^{d-1}_0, x_0$, and $x_1$, and its sides are 
\[
[l^{d-2}_1, x_0, x_1], \dots, [l^{d-2}_d, x_0, x_1], 
[{\cal L}^{d-1}_0, x_0], \mbox{and} \  [{\cal L}^{d-1}_0, x_1].   
\] 
The affine subspace $P^{d+1}$ includes 
a $(d+1)$-simplex $N$ with vertex $v_1$
such that the $d$-simplex in the unit tangent bundle 
at $v_1$ determined by $N$ matches with 
that determined by $\dev(T)$,  
and the side $F_1^N$ opposite to $v_1$ includes 
$[{\cal L}^{d-1}_0, x_1]$. 

\begin{figure}[h] 
\centerline{\epsfysize=3cm
\epsfbox{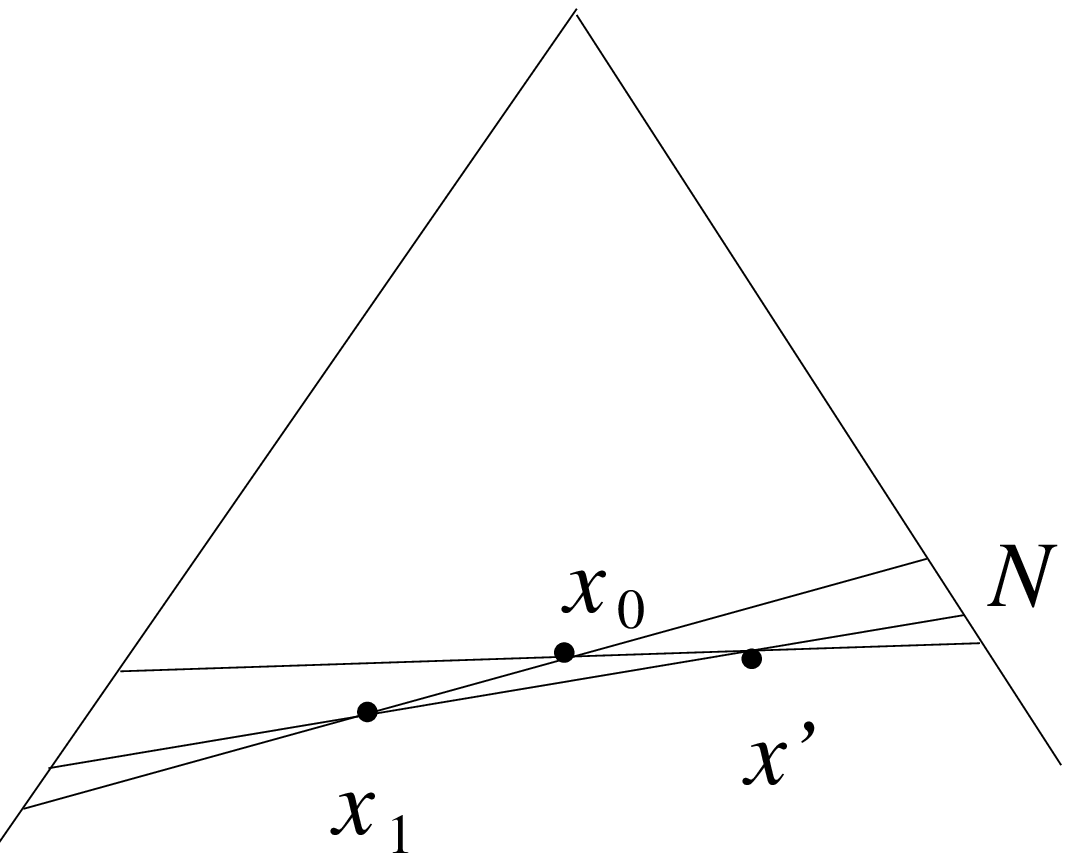}}
\caption{\label{fig:uf25}} 
\end{figure} 

\begin{lem}\label{lem:simplex}
The simplex $[{\cal L}^{d-1}_0, x_0, x_1]$ is 
included in $N$ and includes 
\[N - (\dev(T) - \dev(F_1)) 
- \bigcup_{i=1}^d (T^{d+1}_i - F^d_{1, i}).\]
\end{lem} 
\begin{pf}
We assume that the affine functions in this proof are never constant
functions. The set $\dev(T)$ can be described as 
a set of points of $P^{d+1}$ where the real-valued affine 
functions $f_1, \dots, f_{d+1}$ are simultaneously 
nonnegative. We assume without loss of generality that
the zero set of $f_i$ correspond to the side $\dev(F_i)$ of $\dev(T)$
respectively. Thus, $f_1(x_0) = f_1(x') = 0$, and $f_1(v_1) > 0$, 
although $f_1(x_1) < 0$ 
since $x_1$ does not belong to $\dev(T)$.  

Let $D$ be the set of points of $P^{d+1}$ 
where $f_2, \dots, f_{d+1}$ are simultaneously nonnegative.   
$N$ can be described as the set of points of $D$  
whose values under an affine function $g_1$ are nonnegative. 
Since the boundary of $N$ in $D$ is a zero set of $g_1$
restricted to $D$, $g_1$ is zero on 
the simplex $[{\cal L}^{d-1}_0, x_1]$.  
Thus, $g_1(x_1) = 0$ and $g_1(v_1) > 0$,    
while $g_1(x') = 0$ since ${\cal L}^{d-1}_0$ is a subset of the side
$F_1^N$ of $N$ opposite to $v$ and $x' \in {\cal L}^{d-1}_0$. 

Applying $f_1$  
to the equation $x_0 = \lambda v_1 + \mu x_1 + \nu x'$ 
where $\lambda + \mu + \nu = 1$, we 
obtain $0 = \lambda f_1(v_1) + \mu f_1(x_1)$. 
By equation \ref{eqn:mupos}, we obtain $\lambda > 0$.    
Applying $g_1$ to the above equation again, 
we obtain $g_1(x_0) > 0$, and hence $x_0 \in N^o$. 
Therefore, $[{\cal L}^{d-1}_0, x_0, x_1]$ is a subset of $N$
while $[{\cal L}^{d-1}_0, x_1]$ is a subset of the side $F_1^N$. 

Since $x_1$ does not belong to $\dev(T)$, 
and $[{\cal L}^{d-1}_0, x_0]$ is included in 
$F_1$, $f_1$ is negative in the interior of 
$[{\cal L}^{d-1}_0, x_0, x_1]$.  

Let $h_i$ for each $i$, $i =1, \dots, d$, 
be an affine function nonnegative on 
$T^{d+1}_i$ and have $F^d_{1, i}$ included in the zero set. 
Since $T^{d+1}_i \cap ({\cal L}^{d-1}_0)^o = \emp$
by our choice earlier,  
$h_i$ is negative on $({\cal L}^{d-1}_0)^o$. 
Since $x_0$ and $x_1$ are on $F^d_{1, i}$ for every $i$, 
$h_i$ is negative in the interior of 
$[{\cal L}^{d-1}_0, x_0, x_1]$ for every $i$, $i=1, \dots, d$.   

If $x$ is a point of  
\[N - (\dev(T) - \dev(F_1)) 
- \bigcup_{i=1}^d (T^{d+1}_i - F^d_{1, i}).\]
Then $x \in D$ and $g_1(x) \geq 0$, $f_1(x) \leq 0$, 
and $h_i(x) \leq 0$ for every $i$, $i=1, \dots, d$. 
Suppose that $x$ is not in $[{\cal L}^{d-1}_0, x_0, x_1]$. 
Then choose a point $y$ in the interior of $[{\cal L}^{d-1}_0, x_0, x_1]$
Since $\ovl{xy}^o$ does not meet $F_1^N$, the open line 
$\overline{xy}^o$ meets the interior of a side of 
$[{\cal L}^{d-1}_0, x_0, x_1]$ of the form 
$[l^{d-2}_1, x_0, x_1], \dots, [l^{d-2}_d, x_0, x_1]$, or   
$[{\cal L}^{d-1}_0, x_0]$. Since $f_1, h_1, \dots, h_d$ are affine 
functions, $f_1, h_1, \dots, h_d$ are simultaneously negative on 
$\ovl{xy}^o$. But one of $f_1, h_1, \dots, h_d$ is zero  
at the point of $\ovl{xy}^o$ intersecting the side. 
This is a contradiction. 
\end{pf}        

\begin{rem} In fact, 
$[{\cal L}^{d-1}_0, x_0, x_1]$ is the closure of 
\[N - (\dev(T) - \dev(F_1)) 
- \bigcup_{i=1}^d (T^{d+1}_i - F^d_{1, i})\] in $N$.
\end{rem}

Let $N_t$ for $t \in (0, 1]$ be the convex $(d+1)$-simplex 
that is the image of $N$ under the dilatation 
with center $v_1$ by the magnification factor $t$. 
Let $v'_1$ denote the vertex of $T$ corresponding to $v_1$.
Let us consider the subset $A$ of $(0, 1]$ 
whose element $t$ is such that 
$\tilde M$ includes a $(d+1)$-simplex ${\cal N}_t$
such that $\delta {\cal N}_t \ni v'_1$ 
and $\dev| {\cal N}_t$ is an imbedding onto 
$N_t$. By Remark \ref{rem:magnify}, $A$ is of form $(0, t)$ or $(0, 1]$. 

If $A$ is of form $(0, 1]$, then ${\cal N}_1 \subset \tilde M$ and
$\dev| {\cal N}_1 \cup T$ is an imbedding onto 
$N \cup \dev(T)$ by Proposition \ref{prop:extmap}. 
This means that $x' \in {\cal N}_1$, which 
is a contradiction since $x'$ is a point 
of $\ideal{M} \cap F_1^o$. 
Therefore, $A$ is of form $(0, t)$ for $0 < t \leq 1$.
By Remark \ref{rem:magnify}, 
there exists a $(d+1)$-simplex ${\cal N}_t$ in $\che M$ such that 
$\dev|{\cal N}_t$ is an imbedding onto $N_t$. 
Since $t$ does not belong to $A$, 
we have that ${\cal N}_t \cap \ideal{M} \ne \emp$. 

Again by Proposition \ref{prop:extmap} and an induction, 
$\dev|{\cal N}_t \cup T \cup \hat T^*$ is 
an imbedding onto $N_t \cup \dev(T) \cup T^*$.   
Hence, $\dev({\cal N}_t \cap \ideal{M})$ is 
disjoint from $\dev(T) - \dev(F_1)$ and 
$T^{d+1}_i - F^d_{1, i}$ for every $i$ since
they are images of subsets of $\tilde M$. 
By Lemma \ref{lem:simplex}, 
$\dev({\cal N}_t \cap \ideal{M})$ is a subset of    
$[{\cal L}^{d-1}_0, x_0, x_1]$ since $N_t \subset N$.
Since $[{\cal L}^{d-1}_0, x_0, x_1]$ is 
included in a $\bdd$-ball of radius $\delta/4$,
the $\bdd$-diameter of ${\cal N}_t \cap \ideal{M}$ is less than 
or equal to $\delta/2$.  

Let $T' = {\cal N}_t$. Then it follows that 
\[T'\cap \ideal{M} \subset F_1^{\prime o} \cap \ideal{M} \ne 
\emp\] for a side $F_1'$ and 
the $\bdd$-diameter of $F_1' \cap \ideal{M}$ is less than 
or equal to $\delta/2$. This completes the proof 
of Theorem \ref{thm:dconv}.   

\section{The singular hyperbolic manifolds} 
\label{sec:sihyp}

We prove in this section Theorem \ref{thm:sihyp}.
We will first give definition of hyperbolic manifolds 
with cone-type singular locus, 
and prove Claim \ref{claim}, which implies the Theorem \ref{thm:sihyp}.

In this section, we identify the hyperbolic space $H^3$ with 
the interior of the standard sphere, i.e., the $\bdd$-radius 
$1$ with center $O$, in $\bR^3$ in the real projective space 
$\rp^3$ and the group of hyperbolic isometries with 
the projectivized group $\PSO(1, 3)$ of Lorentz isometric 
linear transformations of $\bR^4 - \{O\}$ in $\GL(4, \bR)$. 

Let $l$ be the line in the hyperbolic space that is 
the intersection of the $z$-axis with the ball $H^3$, 
for convenience, and $N_c(l)$ is the closed hypercyclic 
neighborhood, i.e., given by points of $H^3$ of hyperbolic 
distance less than or equal to $c$ for 
a fixed $c > 0$. Let $\uco{N_c(l)-l}$ denote the 
universal cover of $N_c(l) -l$, which has a hyperbolic 
structure. The developing map $\dev_l$ of $\uco{N_c(l)-l}$ to 
$H^3$ is precisely the covering map of $N_c(l)-l$. 
Since $\uco{N_c(l)-l}$ is the universal cover of $N_c(l)-l$, 
it has coordinates $(r, \theta, z)$ where for a point 
$p$ of $\uco{N_c(l)-l}$, $r(p)$ denotes the hyperbolic distance from 
$l$ to $\dev_l(p)$, $\theta(p)$ the standard 
projected oriented angle from the standard 
lift of the $(x, z)$-half plane given by $x>0$ to $p$ in 
$\uco{N_c(l)-l}$ and $z$ the $z$-coordinate of $\dev_l(z)$.  
The ranges of these variables are given as follows:
$r \in (0, \infty), \theta \in (-\infty, \infty),$ and $z \in (-1, 1)$.
For any self-isometry of $\uco{N_c(l)-l}$, we obtain 
an isometry of $H^3$ fixing $l$. Hence, there exists a subgroup of 
isometries of $\uco{N_c(l)-l}$ isomorphic to 
the group of complex linear maps of $\uco{\bC - \{0\}}$
considering the endpoints of $l$ to be $0$ and $\infty$.
A complex linear map of $\uco{\bC - \{0\}}$ is determined 
by $(r, \theta)$, $r \in \bR^+$, $\theta \in \bR$, 
where $r$ indicate the expansion factor and $\theta$ the rotation amount. 
We will denote the linear map by $s_{r, \theta}$ 
and denote the isometry of $\uco{N_c(l)-l}$ corresponding to 
it by the same notation. 

\begin{lem}\label{lem:comNl}
The completion $\pch{N_c(l)-l}$ of $\uco{N_c(l)-l}$ 
can be identified with the union of 
$\clo(l)$ and $\uco{N_c(l)-l}$ with the metric topology 
such that a Cauchy sequence of  
points $p_i$ in $\uco{N_c(l) -l}$ converges to 
a point $(0, 0, t) \in \clo(l)$ for $t \in [-1, 1]$ if and only if 
$\dev_l(p_i)$ converges to $(0,0,t) \in \clo(l)$.
\end{lem}
\begin{pf} 
Straightforward.
\end{pf}

Our manifold
\[({\pch{N_c(l)-l}} - \{(0, 0, -1), (0, 0, 1)\})
/<s_{r, \theta'}, s_{1, \theta}>\] 
is homeomorphic to a solid torus $T$ if $r \ne 1$. The solid torus 
has a locus of singularity $l'$ corresponding to $l$ which is 
said to be the {\em core singular locus}.
We say that $T$ is a {\em solid-torus with cone-type core-singularity $l'$
of angle $\theta$ and radius $c$}.

A pair of a three-manifold $M$ and a link $L$ in $M$ 
is said to have a {\em singular hyperbolic structure}\/ 
with cone-type singular locus $L$ if $M-L$ has a hyperbolic 
structure and for each component $K$ of $L$, there exists 
a closed neighborhood $N_c(K)$ which is isometric to 
a solid-torus with cone-type singularity $\theta_K$ for $\theta_K > 0$
of radius $c$, $c > 0$.   

Let $p: \uco{M-L} \ra M-L$ denote the universal covering map. 
Since $M-L$ has a hyperbolic structure, $\uco{M-L}$ has 
a developing map $\dev: \uco{M-L} \ra H^3$ and 
a holonomy homomorphism $h: \pi_1(M-L) \ra \PSO(1, 2)$.
Since the standard ball $H^3$ is a subset of an affine patch 
$\bR^3$ in $\rp^3$, it follows that $\dev$ induces 
an affine structure on $\uco{M-L}$. (However, the affine 
structure doesn't descend to that of $M-L$.) 
Theorem \ref{thm:main} and the following claim prove that 
$\uco{M-L}$ is diffeomorphic to $\bR^3$ and prove 
Theorem \ref{thm:sihyp}. 

\begin{claim}\label{claim}
$\uco{M-L}$ is $2$-convex. 
\end{claim}

We begin the proof of the claim. 
Let $N(L)$ denote the neighborhood of the link $L$ that is 
the union of $N_{c_i}(K_i)$, $c_i > 0$, for each component 
$K_i$ of $L$ where $N_{c_i}(K_i)$'s are mutually disjoint.
Then $\uco{M-L}$ equals the union of the connected submanifold 
$\uco{M-N(L)^o}$ the inverse image of $M-N(L)^o$ under $p$ covering 
$M-N(L)^o$ and $\uco{N(L)-L}$, the inverse image under 
$p$ of $N(L)-L$. Each component of $\uco{N(L)-L}$ is 
a cover of $N_{c_i}(K_i)-K_i$, which we denote by $C_i^j$ for some $j$ 
in the countable index set $J_i$. Hence $\uco{N(L)-L}$ consists of 
infinitely many components each isometric to $\uco{N_c(l)-l}$ for
some $c$ with respect to hyperbolic metrics. Hence, there 
exists a projective map $f_i^j: C_i^j \ra \uco{N_c(l)-l}$.
Finally, the intersection of $\uco{M-N(L)^o}$ and 
$\uco{N(L)-L}$ is a union of infinitely many disjoint $2$-cells 
which cover the tori $\delta N(K_i)$. 


Suppose that $\uco{M-L}$ is not $2$-convex. Then there exists 
a $3$-simplex $T$ in $\pch{M-L}$ with $T \cap \lideal{M-L} =
F_1^o \cap \lideal{M-L} \ne \emp$ by Proposition \ref{prop:mconv}.
Hence, there exists an affine geodesic $k$ in $T$ 
starting from the vertex $v_1$ 
of $T$ opposite $F_1$ and ending at a point $y$ of 
$F_1^o \cap \lideal{M-L}$. 
$p| k$ is a geodesic with 
respect to the hyperbolic metric. Suppose that 
$p|l$ is infinitely long under the hyperbolic metric. 
Then $\dev| k$ is an infinitely long geodesic 
 
in $H^3$ and hence $\dev(y)$ is a point of the sphere at
infinity $\delta H^3$. However, note that $\dev(y)$ is 
in the affine convex hull of four points which are images under $\dev$ 
of the vertices of $T$. Since these points are in $H^3$
and $H^3$ is a strictly convex subset of $\bR^3$, 
$\dev(y)$ does not belong to $\delta H^3$.  
Therefore, $p|k$ is finitely long. 

By looking at the geometry of $\dev| k$ and $\dev(C_i^j)$ in $H^3$, 
we see that $k$ cannot leave $C_i^j$ and re-enter the same $C_i^j$.  
Since $k$ is finitely long, $k$ cannot enter infinitely 
many $C_i^j$s. Thus $p| k$ either 
stops at a point of $M-L$ or it enters $N_\eps(K_i)$ for 
a fixed $K_i$ and every $\eps$, $0 < \eps < c_i$. 
The first possibility is impossible since $x$ belongs 
to $\lideal{M-L}$. The second possibility implies that 
$k$ enters a component $C_i^j$ covering $N_{c_i}(K_i)$ and 
never leaves it. Letting $C_{i, \eps}^j$ denote the cover 
of $N_\eps(K_i)$ in $C_i^j$, $k$ enters 
$C_{i, \eps}^j$ for each $\eps > 0$ and never leaves it. 
We assume without 
loss of generality that $k$ is in $C_i^j$ by making $k$ 
smaller if necessary. The indices $i, j$ are fixed from now on.

Since the map $f_i^j: C_i^j \ra \uco{N_{c_i}(l) -l}$ is hyperbolic 
and hence projective, $\dev \circ f_i^j = \vth \circ \dev$ for 
some $\vth$ in $\PSO(3, 1)$. This map $\vth$ is quasi-isometric with 
respect to $\bdd$ on the standard ball $H^3$ in $\bR P^3$.
Since we will be looking at one $C_i^j$, we may assume without 
loss of generality that $\vth = \Idd$ by choosing another $\dev$
if necessary. $C_i^j$ has a path-metric $d'$ induced from the Riemannian 
metric from $\bR^n$ induced by $\dev$. 
The restriction to $C_i^j$ of the induced path metric $\bdd$ on 
$\uco{M -L}$ agrees with $d'$ on $C_i^j$ as the length of path 
in $\uco{M-L}$ with boundary in $C_i^j$ can be 
shortened by a path in $C_i^j$ as we can easily show. 
Hence, the closure $\clo(C_i^j)$ in $\pch{M-L}$ can be identified 
with the completion of $C_i^j$. 
Since the completion of $C_i^j$ is quasi-isometric to 
$\pch{N_{c_i}(l)-l}$, we easily see that $\clo(C_i^j)$ is quasi-isometric 
to $\pch{N_{c_i}(l) -l}$ by an extension map $f$ of $f_i^j$.
Actually $f$ is an isometry.

Note that $\pch{M-L}$ is the union of 
the closure of $\uco{M-N(L)}$ and the closure of $C_i^j$ in $\pch{M-L}$ 
for each $j$.
Let $l^j_i$ be the subset of $\clo(C_i^h)$ corresponding to 
$\clo(l)$ in $\pch{N_{c_i}(l) - l}$ by $f$ (see Lemma \ref{lem:comNl}).
Hence, $\clo(C_i^j)$ is the union of $C_i^j$ and the segment $l^j_i$ 
corresponding to the singular locus, and $x$ belongs to $l^j_i$ 
since $k$ is a subset of $C_i^j$.  

Since $f(k)$ is a subset of $\uco{N_{c_i}(l)-l}$, $\dev(f(k))$ 
is a subset of unique open totally geodesic half-plane $P$ in $H^3$ 
with boundary $l$. $P \cap N_{c_i}(l)$ lifts to a unique subspace $P'$ 
of $\uco{N_{c_i}(l) - l}$ including $k$. Then let $P'' = f^{-1}(P')$.
It follows that $\dev \circ f|P'' = \dev|P''$ is an imbedding onto 
the convex set $P \cap N_{c_i}(l)$. Hence, the closure 
$\clo(P'')$ of $P''$ in $\clo(C_i^j)$ is a tame set, 
a $2$-ball, and obviously $\clo(P'')$ includes the segment $l^j_i$ 
in its boundary. 

Since $\dev(T)$ and $\dev(\clo(P''))$ are convex 
compact subset of $\bR^3$, $\dev(T) \cap \dev(\clo(P''))$ 
is a convex set including $\dev(k)$. 
Since $k \subset T^o$ and $k \subset P''$,
$\dev(T^o) \cap \dev(P'')$ is not empty and 
$\dev(T) \cap \dev(\clo(P''))$ is its closure.
Proposition \ref{prop:extmap} implies that $\dev| T \cup \clo(P'')$ 
is an imbedding onto $\dev(T) \cup \dev(\clo(P''))$.
Since $\clo(l) = \dev(l_i^j)$, this means that 
\[T \cap l_i^j = (\dev| T \cup\clo(P''))^{-1}(\dev(T) \cap \clo(l)).\]
Since $\dev(T)$ is a subset of $H^3$, $\clo(l)$ is a segment 
with endpoints in $\delta H^3$, both sets contain $\dev(x)$, 
it follows that $l$ is tangent to $\dev(F_1^o)$ at $x$. 
Otherwise, $l$ intersects $\dev(T^o)$, which is impossible since 
$l$ is included in the ideal set of $\uco{M-L}$.
Furthermore, $\dev(T) \cap \clo(l)$ 
is not a subset of $\dev(F_1^o)$ by a geometric consideration
that $\dev(T)$ is a subset of $H^3$ and $\clo(l)$ has 
endpoint in the boundary of $H^3$. 
Hence, $T \cap l_i^j$ is not 
a subset of $F_1^o$. This contradicts our assumption. 
Hence, $\uco{M-L}$ is $2$-convex, and the proof 
of Theorem \ref{thm:sihyp} is complete by Theorem \ref{thm:main}.

\section*{Appendix}

The purpose of this section is to prove Lemma \ref{lem:coordfun}.
The essential step to verify 
this lemma is Lemma 4.7 of Milnor \cite[p. 43]{Mil}. 

\begin{lem} \label{lem:coordfun}
Let $\Omega$ be a $2$-disk.
Let $f: \Omega \ra \bR$ be a smooth function 
with a local maximum $y$ and a saddle point $z$
with a $f$-gradient-like vector field $\zeta$.  
Let $\Gamma_1$ be a flow line from $z$ to 
a point $r$ with $f(r) > f(y)$ and 
$\Gamma_2$ a flow line from $z$ to $y$.  
Then there exists a small neighborhood $U$ of 
$\Gamma_1 \cup \Gamma_2$ with coordinates $(u_1, u_2)$ 
so that $r$ has coordinates $(-1, 0),$ $z$ $(0,0),$ 
and $y$ $(1, 0)$\/{\rm ;} $f$ can be written as 
\begin{equation}\label{eqn:gamma}
f(x) = \int_{-1}^{u_1(x)} v_1(t) dt + f(r) - u_2^2;  
\end{equation}
and $(u_1, u_2):U \ra \bR^2$ mapping $\Gamma_1 \cup \Gamma_2$
to the segment $\ovl{(-1, 0)(1, 0)}$.  
\end{lem} 

\begin{figure}[h] 
\centerline{\epsfysize=3.0cm
\epsfbox{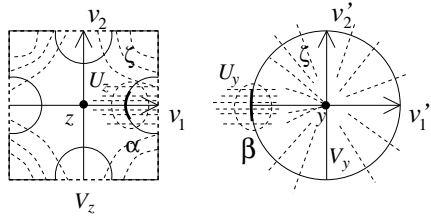}}
\caption{\label{fig:A1} A coordinatized view of 
$V_z, V'_z, \zeta, V_y, V'_y$} 
\end{figure}

\begin{figure}[h] 
\centerline{\epsfysize=3.0cm
\epsfbox{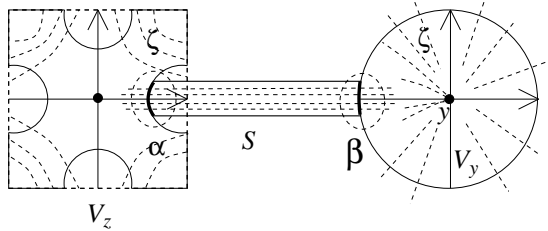}}
\caption{\label{fig:paste} $V_z, V_y$, $A$, and $S$. } 
\end{figure}

The plan to prove this lemma is we do this first for the flow 
line $\Gamma_2$ only. For neighborhoods $V_y$ and $V_z$ of 
the endpoints $y$ and $z$ of $\Gamma_2$ respectively, 
we have the canonical forms of $f$. We modify our gradient-like 
vector field $\zeta$ a little. Then we find a square $S$ obtained by 
the flow including $\Gamma_2$ and attach it to $V_y$ and $V_z$ 
(see Figure \ref{fig:A1} and \ref{fig:paste}). 
We find an integral expression of $f$ on $\Gamma_2$ and extend 
it to $S$ from $V_y$ and $V_z$. This is done using the implicit function 
theorem for infinite dimensional Banach manifolds due to Smale \cite{Smale}.   
It is shown that the integral expression on $S$ extends 
the functions of $V_y$ and $V_z$ by comparing functions 
near arcs in the boundaries of $V_y$ and $V_z$ respectively.
The extension of the integral form to a neighborhood 
of $\Gamma_1 \cup \Gamma_2$ is completely similar 
to the extension to $S$ from $V_y \cup V_z$.

We will show first that there exists 
a neighborhood $U_2$ of $\Gamma_2$ with coordinates 
$(u_1, u_2)$ so that $z$ has coordinates $(0,0)$, and 
$y$ $(1, 0)$; $f$ can be written
\begin{equation}\label{eqn:gamma2}
f(x) = \int_{0}^{u_1(x)} v_1(t) dt + f(z) - u_2(x)^2
\end{equation}
on $U_2$; and $(u_1, u_2): U_2 \ra \bR^2$ sending $\Gamma_2$ to 
$\ovl{(0, 0)(1, 0)}$.  

Choose a neighborhood $V_{z}$ for 
$z$ so that there exists a coordinate system $(v_1, v_2)$ 
where $z$ has coordinates $(0, 0)$, $v_2$ on $\Gamma_2$ 
equals $0$, $f$ is written as $f(z) + v_1^2 - v_2^2$, 
and the gradient-like vector field $\zeta$ 
takes the canonical form $(v_1, -v_2)$, 
and choose $V_y$ for $y$ so that there exists a coordinate system 
$(v'_1, v'_2)$ where $y$ has coordinates $(0, 0)$,
$v'_2$ on $\Gamma_2$ equals $0$,  
$f$ is written as $f(y) - v^{\prime 2}_1 - v^{\prime 2}_2$, 
and $\zeta$ takes the canonical form $(-v'_1, -v'_2)$. 
Actually, we choose $V_z$ and $V_y$ a little bit smaller
so that the respective neighborhoods $V'_z$ and $V'_y$ 
of their closures can be coordinatized as above.  
We assume further that $\clo(V_z)$ and $\clo(V_y)$ are homeomorphic 
to disks, that the image of $\clo(V_z)$ under the coordinate chart
$(v_1, v_2)$ equals a region given by $|v_1| \leq l$ 
and $|v_2| \leq l$ and $|v_1^2 - v_2^2| \leq l'$ for 
some positive numbers $l$ and $l'$ --- thus the boundary 
of the image of $\clo(V_z)$ consists of four arcs in hyperbolas
and eight segments --- and that $\clo(V_y)$ under 
the coordinate chart $(v_1', v_2')$ equals a disk of radius 
$l''$ with center $O$ (see Figure \ref{fig:A1}).  
The boundary of $V_z$ includes an arc segment
$\alpha$ where $f$ is constant and the boundary of $V_y$ includes 
an arc segment $\beta$ where $f$ is constant
such that $\alpha^o$ and $\beta^o$ meet $\Gamma_2$
at unique points. 
We assume without loss of generality that $\alpha$ is 
a subset of a component $\alpha'$ of a level set 
$f^{-1}(f(z) + \kappa)$, and $\beta$ is a subset of 
a component $\beta'$ of $f^{-1}(f(y) - \kappa)$ for 
a small positive number $\kappa$. 
We choose $\alpha$ and $\beta$ so that $v_2$ takes 
values $(-\kappa', \kappa')$ on $\alpha$ and $v'_2$ values 
$(-\kappa', \kappa')$ on $\beta$ for another small positive number 
$\kappa'$. We also assume that the $(v_1, v_2)$-coordinates 
of $\Gamma_2 \cap \alpha^o$ is $(\sqrt{\kappa}, 0)$ 
and the $(v'_1, v'_2)$-coordinates of $\Gamma_2 \cap \beta^o$ 
is $(-\sqrt{\kappa}, 0)$.

By altering $\zeta$ near $\alpha \cap \Gamma_2$, taking 
$\alpha$ sufficiently small, and a partition of unity argument, 
we assume that $\zeta$ is parallel to the vector field $(1, 0)$  
in a neighborhood $U_z$ of $\alpha$ in $V'_z$ and on the $v_1$-axis.  
For simplicity, we let $V'_z = U_z \cup V_z$ from now on.
We also alter $\zeta$ near $\beta \cap \Gamma_2$ using a partition 
of unity and taking $\beta$ sufficiently small so that 
$\zeta$ is parallel to the vector field $(1, 0)$ 
in a neighborhood $U_y$ of $\beta$ in $V'_y$ and on the $v'_1$-axis. 
We let $V'_y = U_y \cup V_z$ from now on.  

Morse theory tells us that there exists a homeomorphism $\psi$ 
from $\alpha'$ to $\beta'$ induced by the flow generated 
by the $f$-gradient-like vector field $\zeta$.
Let $A$ be the component of $\disk - \alpha' - \beta'$ 
not containing $y$ and $z$ homeomorphic to an open annulus. 
By Lemma 4.7 of Milnor \cite[p. 43]{Mil},  
we may alter our gradient-like vector 
field $\zeta$ within a compact subset of $A$, 
and $\psi$ accordingly, so that $\psi(\alpha) = \beta$
and $v'_2 \circ \psi = \pm v_2$ holds on $\alpha$. 
We change the sign of $v'_2$ 
if necessary so that $v'_2 \circ \psi = v_2$ holds on $\alpha$.  
We assume that $\Gamma_2$ is still 
the flow line from $z$ to $y$, which we can achieve 
by a self-diffeomorphism of $A$ supported in a compact subset of $A$
after we obtain a new flow line $\Gamma'_2$ from $z$ to $y$. 
We let $S$ denote the disk that was obtained by 
sweeping $\alpha$ to $\beta$ using the flow.
(See Figure \ref{fig:paste}.) 

Introduce a coordinate $u_1$ on $\Gamma_2$ 
so that it agrees with $v_1$ in $V_{z} \cap \Gamma_2$ 
and with $v'_1$ in $V_y \cap \Gamma_2$; 
we may have to change the coordinate function $v'_1$
in $V_y$ by a constant if necessary
so that $f(x) = f(y) - (v'_1 - c)^2 - (v'_2)^2$ 
on $V_y$ for a constant $c$, $c >0$. 
We chose $u_1$ so that $u_1(z) = 0$,
and then we see that $u_1(y) = c$.  
Restricted on the arc in $\Gamma_2$ between $z$ and $y$, $f$ is strictly 
increasing, and $f$ is decreasing elsewhere on $\Gamma_2$. 
We may write 
\begin{equation}\label{eqn:intf}
f(x) = 2\int_0^{u_1(x)} v(t) dt + f(z) 
\end{equation}
on $\Gamma_2$ where $v$ is a real-valued 
function with $v(t) > 0$ on $(0, c)$, $v(0) = v(c) = 0$, 
negative elsewhere, and equals a negative constant $C$ outside
a bounded interval of $\bR$.
Let $\eps$ be the supremum of $u_1$ on $\Gamma_2 \cap V'_z$ 
and $\eps'$ that of $c - u_1$ on $\Gamma_2 - V'_y$.
For $0 < t < \eps$, $v(t) = t$, and for 
$c - \eps' < t < c$, $v(t) = -t + c$ so that 
our equation \ref{eqn:intf} extends those of 
the local forms of $f$ in $V_z \cap \Gamma_2$ and 
$V_y\cap \Gamma_2$ respectively.

Let $\Gamma'_{k}$ for $k, |k| < \kappa$ denote the flow line from  
the point of $\alpha$ with $v_2$-coordinate $k$ to that of 
$\beta$ with $v'_2$-coordinate $k$. We introduce a coordinate function 
$u_2$ on the disk $S$ by letting $v_2(x)$ for $x \in S$ 
equal to $k$ if $\Gamma'_k$ passes through $x$. 
Then $u_2$ extends the coordinates $v_2$ on $V'_{z}$ and 
$v'_2$ on $V'_y$ smoothly by the third paragraph above. 
We let $u_2$ denote the extended function.
We extend $u_1$ on $S$ by letting $u_1(x)$ equal to 
the value of $u_1$ on the point of $\Gamma_2$ 
intersected with the level curve passing through $x$.

We claim that for $\kappa$, $\kappa > 0$, sufficiently small, 
if $|k|< \kappa$, we can write for $x \in \Gamma'_k$  
\begin{equation}\label{eqn:intf2}
f(x) = 2\int_0^{u_1(x) + \eta(x)} v(t) dt - u_2(x)^2 + f(z)
\end{equation}
for a small $C^\infty$-function $\eta$ defined on $S$. 

Let $C^0[\eps, c-\eps']$ be the Banach space of 
continuous real-valued functions defined on $[\eps, c-\eps']$ 
with the sup norm $||\quad ||_\infty$. Let $\cal O$ be the open subset 
of $C^0[\eps, c-\eps']$ of functions $g$ on $[\eps, c-\eps']$ 
such that $|g(r)|$ is less than $\min\{\eps/2, \eps'/2\}$
each $r$. 

We define a continuous functional 
${\cal F}: {\cal O} \ra C^0 [\eps, c-\eps']$ 
for a real valued function 
$g$ defined on $[\eps, c-\eps']$ by letting 
the function ${\cal F}(g): [\eps, c -\eps'] \ra \bR$ 
of variable $r$ equal to 
\[k(r) = 2\int_0^{r + g(r)} v(t) dt - f(z). \]

We claim that $\cal F$ is a $\sigma$-proper map 
(see Elworthy-Tromba \cite[Section 3]{ElTr}).
That is, the domain $\cal O$ can be written as a countable 
union of closed sets on each of which $\cal F$ is a proper map.  
$\cal O$ equals $\bigcup_{n=1}^\infty {\cal O}_n$ where 
${\cal O}_n$ equals the subset consisting of elements $g$ of 
$C^0[\eps, c-\eps']$ with 
$||g||_\infty \leq \min\{\eps/2, \eps'/2\} (1 - 1/2^n).$ 
The following calculations shows us that 
${\cal F}|{\cal O}_n$ is a proper map for each $n$:
Let $k_1$ and $k_2$ correspond to $g_1$ and $g_2$ in ${\cal O}_n$ 
under $\cal F$ respectively and observe that 
\begin{eqnarray}\label{eqn:ineq1}
|k_1(r) - k_2(r)| &= & 2| \int_0^{r + g_1(r)} v(t) dt 
- \int_0^{r + g_2(r)} v(t) dt| \nonumber \\
&= & 2| \int_{r+g_1(r)}^{r +g_2(r)} v(t) dt | \nonumber \\
&\geq & | g_1(r) - g_2(r)| \inf_{[\eps/2, c-\eps'/2]} |v| 
\end{eqnarray}
 
This implies that if $k_i \ra k$ in ${\cal F}({\cal O}_n)$, 
then $g_i$ is a Cauchy sequence for any $g_i$ in ${\cal O}_n$ 
satisfying $k_i = {\cal F}(g_i)$. Hence, 
$({\cal F}|{\cal O}_n)^{-1}(K)$ for 
a sequentially compact set $K$ is compact.

The differential of this functional $\cal F$ at $0$ is given by 
$g \mapsto vg$. Since $v$ is never zero in 
$(0, c)$, this map is a linear isomorphism of the space of
real-valued smooth functions over $[\eps, c-\eps']$. 
We define $f_k: [\eps, c-\eps'] \ra \bR$ by letting $f_k(r)$ 
equal to $f(x)$ for $x$ satisfying $u_1(x) = r$ and $u_2(x) = k$.  
Since $\cal F$ is a $\sigma$-proper-map, 
by the infinite dimensional Sard theorem of Smale \cite{Smale}, 
we can always solve the equation 
\[f_k(r) = 2\int_0^{r + \eta_k(r)} v(t) dt - k^2 + f(z)\]
for a continuous function $\eta_k$ defined on $[\eps, c-\eps']$
if $|k|$ is sufficiently small. 

Let $S_\kappa$ denote the union of $\Gamma_l$ for $|l| < \kappa$.
Assume $\kappa$ is sufficiently small.
Our equation \ref{eqn:intf2} would follow if we define $\eta$ by 
letting $\eta(x) = \eta_{u_2(x)}(u_1(x))$ for $x \in S_\kappa$ 
and show $\eta$ to be $C^\infty$. The fact that $\eta$ is smooth 
will follow if the two-variable function 
$\eta': [\eps, c-\eps'] \times [-\kappa, \kappa]$ defined by 
letting $\eta'(r, k)$ to be $\eta_k(r)$ is smooth.
Since the function $f'$ defined by letting $f'(r, k)$ equal to 
$f_k(r)$ is smooth on $S_\kappa$, $k \mapsto f_k$ is continuous. 
Hence, the inequality \ref{eqn:ineq1} shows
that $k \mapsto \eta_k$ is also a continuous 
function, and the continuity of $\eta'$ follows. 
The smoothness of $\eta'$ follows by partial differentiating 
with respect to variables $r$ and $k$ the function  
\[f'(r, k) = 2\int_0^{r + \eta'(r, k)} v(t) dt - k^2 + f(z)\]
 
and a bootstrap argument.

For $\kappa$ sufficiently small, we let $u'_1 = u_1 + \eta$ for each $x$ 
in $S_\kappa$. Then $u'_1$ is a coordinate function on $S_\kappa$.  
Near $U_z$, $f$ is of form $f(z) + v_1^2 - u_2^2$
and of form $f(x) = 2\int_0^{u'_1(x)} v(t)dt - u_2(x)^2 + f(z)$ 
in $U_z \cap S_\kappa$ where $v(t) = t$ for $t \in (0, \eps'')$
for some $\eps'' > \eps$.
Since the value of $u'_1(x)$ is near $\eps$ on $\alpha$ 
if $\kappa$ is sufficiently small, an integration shows that 
$f(x) = (u'_1(x))^2 - u_2(x)^2 + f(z)$ for 
$x$ in $S_\kappa$ sufficiently near $\alpha$ since $\eps < \eps''$. 
Hence, for $x$ in $S_\kappa \cap V'_z$ sufficiently near $\alpha$, 
we have $(u'_1(x))^2 = (v_1(x))^2$.
Since $u'_1$ and $v_1$ are positive sufficiently near $\alpha$, 
$u'_1$ extends $v_1| V_z$ to $S_\kappa$ smoothly. 

Near $\beta$, $f$ is of form $f(y) - (v'_1 -c)^2 - (u_2)^2$  
and of form $f(x) = 2\int_0^{u'_1(x)} v(t)dt - u_2(x)^2 + f(z)$
in $S_\kappa$. Since $v(t) = -t + c$ for $t \in (c - \eps''', c)$
for some $\eps'''$, $\eps''' > \eps'$, 
and $f(y) - f(z) = 2 \int^c_0 v(t) dt$, it follows that 
\begin{equation}\label{eqn:intf3}
f(y) - f(z) = 2\int_0^{c-\eps'''} v(t) dt + (\eps''')^2. 
\end{equation}
If a point $x$ in $S_\kappa$ is sufficiently near $\beta$,
then $u'_1(x) > c- \eps'''$. For such $x$, 
since $v(t) = c - t$ for $c- \eps''' < t < c$, we can write
\begin{eqnarray}\label{eqn:intf4}
f(x) &=& 2\int^{u'_1(x)}_{c- \eps'''} v(t) dt 
+ 2 \int^{c- \eps'''}_0 v(t) dt - (u_2(x))^2 + f(z) \nonumber \\ 
& = & - (u'_1(x) -c)^2 + (\eps''')^2 + 2\int^{c-\eps'''}_0 v(t) dt 
- (u_2(x))^2 + f(z). 
\end{eqnarray}
In $V'_y$, $f$ is of form 
\[f(x) = -(v_1'(x) - c)^2 - (v'_2(x))^2 + f(z) + f(y) - f(z).\]
By equations \ref{eqn:intf3} and \ref{eqn:intf4} and this equation, 
we obtain $(v_1'(x) -c)^2 = (u'_1(x) - c)^2$. 
Since $v_1'(x) < c$ near $\beta$ and so is $u'_1(x)$, 
we obtain $v_1'(x) = u'_1(x)$ near $\beta$ in $S_\kappa \cap V'_y$
assuming that $\kappa$ is chosen sufficiently small. 
Thus $u'_1| S_\kappa$ extends $v'_1| V_y$ smoothly. 
Let us denote by $u'_1$ the coordinate function extending 
$v_1$ and $v'_1$ to $V_z \cup V_y \cup S_\kappa$. 
Then $(u'_1, u_2)$ forms a global coordinate system of 
$V_z \cup V_y \cup S_\kappa$.  
Since $f$ takes canonical form in $V_z$ and $V_y$,  
it is easy to see that $f$ can be written of the form 
\[ f(x) = \int_{0}^{u'_1(x)} v(t) dt + f(z) - (v_2)^2\]
in $V_z \cup V_y \cup S_\kappa$.
We let $U = V_z \cup V_y \cup S_\kappa$ 
and $u_1 = u'_1/c$ and we let our new  $v$ be 
$c$ multiplied by $v$. Then $z$ has coordinates $(0, 0)$ 
and $y$ $(1, 0)$, $\Gamma_2$ is on the $u_1$-axis, 
and $f$ has the form of equation \ref{eqn:gamma2}.

To prove Lemma \ref{lem:coordfun} fully, we need to extend our 
integral expression \ref{eqn:gamma2} to a neighborhood of $\Gamma_1$ to 
obtain the expression \ref{eqn:gamma}. However, this 
can be done by the identical method by patching coordinates together.


\begin{thebibliography}{}

\bibitem{Abr} R. Abraham and J. Marsden, {\em Foundations of 
Mechanics}, Addison-Wesley 1987.

\bibitem{Bens} Y. Benoist, 
{\em Nilvari\'et\'es projectives}, 
Comment. Math. Helv. {\bf 69} (1994), 447--473.

\bibitem{Ben} J. P. Benz\'ecri, {\em Vari\'et\'es localement affines}
in Sem. Topologie et G\'eom. Diff., Ch. Ehresmann (1958--60), 
No. 7 (mai 1959). 

\bibitem{Ben2} J. P. Benz\'ecri, 
{\em Vari\'et\'es localement affines et projectives}, 
Bull. Soc. Math. France {\bf 88} (1960), 229--332.

\bibitem{B} M. Berger, {\em Geometry} I, 
Springer, 1987.  

\bibitem{Car} Y. Carri\`ere, {\em Autour de la 
conjecture de L. Markus sur les vari\'et\'es affines},
Invent. Math. {\bf 95} (1989), 615--628. 

\bibitem{Car2} Y. Carri\`ere, oral communication, 1989.

\bibitem{Car3} Y. Carri\`ere, {\em Questions ouvertes sur 
les vari\'et\'es affines} in S\'eminaire Gaston Darboux 
de G\'eom\'etrie et Topologie Diff\'erentielle (1991--1992)
Univ. Montpellier 1993, pp. 69--72. 

\bibitem{Casson} A. Casson and S. Bleiler,
{\em Automorphisms of surfaces after Nielsen and 
Thurston}, Cambridge University Press 1988.

\bibitem{rdaf} Y. Chae, S. Choi, and C. Park, {\em Real projective manifolds 
developing into an affine space}, International J. Math. {\bf 4} (1993),
179--191.

\bibitem{cdcr1} S. Choi, 
{\em Convex decomposition of real projective 
surfaces I}\/: {\em $\pi$-annuli and convexity}, 
J. Differential Geom. {\bf 40} (1994), 165--208.

\bibitem{cdcr2} S. Choi, {\em Convex decompositions of 
real projective surfaces II}\/: {\em Admissible decompositions}, 
J. Differential Geom. {\bf 40} (1994), 239--283.

\bibitem{iconv} S. Choi, {\em $i$-convexity of manifolds with 
real projective structures}, Proc. Amer. Math. Soc. {\bf 122}
(1994), 545--548.

\bibitem{cdcr3} S. Choi, {\em Convex decompositions of 
real projective surfaces III}\/: {\em For closed or nonorientable 
surfaces}, J. Korean Math. Soc. {\bf 33} (1996), 1139--1171.

\bibitem{uaf1} S. Choi, {\em The universal covers of affine $3$-manifolds 
of discompactness two} in Proceedings of 
the First Brazil-U.S. conference on Geometry, Topology, and 
Physics (Campinas, Brazil 1996), ed. B. Apanasov et al, 
de Gruyter 1997, pp. 107--118.

\bibitem{psconv} S. Choi, {\em Convex and concave decomposition of manifolds 
with real projective structures}, 
to appear in Mem. Soc. Math. France, math.GT/984110 at xxx.lanl.gov.


\bibitem{CG} S. Choi and W. Goldman, {\em The classification 
of real projective structures on compact surfaces},  
Bull. Amer. Math. Soc. {\bf 34} (1997), 161--171.

\bibitem{ElTr} K. Elworthy and A. Tromba, {\em Differential structures 
and Fredholm maps on Banach manifolds} in Global Analysis 
(Proc. Symp. Pure. Math. Vol 15), ed. S.S. Chern
and S. Smale, Amer. Math. Soc., 1970, pp. 45--94.

\bibitem{FH} W. Flyod and A. Hatcher, {\em Incompressible surfaces 
in punctured torus bundles}, Topology App. {\bf 13} (1982), 
261--282.

\bibitem{F} D. Fried, {\em Distality, completeness and 
affine structures}, J. Differential Geom. {\bf 24} (1986),
265--273.

\bibitem{FG} D. Fried and W. Goldman, {\em Three-dimensional 
affine crystallographic groups}, Adv. Math. {\bf 47} (1983), 
1--49.

\bibitem{FGH} D. Fried, W. Goldman, and M. Hirsch,  
{\em Affine manifolds with nilpotent holonomy}, Comment. 
Math. Helv. {\bf 56} (1981), 487--523 

\bibitem{GH} W. Goldman and M. Hirsch, {\em The radiance 
obstruction and parallel forms on affine manifolds},
Trans. Amer. Math. Soc. {\bf 286} (1984), 629--649.

\bibitem{G1} W. Goldman, {\em Projective structures with 
Fuchsian holonomy}, J. Differential Geom. {\bf 25} (1987), 297--326.

\bibitem{G2} W. Goldman, {\em Convex real projective structures 
on compact surfaces}, J. Differential Geom. {\bf 31} (1990), 
791--845.

\bibitem{Gromov} M. Gromov, {\em Stable mappings of foliations 
into manifolds}, Math. USSR-Izvestija {\bf 3} (1969), 671--694.

\bibitem{HT} A. Hatcher and W. Thurston, {\em Incompressible 
surfaces in $2$-bridge knot complements}, Invent. Math. 
{\bf 79} (1985), 225--246.

\bibitem{Helg} S. Helgason, {\em Differential 
geometry, Lie groups, and symmetric spaces}
(Pure and applied mathematics Vol. 80),
Academic Press, 1978. 

\bibitem{Hemp} J. Hempel, {\em $3$-manifolds} (Annals of 
Mathematics Studies No. 86), Princeton University Press, 1976  

\bibitem{Hirsch} M. Hirsch, {\em Differential Topology} 
(Graduate Texts in Math. 33), Springer, 1976. 

\bibitem{Hodgson} C. Hodgson, Degeneration and regeneration of 
geometric structures on three-manifolds, Ph.D. Thesis, 
Princeton University, 1986. 

\bibitem{KT} Y. Kamishima and S. Tan, 
{\em Deformation spaces on geometric structures}, in
Aspects of Low Dimensional Manifolds 
(Advanced Studies in Pure Math. Vol. 20), 
ed. Y. Matsumoto and S. Morita, Kinokuniya, 1992, pp. 263--299.  

\bibitem{Koj1} S. Kojima, {\em Nonsingular parts of 
hyperbolic $3$-cone-manifolds}, in 
Proceedings of the 37th Taniguchi Symposium on 
topology and Teichm\"uller spaces (Finland, July 1995), 
ed. S. Kojima et al, World Scientific Pub. Co. 1996, pp. 115--122.

\bibitem{Koj2} S. Kojima, {\em Deformations of hyperbolic 
$3$-cone-manifolds}, J. Differential Geometry {\bf 49} (1998), 
469--516.

\bibitem{Kuiper} N. Kuiper, {\em On compact conformally euclidean 
spaces of dimension $> 2$}, Ann. Math. {\bf 52} (1950), 487--490.

\bibitem{Mark} L. Markus, 
{\em Cosmological models in differential geometry},
mimeographed notes, University of Minnesota, 1962.  

\bibitem{Mol} E. Moln\'ar, {\em The projective interpretation of
the eight $3$-dimensional homogeneous geometries},
Beitr\"age zur Algebra und Geometrie
(Contribution to Algebra and Geometry), {\bf 38} 
(1997), 262--288

\bibitem{Mil3} J. Milnor, {\em On the existence of 
a connection with curvature zero}, Commentarii Mathematici 
Helvetici {\bf 32} (1958), 215--223.

\bibitem{Mil} J. Milnor, {\em Lectures on the h-cobordism 
theorem} (Princeton Mathematical Series), 
Princeton University Press, 1965. 

\bibitem{Mil2} J. Milnor, {\em Morse theory} 
(Ann. of Math. Studies No. 51),
Princeton University Press, 1969.

\bibitem{NY} T. Nagano and K. Yagi, {\em The affine structures 
on the real two-torus I}, Osaka J. Math. {\bf 11} (1977), 178--187.

\bibitem{Oni} A. Onishchick and E. Vinberg, 
{\em Lie groups and algebraic groups}, 
Springer, 1990.  

\bibitem{Pal} C. Palmeira, {\em Open manifolds foliated by 
planes}, Annals of Mathematics {\bf 107} (1978), 109--131.

\bibitem{Ratc} J. Ratcliff, {\em Foundations of hyperbolic manifolds}
(Graduate Texts in Mathematics, vol 149), Springer, 1994.

\bibitem{Rudin} W. Rudin, {\em Principles of Mathematical 
Analysis}, MacGraw-Hill, 1964. 

\bibitem{Sim} J. Smillie, {\em Affinely flat manifolds}, 
Doctoral Dissertation, University of Chicago, 1977.

\bibitem{Sim2} J. Smillie, {\em An obstruction to the existence 
of affine structures}, Invent. Math. {\bf 64} (1981), 411--415.

\bibitem{Smale} S. Smale, {\em An infinite dimensional version of 
Sard's theorem}, Amer. J. Math. {\bf 87} (1965), 861--866.

\bibitem{ST} D. Sullivan and W. Thurston, {\em Manifolds with 
canonical coordinate charts}\/: {\em some examples},
L'Enseignement Math. {\bf 29} (1983), 15--25.

\bibitem{Thiel} B. Thiel, Einheitliche beschreibung der 
acht Thurstonschen geometrien, Diplomarbeit, Universit\"at 
zu G\"ottingen, 1997.

\bibitem{Yosh} K. Yoshida, {\em Functional Analysis}\/ 
(Grundlehren Math. Wiss. Bd. 123), Springer, 1980. 

\end{thebibliography}
\end{document}